\documentstyle[12pt]{article}

\makeatletter
\newcommand{\l@a}[2]{\hbox to\textwidth{#1\leaders\hbox to
0.78em{\hfil .\hfil}\hfil #2}}

\renewcommand{\l@section}{\@dottedtocline{1}{0em}{2em}}
\makeatother

\addcontentsline{toc}{section}{Introduction}
\newcounter{lm}[section]
\newcounter{thm}[section]
\newcounter{prop}[section]
\newcounter{rem}[section]
\newcounter{cor}[section]

\newcounter{ex}[section]

\newcounter{defin}[section]
\begin{document}
\newcommand{\thm}{\refstepcounter{thm} {\bf Theorem \arabic{section}.%
\arabic{thm}.} }

\renewcommand{\thethm}{\thesection.\arabic{thm}}

\newcommand{\nthm} [1] {\refstepcounter{thm} {\bf Theorem 
\arabic{section}.%
\arabic{thm}:} #1}

\newcommand{\nnthm} [1] {\refstepcounter{thm} {\bf Theorem 
\arabic{section}.%
\arabic{thm}} #1}

\newcommand{\nprop} [1] {\refstepcounter{prop} {\bf Proposition 
\arabic{section}.%
\arabic{prop}:} #1}

\newcommand{\prop} {\refstepcounter{prop}{\bf Proposition \arabic{section}.%
\arabic{prop}.} }

\renewcommand{\theprop}{\thesection.\arabic{prop}}

\newcommand{\cor} {\refstepcounter{cor}{\bf Corollary \arabic{section}.%
\arabic{cor}.} }

\newcommand{\ncor} [1] {\refstepcounter{cor} {\bf Corollary
\arabic{section}.%
\arabic{cor}:} #1} 

\renewcommand{\thecor}{\thesection.\arabic{cor}}

\newcommand{\lemma}{\refstepcounter{lm}{\bf Lemma \arabic{section}.%
\arabic{lm}.} }

\newcommand{\nlemma} [1] {\refstepcounter{lm} {\bf Lemma
\arabic{section}.%
\arabic{lm}:} #1}

\renewcommand{\thelm}{\thesection.\arabic{lm}}

\newcommand{\ex}{\refstepcounter{ex}{\bf Example \arabic{section}.%
\arabic{ex}.} }

\newcommand{\nex} [1] {\refstepcounter{ex} {\bf Example
\arabic{section}.%
\arabic{ex}:} #1} 

\renewcommand{\theex}{\thesection.\arabic{ex}}

\newcommand{\defin}{\refstepcounter{defin}{\bf Definition \arabic{section}.%
\arabic{defin}.} }

\renewcommand{\thedefin}{\thesection.\arabic{defin}}

\newcommand{\rem}{\refstepcounter{rem}{\bf Remark \arabic{section}.%
\arabic{rem}.} }

\newcommand{\nrem}{\refstepcounter{rem}{\bf \arabic{section}.%
\arabic{rem}.} }

\renewcommand{\therem}{\thesection.\arabic{rem}}

\newcommand{\rems}{{\bf Remarks }}

\newcommand{\be}{\begin{equation}}
\newcommand{\ee}{\end{equation}}

\def\qed{\hfill $\Box$}
\newcommand{\vi}{{\varphi}}
\newcommand{\var}{{\varphi}}
\newcommand{\si}{{\sigma}}
\newcommand{\om} {{\omega}}
\newcommand{\G}{{\Gamma}}
\newcommand{\gi}{{G_{\infty }}}
\newcommand{\is}{{\simeq}}
\newcommand{\pe}{{\pi_1}}
\newcommand{\supp}{{\rm supp ~}}
\newcommand{\const}{{\rm const}}

\newcommand{\C}{{\bf C}}
\newcommand{\F}{{\bf F}}
\newcommand{\N}{{\bf N}}
\newcommand{\pr}{{\bf {P}}}
\newcommand{\Q}{{\bf Q}}
\newcommand{\sph}{{\bf S}}
\newcommand{\T}{{\bf T}}
\newcommand{\R}{{\bf R}}
\newcommand{\Z}{{\bf Z}}
\newcommand{\id}{{\bf 1}}
\newcommand{\z}{{\bf 0}}

\newcommand{\talpha}{{\tilde {\alpha}}}
\newcommand{\tvar}{{\tilde {\varphi}}}
\newcommand{\tdelta}{{\tilde {\delta}}}
\newcommand{\tsi}{{\tilde {\sigma}}}
\newcommand{\tpsi}{{\tilde {\psi}}}
\newcommand{\tchi}{{\tilde {\chi}}}
\newcommand{\tpartial}{{\tilde {\partial}}}
\newcommand{\tf}{{\tilde {f}}}
\newcommand{\tg}{{\tilde {g}}}
\newcommand{\tp}{{\tilde {p}}}
\newcommand{\tx}{{\tilde {x}}}
\newcommand{\tu}{{\tilde {u}}}
\newcommand{\te}{{\tilde {E}}}
\newcommand{\tA}{{\tilde {A}}}
\newcommand{\tC}{{\tilde {C}}}
\newcommand{\tH}{{\tilde {H}}}
\newcommand{\tI}{{\tilde {I}}}

\newcommand{\bk}{{\bar {k}}}
\newcommand{\br}{{\bar {r}}}
\newcommand{\bu}{{\bar {u}}}
\newcommand{\bv}{{\bar {v}}}
\newcommand{\bx}{{\bar {x}}}
\newcommand{\by}{{\bar {y}}}
\newcommand{\bz}{{\bar {z}}}
\newcommand{\bE}{{\bar {E}}}
\newcommand{\bF}{{\bar {F}}}
\newcommand{\bH}{{\bar {H}}}
\newcommand{\bR}{{\bar {R}}}

\newcommand{\ol} {\overline}
\newcommand{\ok} {{\overline k}}
\newcommand{\ox} {{\overline x}}
\newcommand{\oy} {{\overline y}}

\newcommand{\hvar}{{\widehat {\varphi}}}
\newcommand{\hpsi}{{\widehat {\psi}}}
\newcommand{\halpha}{{\widehat {\alpha}}}
\newcommand{\hsi}{{\widehat {\sigma}}}
\newcommand{\hmu}{{\widehat {\mu}}}
\newcommand{\widehatpar}{{\widehat {\partial}}}
\newcommand{\ha}{{\widehat {a}}}
\newcommand{\hf}{{\widehat {f}}}
\newcommand{\hg}{{\widehat {g}}}
\newcommand{\hh}{{\widehat {h}}}
\newcommand{\hp}{{\widehat {p}}}
\newcommand{\hr}{{\widehat {r}}}
\newcommand{\hx}{{\widehat {x}}}
\newcommand{\hz}{{\widehat {z}}}
\newcommand{\hy}{{\widehat {y}}}
\newcommand{\hA}{{\widehat {A}}}
\newcommand{\hc}{{\widehat {C}}}
\newcommand{\hE}{{\widehat {E}}}
\newcommand{\hG}{{\widehat {G}}}
\newcommand{\hH}{{\widehat {H}}}
\newcommand{\hI}{{\widehat {I}}}
\newcommand{\hM}{{\widehat {M}}}
\newcommand{\hR}{{\widehat {R}}}
\newcommand{\hX}{{\widehat {X}}}
\newcommand{\hY}{{\widehat {Y}}}

\newcommand{\dM}{{\dot {M}}}
\newcommand{\dD}{{\dot {D}}}
\newcommand{\dE}{{\dot {E}}}
\newcommand{\dX}{{\dot {X}}}

\newcommand{\dimc}{{{\dim}_{\C}\,}}
\newcommand{\De}{{Deg \,}}
\newcommand{\fr} [1] {{F^{#1}A}}
\newcommand{\po}{{\pi_0}}
\newcommand{\tors}{{\rm Tors ~}}
\newcommand{\Tor}{{\rm Tor ~}}
\newcommand{\Pic}{{\rm Pic ~}}
\newcommand{\Ker}{{\rm Ker ~}}
\newcommand{\ML}{{\rm ML}}
\newcommand{\Dk}{{\rm Dk}}
\newcommand{\LND}{{\rm LND}}
\newcommand{\LNDG}{{\rm LND_{gr}}}
\newcommand{\GL}{{\rm GL}}
\newcommand{\Gr}{{\rm Gr ~}}
\newcommand{\gr}{{\rm gr ~}}
\newcommand{\Div}{{\rm Div ~}}

\newcommand{\mod} {{\rm mod ~}}
\newcommand{\chr} {{\rm char ~}}
\newcommand{\grad} {{\rm grad ~}}
\newcommand{\spec} {{\rm spec ~}}
\newcommand{\reg} {{\rm reg ~}}
\newcommand{\dt} {{\rm det ~}}

\newcommand{\dg} {{\rm deg ~}}
\newcommand{\trdg} {{\rm tr.deg ~}}
\newcommand{\dgv} {{\rm deg}_{\varphi} ~}
\newcommand{\dgy} {{\rm deg}_y ~}
\newcommand{\dge} {{\rm deg}_1 ~}
\newcommand{\dgp} {{\rm deg}_\partial ~}

\newcommand{\ldm} {{\rm L-dim ~}}
\newcommand{\kdm} {{\rm K-dim ~}}
\newcommand{\dm} {{\rm dim ~}}
\newcommand{\cdm} {{\rm codim ~}}

\title
%[~Affine modifications]
{Affine modifications and affine hypersurfaces with a very 
transitive automorphism group}

\vspace{.5cm}

\author{~Sh. ~Kaliman\thanks{Partially supported by the 
NSA grant MDA904-96-01-0012.},
%\address 
%Department of Mathematics and Computer Science \\
%University of Miami \\ 
%Coral Gables, FL  33124, U.S.A.  
%\email kaliman@math.miami.edu  
%\and
~M. ~Zaidenberg
%\address 
%Institut Fourier de Math{\'e}matiques, UMR 5582 \\
%Universit{\'e} Grenoble I\\
%BP 74, 38402 St. Martin d'H{\`e}res--c{\'e}dex, France
%\email zaidenbe@ujf-grenoble.fr
}

\date{}
\maketitle

\begin{abstract} 
We study a kind of modification of an affine domain which produces 
another affine domain. First appeared in passing in the 
basic paper of O. Zariski \cite{Zar}, it was further considered by 
E.D. Davis \cite{Da}. In \cite{Ka 1} its geometric 
counterpart was applied to construct contractible smooth affine varieties non-isomorphic to Euclidean spaces.  
Here we provide certain conditions (more general then those in 
\cite{Ka 1}) which guarantee preservation of the topology under a modification.
 
As an application, we show that the
group of biregular automorphisms of the 
affine hypersurface $X \subset \C^{k+2}$ given by the equation 
$uv=p(x_1,\dots,x_k)$ where $p \in \C[x_1,\dots,x_k],$ acts 
$m-$transitively on the smooth part reg$\,X$ of $X$ for any $m \in \N.$ 
We present examples of such hypersurfaces diffeomorphic to Euclidean spaces. 
\end{abstract}

\tableofcontents

\section*{Introduction}

It is well known (and elementary) that for $n > 1$ the automorphism group Aut$\,\C^n$ of the 
affine space $\C^n$ (or, which is the same, of the polynomial ring in $n$ variables $\C^{[n]} := \C[x_1,\dots,x_n]$) acts 
$m-$transitively on $\C^n$ for any 
$m \in \N.$ That is, the diagonal action of the group Aut$\,\C^n$ on the
$m-$th symmetric power $S^m\C^n,$ and even on the
$m-$th Cartesian power $(\C^n)^m,$ is transitive outside of the diagonals for any\footnote{See \cite{RoRu} for an 
$\infty-$transitivity property of the group of analytic automorphisms of $\C^n.$} $m \in \N.$ Clearly, every Zariski open 
subset of the form $\C^n \setminus K$ where $K \subset \C^n$ 
is a finite set of points, possesses this property. 

Let $V$ be a compact complex space. By a generalized 
Bochner-Montgomery Theorem \cite[(2.3)]{Akh}, 
the automorphism group Aut$\,V$ is a complex Lie group. Therefore, it cannot be 
$m-$transitive on $V$ for $m > $dim Aut$\,V.$ Thus, the question arises:

%%\smallskip

{\it For which complex manifolds or, at least, for which quasi-projective 
varieties $X$ the group Aut$\,X$ is $m-$transitive on $X$ for any $m \in \N?$ }

%%\smallskip

\noindent Or, more restrictively (to exclude the above examples of type 
$\C^n \setminus K$) 

%%\smallskip

{\it Which affine algebraic varieties possess this property?}

%%\smallskip

\noindent In sections 4 and 5 we describe a class 
of affine hypersurfaces 
non-homeomorphic, in general, to the affine spaces, 
with very transitive automorphism groups (see The Transitivity Theorem 
in sect. 5). 

One more remark is in order. One might specify the above problem by asking
whether, for $m$ given, there exists an $m-$transitive algebraic
group action on $X.$ For $m$ sufficiently large with respect to the dimension
of $X,$ the class of such varieties $X$ seems to be rather poor. Actually, 
already the affine plane $\C^2$ does not admit an $m-$transitive algebraic
group action for $m \ge 3.$ Indeed, by the Jung-van der Kulk 
Theorem \cite{Ju, vdK}, the group
Aut$\,\C^2$ can be represented as the amalgamated product of 
the affine group 
Aff$\,\C^2$ and the Jonqui{\`e}re subgroup 
$J(\C^2) \subset$ Aut$\,\C^2$ of the triangular transformations 
$(x,\,y) \longmapsto (x,\,y+p(x))$ where 
$p \in \C[x].$ 
By the theorem of Serre \cite{Se}, any subgroup of finite length
$G \subset $Aut$\,\C^2$ is conjugate either with a subgroup of Aff$\,\C^2$ 
or with a subgroup of $J(\C^2).$ Since the polynomial $x$ is an invariant 
of the Jonqui{\`e}re subgroup $J(\C^2),$ in the latter case $G$ possesses a
non-constant invariant function, and therefore, it is not even transitive on
$\C^2.$ In the former case, $G$ is at most $2-$transitive; indeed, the 
collinearity of a triple of points is preserved by the action of 
the affine group Aff$\,\C^2.$ It remains to observe following 
D. Wright \cite{Wr 1} 
that any algebraic subgroup $G \subset $Aut$\,\C^2$ is of finite length. 

%%\smallskip

The proof of Theorem \ref{thm6.1} on $\infty-$transitivity of the 
automorphism group of a smooth hypersurface $X \subset \C^{k+2},\,\,k \ge 2,$ 
given by the equation $uv-p({\ol x})=0$ goes as follows. 
This hypersurface $X$
can be naturally presented in two different ways as an affine 
modification of $\C^{k+1}$ (see below). Due to Corollary \ref{cor2.2}, 
one can lift to $X$ 
those automorphisms of the affine space $\C^{k+1}$ which preserve the locus of
at least one of these modifications. 
In such a way we obtain two subgroups $\hG_1$ and $\hG_2$ 
of the automorphism group Aut$\,X.$ It turns out that  
the subgroup $\hG$ of Aut$\,X$ generated by $\hG_1$ and $\hG_2$ 
acts $m-$transitively on $X$ for any $m \in \N.$ 

While a general hypersurface $X\subset \C^{k+2}$ as above has quite a 
rich topology, in sections 4 and 6 we give 
a series of examples of hypersurfaces 
of this type which are diffeomorphic to the affine space 
$\C^{k+1}, \,\,k \ge 3.$ 
Presumably, they are not, in general, isomorphic to $\C^{k+1},$ and so, 
they should present new exotic algebraic structures on $\C^{k+1}$ 
(see e.g. \cite{Za 2}). As well, this would show that the Miyanishi's
characterization of the affine 3-space ${\bf A}^3$ \cite{Miy}
cannot be applied to 
${\bf A}^n$ for $n \ge 4$  (see Remark \ref{rm5.3}). 
But if one of them were isomorphic to $\C^{k+1}$
this would answer in negative alternatively, 
either to the Zariski Cancellation Problem 
or to the Abhyankar-Sathaye Embedding Problem (see below). 

Sections 1-3 are devoted to a kind of modification which acts on the affine
rings producing new affine rings; we call it {\it affine modification}. 
It appeared in passing in the classical Zariski paper \cite{Zar} 
(see also \cite[III.2]{Hiro}), but apparently, the first proper study of 
this transform  
was done in \cite{Da}. 
A geometric counterpart of the affine modification occurred to be useful 
in constructing exotic algebraic structures on the affine spaces,
as it was done in \cite{Ka 1} (see also \cite{Za 2}).   

In Corollary \ref{cor2.2} we lift automorphisms to an affine modification. 
In section 3 we provide certain conditions (more general then those in 
\cite{Ka 1}) which guarantee preservation 
of the topology under a modification.

%%\smallskip

Recall the Abhyankar-Sathaye Embedding Problem: 

%%\smallskip

\noindent {\it Is every closed embedding
of $\C^k$ into $\C^n$ rectifiable, i.e. equivalent to a linear one up 
to the actions of the automorphisms groups Aut$\,\C^k$ resp. Aut$\,\C^n?$}

%%\smallskip

\noindent 
In section 7 we give a generalization of the Sathaye-Wright Theorem 
\cite{Sat, Wr 2} which guarantees rectifiability of the special embeddings 
$\C^2\hookrightarrow \C^3$ given by the equations of the form 
$f(x,\,y)z^n + g(x,\,y) = 0.$ Notice that for $n = 1$ such a surface is 
an affine modification of the affine plane $\C^2$ along the divisor $D_f=(f)$
with center at the ideal $I = (f,\,g) \subset \C[x,\,y].$ 
More generally, in Theorem \ref{thm4.2} 
we show that any smooth acyclic surface 
$X \subset \C^3$ given by the equation  $f(x,\,y)z^n + g(x,\,y) = 0$ is
isomorphic to $\C^2$ (and can be rectified). This is no longer true for 
$\Q-$acyclic surfaces; see Example \ref{ex4.1}. 

\section{Affine modifications}

Consider a triple $(A,\,I,\,f)$ consisting of a commutative  
ring A with unity, an ideal $I$ of $A,$ and an element $f$ of $I$ 
which is not a zero divisor.
We call it a {\it Noetherian triple} if $A$ is Noetherian, and  an 
{\it affine triple} (over $\C$) if $A$ is an affine domain /$\C,$ that is, 
a finitely generated commutative $\C-$algebra which is a domain. 
With $t$ being a new symbol, denote by
$$A[It] = A \oplus \bigoplus_{n=1}^{\infty} (It)^n \simeq 
A \oplus I\oplus I^2 \oplus \dots = Bl_I(A) 
%\simeq \mbox{Sym}_A I
$$
the blowup algebra, or the Rees algebra \cite[\S 5.2]{Ei}, \cite{Va}.

\smallskip

\noindent \defin By the {\it affine modification $\Sigma_{I,\,f}(A)$ 
of the ring $A$ along $(f)$ with center $I$} 
(or, shortly, with the locus $(I,\,f)$) we mean the quotient of the
blowup algebra $A[It]$ by 
the principal ideal generated by the element $1 - ft \in A[It]:$
$$\Sigma_{I,\,f}(A) =  A[It]/(1 - ft)\,.$$
When $f$ and $I$ are fixed, without abuse of notation we denote 
$A' = \Sigma_{I,\,f}(A).$ Clearly, if $(A,\,I,\,f)$ is a Noetherian triple, 
then 
$A'$ is again a commutative Noetherian ring with unity. 

\label{def1.1}
\smallskip

Our purposes in this paper are mainly of geometric 
nature\footnote{By this reason, we do not consider here 
possible generalizations, for instance,
such as replacing the filtration of $A$ by the powers 
$\{I^n\}_{n \in \N}$ of the ideal $I$ (resp. the multiplicative 
system $\{f^n\}_{n \in \N} \subset I$) by a more general one.}. 
Thereby, we adopt the following 

\smallskip 

\noindent {\bf Conventions.} 
Hereafter $(A,\,I,\,f)$ is assumed to be an affine triple /$\C.$
Besides, we fix two systems of generators
$a_1,\dots,a_r$ of the algebra $A$ resp. $b_0 = f,\,b_1,\dots,b_s$ 
of the ideal $I.$ Denoting $\C^{[r]}$ the polynomial ring in $r$ 
variables, consider 
the surjective homomorphisms 
$$\var\,:\,\C^{[r]}:=\C[x_1,\dots,x_r] \longrightarrow\!\!\!\to A,\,\,\,x_i 
\longmapsto a_i,\,\,\,i=1,\dots,r,$$ resp.
$$\var_I\,:\,\C^{[r+s+1]}=\C[x_1,\dots,x_r,\,y_0,\dots,y_s] \to A[It]\,,$$ 
$$x_i \longmapsto a_i,\,\,\,y_j \longmapsto b_jt,
\,\,i=1,\dots,r,\,\,j=0,\dots,s.$$
Let $X$ resp. $Y$ be the image of the associated closed embedding
$\hvar\,:\, $spec$\,A \hookrightarrow \C^r$ resp. 
${\widehat \var_I}\,:\, $spec$\,A[It] \hookrightarrow \C^{r+s+1}.$
Then $X$ and $Y$ are reduced irreducible affine varieties
\footnote{We do not suppose the divisor $D_f$
resp. the affine subscheme spec$\,A/I$ in $X$ being reduced or irreducible;
that is, the ideals $(f)$ and $I$ are not assumed being radical or primary.}     
(indeed, $A[It] \subset A[t]$ is an affine domain /$\C$). The principal ideal 
$(1-ft) \subset A[It]$ being prime\footnote{Indeed, the principal ideal 
generated by the regular function $1-ft$ in the algebra 
$A[t] = \C[X \times \C]$ is prime. That is, $p(t)q(t) = (1-ft)r(t)$ where
$p,\,q,\,r \in A[t],$ implies that $(1-ft)$ divides $p(t)$ or 
$(1-ft)$ divides $q(t)$ in
$A[t].$ But if, say, $(1-ft)$ divides $p(t)$ in $A[t],$ and 
$p \in A[It],$ then, 
as it is easily seen, $(1-ft)$ divides $p(t)$ in $A[It].$}, 
the $\C-$algebra $A'$ is an affine domain /$\C,$ too. Thus, spec$\,A'$ is a 
reduced irreducible affine variety isomorphic to the hypersurface 
$V(1-ft) \subset $spec$\,A[It]\simeq Y.$

Denote by $\rho\,:\,A[It] \longrightarrow\!\!\!\to A' = A[It]/(1-ft)$ 
the canonical surjection. 
Since $\var_I(1-ft) = 1- y_0,$ the image $X'$ of the closed embedding 
${\widehat {\var_I\circ \rho}}\,:\, $spec$\,A' \hookrightarrow \C^{r+s+1}$
coincides with the hyperplane section $X' = H' \cdot Y \simeq V(1-ft)$ 
of $Y$ defined by the affine hyperplane $H' := \{y_0 = 1\} \simeq \C^{r+s}.$ 

\smallskip

\noindent \defin We call $X'\subset \C^{r+s}$ as above the 
{\it affine modification of the affine variety 
$X\subset \C^{r}$ with the locus} $(I,\,f)$
(or in other words, {\it along the divisor $D_f$ with center $I$});
we denote $X' = \Sigma_{I,\,f} (X).$  

%%\smallskip

Thus, this definition takes into account the distinguished
systems of generators $a_1,\dots,a_r$ of the algebra $A$ 
resp. $b_0 = f,\,b_1,\dots,b_s$ of the ideal $I,$ that is, 
the closed embeddings 
spec$\,A \simeq X \hookrightarrow \C^{r}$ and 
spec$\,A' \simeq X' \hookrightarrow \C^{r+s}.$

\label{def1.2}
\smallskip
{%\footnotesize
\noindent \rems 

\noindent \nrem For a Noetherian triple $(A,\,I,\,f)$ 
and a fixed 
system of generators $b_0 = f,\,b_1,\dots,b_s$ of the 
ideal $I$ one may
consider the closed embedding 
${\widehat \beta}\,:\,$ spec$\,A' \hookrightarrow 
$ spec$\,A \times \C^s$ associated with the surjective 
homomorphism of the polynomial algebra 
$A^{[s]} := A[y_1,\dots,y_s]$ over $A :$
$$\beta\,:\,A^{[s]} \longrightarrow\!\!\!\to A',\,\,\,
\beta(y_i) = \rho(b_it),\,\,\,i=1,\dots,s\,.$$
For the affine triple $(A,\,I,\,f)$
the embedding ${\widehat \var_I}\,:\, $spec$\,A' 
\hookrightarrow \C^{r+s}$
as above is the composition 
${\rm spec}\, A' \stackrel{{\widehat \beta}}{\hookrightarrow}   
{\rm spec}\, A \times \C^s \stackrel{\hvar \times 
{\rm id}}{\hookrightarrow} \C^{r} \times \C^{s}\,.$ 
\label{rm1.1}
%$$h_1({\overline x},\,{\overline y}),\dots,
%h_N({\overline x},\,{\overline y}) 
%\in \C^{[r+s+1]}$$ 
%which are polynomials homogeneous in ${\overline y}.$ 
%Then the ideal $I(X') \subset \C^{[r+s]}$ 
%is generated by the polynomials 
%${\widetilde h}_1({\overline x},\,{\overline y}'),
%\dots,{\widetilde h}_N({\overline x},\,{\overline y}'),$ 
%where
%$${\widetilde h}_i(x_1,\dots,x_r,\,y_1,\dots,y_s) := 
%h_i(x_1,\dots,x_r,1,\,\,y_1,\dots,y_s).$$ 

\nrem \label{rm1.2'} 
Denote $Y^* = Y \setminus V_Y$ where 
$Y \subset X \times \C^{s+1} \subset \C^{r+s+1}$
is the affine cone introduced above 
with the vertex set $V_Y:=X \times \{{\overline 0}\}.$   
Then the blow up 
$$Z = Bl_I X = {\rm Proj}_A Bl_I A = 
{\rm Proj}_A A[It] \subset X \times \pr^s \subset \C^r\times \pr^s$$ 
of $X$ with center $I$ is the canonical projection of $Y^*,$ that is, 
$Z = Y^*/\C^*$ \cite[\S 5.2]{Ei}. 
This projection restricted to the hyperplane section $X' = H' \cap Y$ 
yields an isomorphism
$X' \stackrel{\simeq}\longrightarrow Z \setminus H_0$ where 
$H_0 \subset \pr^s$ is the coordinate hyperplane $y_0=0.$  

We denote by $E$ the exceptional divisor\footnote{By abuse of language, 
usually by the exceptional divisor we mean its support.} 
$\si_I^{-1}(C)\subset Z$ where 
$C:=V(I) \subset X$ and $\si_I\,:\,Z \to X$ is the blowup morphism,  
and by $E'$ its affine part $E'=E \setminus H_0 \subset X'.$    

%%\smallskip

\noindent \nrem \label{rm1.3'}
In the sequel we distinguish between two notions of proper transform 
of the divisor $D_f;$ for that we use two terms, proper transform
and strict transform, respectively. 
Namely, by the {\it proper transform} $D_f^{pr}$ of $D_f$ in $Z$ we mean  
the set of homogeneous prime ideals $p \in Z=$Proj$\,_A A[It]$ such that 
$ft \in p.$ Thus,
supp$\,D_f^{pr} = \{y_0=ft=0\} = H_0 \cap Z,$ 
and so,\footnote{By abuse of notation, 
we write $Z \setminus D_f^{pr}$ 
instead of $Z \setminus $supp$\,D_f^{pr}.$} 
$X' \simeq Z \setminus H_0 = Z  \setminus D_f^{pr}.$ 
Whereas the {\it strict transform} $D_f'$ of $D_f$ in $Z$ 
we understand as usually, i.e. as the closure in $Z$ of 
the preimage $\si_I^{-1}(D \setminus V(I)).$ In general, 
these two transforms are different;
see Examples \ref{ex1.1} and \ref{ex3.4} below. 

%%\smallskip

\noindent \nrem \label{rm1.5}
In the case when $D := D_f$ is a reduced effective divisor
and $I$ is the (radical) ideal of a closed reduced subvariety 
$C = V(I) \subset {\rm reg} D$ in $X,$ the affine modification 
$X' = \Sigma_{I,\,f}(X) =: \Sigma_{C,\,D}(X)$
coincides with the one considered in \cite{Ka 1}.}

\smallskip

\noindent {\bf Notation.} Denote ${\rm Frac}\, A$ 
the field of fractions of a domain $A.$ For an affine  
triple $(A,\,I,\,f),$ 
$A[I/f]$ denotes the subalgebra of the field ${\rm Frac}\, A$ 
generated over
$A$ by the elements $g/f$ where $g$ runs over the ideal $I.$ 
Under the above 
conventions we have 
$$A[I/f] = A[b_1/f,\dots,b_s/f] = \{a/f^k \in {\rm Frac} \,
A\,|\,a \in I^k,\,\,k=0,\,1,\dots\}\,.$$ 
For an affine variety $X$ we denote by $\C[X]$ the algebra 
of regular functions on $X,$ and by $\C(X)$ the rational 
function field Frac$\,\C[X]$ of $X.$

\smallskip

\noindent \nprop{{\bf The first properties of the affine modification.}} 

\noindent {\it (a) In the notation as above,
the affine modification $X'= \Sigma_{I,\,f} (X)$ is isomorphic to the 
complement $Z \setminus D_f^{pr}.$

%%%\smallskip

\noindent (b) 
There is a canonical isomorphism 
$\alpha\,:\,A' \stackrel{\simeq}\longrightarrow 
A[I/f]= A[b_1/f,\dots,b_s/f]$ which sends $\rho(A)\simeq A$ 
isomorphically onto $A$ and $\rho(b_it)$ into $b_i/f,\,\,i=1,\dots,s.$ 
Thus, $X' = {\rm spec}\ A' \simeq {\rm spec}\ A[I/f].$

%%%\smallskip

\noindent (c)} (E.D. Davis \cite{Da}) {\it Consider the 
surjective homomorphism $$\beta\,:\,A^{[s]}=A[y_1,\dots,y_s] 
\longrightarrow\!\!\!\to A[I/f] = A[b_1/f,\dots,b_s/f] \simeq A'$$ 
where $\beta(y_i)= 
b_i/f,\,\,\,i=1,\dots,s.$
Denote by $I'$ the ideal of the polynomial algebra $A^{[s]}$ 
generated by the elements $L_1,\dots,L_s \in {\rm ker}\,\beta$ 
where $L_i = fy_i - b_i.$
Then ${\rm ker}\,\beta=I'$ (that is, the subvariety 
$X' \subset X \times \C^s$
is defined by the equations $fy_i - b_i = 0,\,\,i=1,\dots,s$) 
iff $I'$ is a prime ideal, i.e. spec$\,A^{[s]}/I'$ is a reduced 
irreducible 
subvariety in $X \times \C^s.$ The latter is true, for instance, if
the system of generators $b_0 = f,\,b_1,\dots,b_s$ of the ideal $I$ 
is regular\footnote{I.e. for each $i=1,\dots,s$ the image of $b_i$ 
is not a zero divisor in $A/(b_0,\dots,b_{i-1}).$}. 

%%%\smallskip

\noindent (d) 
Furthermore, if ${\rm ker}\,\beta=I'$
then we have an isomorphism $E' \simeq C \times \C^s$ 
where $E'\subset X'$ is the exceptional divisor and $C = V(I).$}

\label{pr1.1}
\smallskip

\noindent {\it Proof.~}  For the proof of (a) see Remarks \ref{rm1.2'} 
and \ref{rm1.3'} above,
and for that of (c) see \cite[Prop. 2]{Da}, \cite[(A.6.1)]{Ful} or 
also \cite{Mik}. The statement (d) immediately follows from (c). 
Indeed, by (c), we have 
$$E' = \si_I^{-1}(C) = \{(x,\,\by) \in X\times \C^s\,|\,b_i(x) =
 0,\,\,i=0,\dots,s\}=C \times \C^s \subset X \times \C^s\,.$$
(b) On $X' = H' \cap Y$ we have $t = 1/f,$ and 
$y_i = b_it,\,\,i=1,\dots,s.$ Thus, 
$A' = \C[X'] = A[1,\,y_1,\dots,y_s]\,|\,X' = A[b_1/f,\dots,b_s/f],$ 
as stated. 
\qed

\smallskip

It is well known (see e.g. \cite[7.17]{Ha}) that every birational
projective morphism of quasiprojective varieties $Y \to X$ 
is a blow up of $X$ with center at a subsheaf of ideals
${\cal I} \subset {\cal O}_X.$ Similarly, we have the following theorem.

\smallskip

\noindent \nthm{{\bf Birational morphisms as affine modifications.}} 

\noindent {\it Any birational morphism $q\,:\,Y \to X$ 
of affine varieties is an affine modification. More precisely, 
there are an ideal $I \subset A = \C[X],$ an element $f \in I$ 
and an isomorphism 
$\alpha\,:\,Y \stackrel{\simeq}\longrightarrow X' := \Sigma_{I,\,f}(X)$ 
such that $q = \sigma_I \circ \alpha$ where $\sigma_I\,:\,X' \to X$ 
is the blowup morphism with center $I.$ }

\label{thm1.1}
\smallskip

\noindent {\it Proof.~}  Denote $A_1 = (q^*)^{-1}(\C[Y]) \subset \C(X) = $Frac$\,A,$
where $(q^*)^{-1}\,:\,\C(Y) \to \C(X)$ is the induced isomorphism of 
the rational function fields. Each function $a \in A=\C[X]$ 
being lifted by 
$q^*$ to a function in $\C[Y]$ comes back under the inverse 
homomorphism $(q^*)^{-1}.$
Therefore, $A \subset A_1 \subset $Frac$\,A.$ It is enough to show 
that $A_1 = A[I/f]$ for some ideal $I$ of $A$ and some $f \in I.$ 

The affine domain $A_1$ being  finitely generated
we have $$A_1 = A[a_1/f_1,\dots,a_k/f_k] =  A[b_1/f,\dots,b_k/f] = A[I/f]$$ 
for some $a_i,\,f_i \in A$ where $f = f_1 \cdot\dots\cdot f_k \in A,\,\,
b_i=fa_i/f_i \in A$ and $I = (f,\,b_1,\dots,b_k) \subset A,$ 
as desired.
\qed

\smallskip

%\noindent \rems 
%
%\noindent \nrem Recall that the blowup morphism 
%$\sigma_I\,:\,Z \to X$ coincides with the restriction to $Z$ of 
%the first projection pr$_1\,:\,X \times \pr^s \to X.$ The composition  
%$$X' \hookrightarrow
%Z \subset X \times \pr^s \stackrel{{\rm pr}_1}\longrightarrow X\,$$ 
%associated with the composition of homomorphisms 
%$A \stackrel{i}\hookrightarrow A[It] \stackrel{\rho}\longrightarrow A'\,$
%yields a birational morphism $\sigma_I\,:\,X' \to X.$  
%It coincides with the restriction to $X'$ of the first projection 
%pr$_1\,:\,\C^r \times \C^s \to \C^r.$
%The induced isomorphism of the rational function fields 
%$\C(X) = $Frac$\, A \stackrel{\simeq}\longrightarrow $Frac$\, A' = \C(X')$ 
%canonically identifies the algebra $A'$ with the subalgebra 
%$A[I/f] \subset\C(X)$ (cf. Proposition \ref{pr1.1}(b)).
%
%The isomorphism 
%$$\sigma_I^{-1}\,:\, X \setminus D_f \stackrel{\simeq}
%\longrightarrow X' \setminus E' \hookrightarrow X' \hookrightarrow 
%C^r \times \C^s$$ 
%is given by the formulas 
%$x_i = a_i,\,\,y_j=b_j/f,\,\,i=1,\dots,r,\,\,j=1,\dots,s.$ 
%
%\label{rm1.3}
%%%\smallskip
%
%\noindent \nrem Actually, the prime ideal ker$\,\beta$ as in 
%Proposition \ref{pr1.1}(c) which defines the affine 
%subvariety $X' \subset X \times \C^s$ coincides with the radical rad$\,I',$ 
%and $I' \supset f^r$ker$\,\beta$ for a sufficiently large $r$ \cite{Da}.
%
%\label{rm1.4}
%\smallskip

%An affine modification may possess the following decomposition.

%\smallskip 

\noindent \nprop {{\bf A decomposition of an affine modification.}}
{\it Let $f = f_1f_2 \in I.$ Denote $I_1 = (I,\,f_1)$ and 
$A_1 = \Sigma_{I_1,\,f_1} (A).$ Consider the ideal $I_2$ of
the algebra $A[I_1/f_1]\simeq A_1$  generated by the subspace $I/f_1.$
Then we have 
$A' := \Sigma_{I,\,f} (A) \simeq A_1':=\Sigma_{I_2,\,f_2} (A_1).$

Respectively, denoting $X = $spec$\,A,\,\,\,X' = $spec$\,A'$ and 
$X_1 = $spec$\,A_1$ we obtain the decomposition $\sigma_I = \sigma_{I_2}
\circ \sigma_{I_1}$ where $\sigma_I\,:\,X' \to X,\,\,\,
\sigma_{I_1}\,:\,X_1 \to X$ and $\sigma_{I_2}\,:\,X' \to X_1$ 
are the corresponding blowup morphisms.}

\label{pr1.2}
\smallskip 

\noindent {\it Proof.~}  Let $I=(f,\,b_1,\dots,b_s).$ 
By Proposition \ref{pr1.1}(b), we have: 
$$A' \simeq A[I/f] = A[b_1/f,\dots,b_s/f] = A[b_1/f_1f_2,\dots,b_s/f_1f_2]
\subset {\rm Frac}\, A\,,$$
$$A_1 \simeq A[I_1/f_1] =
A[b_1/f_1,\dots,b_s/f_1] \subset A[I/f]\subset {\rm Frac}\, A\,,$$ and 
$I_1 = (f_2,\,b_1/f_1,\dots,b_s/f_1)\subset A_1.$ Hence, 
$$A_1' \simeq A_1[I_2/f_2]= A_1[b_1/f_1f_2,\dots,b_s/f_1f_2] \subset 
{\rm Frac}\, A\,.$$ 
Clearly, 
the latter subalgebra coincides with $A[I/f] \simeq A',$
as claimed. \qed

\smallskip 

Next we give several examples of affine modifications.

\smallskip

{%\footnotesize

\noindent \nex {{\bf Deleting a divisor.}} \label{ex1.1}
If $I = A,$ that is, the height ht$\,I = 0,$ then we have
$A' \simeq A[1/f],\,\,\, Z \simeq X$ and 
$X' \simeq X \setminus D_f \simeq Z \setminus D_f^{pr}.$ 

If $I=(f)$ is a principal ideal (and so, ht$\,I = 1$) 
we obtain $A[I/f] = A,$ 
whence, $X' = X \simeq Z = X \times \pr^0$ whereas
supp$\,D_f^{pr} = \emptyset.$ Thus, 
the equality $X' = Z \setminus D_f^{pr}$ 
still holds. 

Furthermore, even if ht$\,I \ge 2$ 
it may happen that supp$\,D_f^{pr} = \sigma_I^{-1}($supp$\,D_f)$ and so, 
$X' \simeq X \setminus D_f,$ i.e. 
that the hyperplane section $H_0 \cap Z$ contains the exceptional divisor
$E$ of the blow up $\sigma_I.$ For instance, take 
$A = \C^{[2]}=\C[x,\,y],$ that is, $X = \C^2,\,\,I = (x,\,y)$ and 
$f = x^2 \in I.$ Then we have
$Z = \{((x,\,y),\,(u:v)) \in \C^2 \times \pr^1\,|\,xv = yu\},$ and 
the curve supp$\,D_f^{pr}$ given in $Z$ by the equation $ux = 0$ consists of
two lines, the first one $l':=\{u = 0\}$ being the strict transform in $Z$ 
of the affine line $l:=$supp$\,D_f=\{x=0\}$ 
and the second one $\{x=y = 0\}$ being the exceptional divisor $E \subset Z.$ 
Thus, $X' = Z \setminus (l' \cup E) \simeq
X \setminus l \simeq \C^* \times \C \subset \C^2 \simeq X.$ While for 
the affine triple $(A,\,I,\,x)$ we obtain 
$X' = \Sigma_{I,\,x}(X) \simeq \C^2.$ 

\smallskip

\noindent \nex {{\bf A singular modification of the affine plane.}} \label{ex1.5}
Set $A = \C[x,\,y]$ (that is, $X = \C^2$), $f = x,\,\,\,I=(x,\,y^2).$ 
The system of generators $b_0 := x,\,\,b_1 = y^2$ of the ideal $I$ is regular. 
Hence, by Proposition \ref{pr1.1}(c), $X'$ is the affine 
surface $xz = y^2$ in $\C^3.$ It is isomorphic to 
the quotient of the affine plane $\C^2$ by the involution 
$(x,\,z) \longmapsto (-x,\,-z).$

\smallskip

\noindent \nex 
{{\bf Modification of an affine space along a hyperplane 
with center at a point.}} \label{ex1.4}
Let $A = \C[x_1,\dots,x_n,\,y],$ $f=y,$ 
and let $I=(x_1,\dots,x_n,\,y)$ be the maximal ideal corresponding 
to the origin of $X = \C^{n+1}.$ 
Then 
$A' \simeq \C[x_1y^{-1},\dots,x_ny^{-1},\,y]\subset \C(x_1,\dots,x_n,\,y).$ 
By Proposition \ref{pr1.1}(c), 
$X' \simeq \C^{n+1}$ is given in $\C^{2n+1}$ by the equations 
$x_i = yy_i,\,\,\,i=1,\dots,n,$ and the exceptional divisor 
$E' = E \setminus D_f^{pr} \simeq \C^n$ is given in $X'$ as $y=0.$ 
In the coordinates $(y_1,\dots,y_n,\,y)$ of $X' \simeq \C^{n+1}$ 
the blowup morphism
$\sigma_I\,:\,X' \to X$ is given as $\sigma_I (y_1,\dots,y_n,\,y) =
(yy_1,\dots,yy_n,\,y).$

A modification of $\C^{n+1}$ along a 
coordinate hyperplane $H=D_y \simeq \C^n$ 
with center at a coordinate subspace 
$C=\C^k \subset \C^n$ can be described
in a similar way. 

\smallskip 

\noindent \nex {{\bf Modification of an affine space along a divisor 
with center at a codimension two complete intersection.}} \label{ex1.6}
Set $A = \C^{[r]},$ that is, $X = \C^r,$ and $I =(f,\,g),$
where $f,\,g \in A$ are non-constant polynomials without common factor. 
Then $f,\,g$ form a regular system of generators of the 
ideal $I.$ Thus, by Proposition \ref{pr1.1}(c), the affine modification 
$X'= \Sigma_{I,\,f} (X)$ is the hypersurface in $\C^{r+1}$ with the equation 
$f({\overline x})y - g({\overline x}) = 0$ where ${\overline x} = (x_1,\dots,x_r).$ 
The blowup morphism $\sigma_I\,:\,X' \to X$ is the restriction to $X$ 
of the projection 
$\C^{r+1} \to \C^{r},\,\,\,({\ol x},\,y) \to {\ol x}.$ The exceptional divisor 
$E' \subset X'$ is given in $\C^{r+1}$ by the equations 
$f({\overline x})=g({\overline x}) = 0;$ so, 
$E' \simeq V(I) \times \C.$

\smallskip

\noindent \nex {{\bf The Russell cubic
threefold.}} \label{ex1.7}
In particular, set $A = \C[x,\,z,\,t]$ (that is, $X = \C^3$), 
$f = -x^2$ and $I =(f,\,g)$ where $g=x + z^2 + t^3.$ 
Then the affine modification 
$X' = \Sigma_{I,\,f} (X)$ is the smooth
3-fold $x + x^2y + z^2 + t^3=0$ in $\C^4.$ We call it
the {\it Russell cubic} (see \cite{Ru 1}). 
It birationally dominates the affine space $\C^3$ via the blowup morphism 
$\si_I \,:\,X' \to X \simeq \C^3,\,\,\,\si_I \,:\, (x,\,y,\,z,\,t) 
\longmapsto (x,\,z,\,t).$

In turn, the Russell cubic
threefold $X' \subset \C^4$ is birationally 
dominated\footnote{This observation is due to P. Russell.} by $\C^3.$
Indeed, for any $x\ne 0,\,$ $y$ is expressed
in terms of $z$ and $t;$ whence, the part $\{x\ne0\}$ of the threefold
$X'$ is isomorphic to $\C^2\times \C^{*}.$
The `book-surface' $B:=\{x=0\} \subset X'$ is the product
$\C\times \Gamma_{2, \,3}$ where
$\Gamma_{2, \,3} :=\{z^{2}+t^{3}=0\}\subset \C^2.$ 
Fix a smooth point
$\rho\in \Gamma_{2,\, 3},$
and perform the affine modification
$\sigma'\,:\,X''\longrightarrow X'$ of
$X'$ along $B$ with the center $C:=\C\times\{\rho\}.$
In this way we replace $B$ by a smooth surface $E'\simeq \C^2$
and replace the function $x$ by a function
$h \,:\,X''\longrightarrow \C$ such that all the fibers of
$h$ are smooth reduced surfaces isomorphic to $\C^2.$ 
Using an explicit presentation of $X''$ it can 
be checked that $X''\simeq \C^3$
(see Example \ref{ex7.3} below), 
and so, $\sigma'\,:\,\C^3\simeq X'' \longrightarrow X'$
is a birational (whence, dominant) morphism. 

It is known that the Russell cubic is diffeomorphic to $\C^3$
(see e.g. \cite{Ru 1, Ka 1, Za 2} or Example \ref{ex3.2} below). 
However, by a theorem of Makar-Limanov \cite{ML 1} 
(or also \cite{De, KaML 1, Za 2}), it is not isomorphic to $\C^3.$ 
A smooth affine variety which is diffeomorphic 
but non-isomorphic to $\C^n$ 
is called an {\it exotic} $\C^n.$ Thus, the 
Russell cubic $X'$ is an exotic 
$\C^3$ {\it of sandwich type}, that is, there 
are birational morphisms 
$\C^3 \to X' \to \C^3.$

\smallskip 

\noindent \nex {{\bf Decomposing an affine modification.}} \label{ex1.8}
Let $A = \C^{[2]} = \C[x,\,y],\,\,\,I = (f,\,g)$ where
$f,\,g \in \C[x,\,y]$ and $f = f_1f_2.$ Set $A_1 = \C^{[3]}/(f_1z-g)$ 
(see Example \ref{ex1.6}), and let
$\si_I $ be the restriction of $ \si_{I_2}\circ \si_{I_1}$
to the surface $fz-g=0$
where $\si_{I_1}\,:\,\C^3 \to \C^2,\,\,
\si_{I_1}\,:\,(x,\,y,\,z) \longmapsto (x,\,y),$ 
and $\si_{I_2}\,:\,\C^3 \to \C^3,\,\, (x,\,y,\,z) 
\longmapsto (x,\,y,\,f_2(x,\,y)z)$ 
is the affine modification 
of $\C^3$ along the divisor 
$D_{f_2} \times \C$ with center 
$D_{f_2} \simeq $ spec$\,\C^{[3]}/(f_2, \,z).$}

\section{The universal property of affine modifications}

Here we show that the universal property of the blow ups 
(see e.g. \cite{Hiro})
is still valid for the affine modifications. 

\smallskip

\noindent \nprop {{\bf Lifting a homomorphism to affine modifications.}}
{\it Consider two affine triples $(A,\,I,\,f)$ and 
$(A_1,\,I_1,\,f_1)$ 
and their affine modifications 
$A' = \Sigma_{I,\,f} (A) \simeq A[I/f]$ resp. 
$A_1' = \Sigma_{I_1,\,f_1} (A_1)\simeq A_1[I_1/f_1].$  
Let $\mu\,:\,A \to A_1$ be a homomorphism  
such that $\mu (I) \subset I_1$ and $\mu (f) = \alpha f_1$ where $\alpha \in A_1$ 
is an invertible element.  
Then there exists a unique homomorphism $\mu'\,:\,A' \to A_1'$
which extends $\mu;$ it can be defined as follows: 
$$\mu'\,:\, A[I/f] \ni a/f^k  \longmapsto \alpha^{-k}\mu(a)/f_1^k \in A_1[I_1/f_1].$$
If, in addition, $\mu(A) = A_1$ and 
$\mu (I) = I_1,$ then  $\mu'(A') = A_1'.$}

\label{pr2.1}
\smallskip 

\noindent {\it Proof.~}  
By Proposition \ref{pr1.1}(b) (see also Remark \ref{rm1.3'} above), 
there are canonical isomorphisms 
$A' \simeq A[I/f]$ resp. $A_1' \simeq A_1[I_1/f_1].$ 
Thus, it is enough to show that $\mu$ admits a unique extension
$\mu'\,:\,A[I/f] \to A_1[I_1/f_1],$ which is surjective if so are $\mu$ and
$\mu\,|\,I\,:\,I \to I_1.$ Furthermore, since 
$A_1[I_1/f_1] = A_1[I_1/\alpha f_1],$ 
replacing $\alpha f_1$ by $f_1$ we may suppose in the sequel 
that $\alpha = 1,$ i.e. that $f_1 =\mu(f).$

Consider the natural extension 
${\tilde \mu}\,:\,A[t] \to A_1[\tau],\,\,\,\,\,\,
{\tilde \mu}\,|\,A = \mu,\,\,\,{\tilde \mu}(t)=\tau.$
Since ${\tilde \mu}$ sends the principal ideal $(1-ft) \subset A$ into 
the principal ideal $(1-f_1\tau) \subset A_1,$ it induces a homomorphism
$\hmu\,:\,A[t]/(1-ft) \to A_1[\tau]/(1-f_1\tau).$ In this way, via the 
canonical isomorphisms $A[t]/(1-ft) \simeq A[1/f] \subset {\rm Frac}\, A$ resp. 
$A_1[\tau]/(1-f_1\tau) \simeq A_1[1/f_1]\subset {\rm Frac}\, A_1,$ we obtain
an extension $\hmu\,:\, A[1/f] \to  A_1[1/f_1]$ of $\mu.$ Any such extension 
sends a generic element $a/f^k \in  A[1/f]$ into the element 
$\mu(a)/ \mu(f^k)= \mu(a)/ f_1^k\in  A_1[1/f_1].$ Therefore, 
$\hmu$ is uniquely defined on $A[1/f]$ by the formula 
$\hmu(a/f^k) =\mu(a)/ f_1^k.$ 

Observe that an element $a/f^k \in A[1/f]$ belongs to $A[I/f]$ iff
$a \in I^k.$ Since by assumption, 
$\mu(I^k) \subset I_1^k,\,\,\,k=0,\,1,\dots,$ we have 
$\hmu(A[I/f]) \subset
 A_1[I_1/f_1].$ So, $\mu':=\hmu\,|\,A[I/f]$ is a desired extension, 
 and, clearly, it is unique. 

Now, if $\mu(A) = A_1$ and $\mu (I) = I_1,$ then 
$\mu(I^k) = I_1^k,\,\,\,k=0,\,1,\dots.$ 
Hence, any element $a_1/f_1^k \in A_1[I_1/f_1]$ 
can be written as $\mu'(a/f^k)$ with $a \in I^k$. Therefore, 
in this case $\mu'$ is a surjection. 
\qed

\smallskip 

\noindent \defin
Let $J \subset A$ be a prime ideal such that $f \notin J$. 
We define 
the {\it strict transform} $J^{\rm st}$ of $J$ in $A'=A[I/f]$ 
as follows:
$$J^{\rm st} = \{ a' \in A' \,|\, f^k a' \in J \,\,\, 
{\rm for \, some} \,\,\, k \in \N\}\,.$$
It is easily seen that $J^{\rm st}\subset A'$ is a prime ideal containing $J.$ 

\label{df2.1}
\smallskip

The following statements complete Proposition \ref{pr2.1}. 

\smallskip

\noindent \nprop {{\bf How does a modification affect a subvariety.}} 

\noindent {\it (a) Let $(A,\,I,\,f)$ be an affine triple, and
let $J$ be a prime ideal in $A$
such that $ f\notin J.$ Denote by $\mu\,:\,A \longrightarrow\!\!\!\to A_1:= A/J$ 
the canonical surjection. Set $I_1 = \mu (I),\,\,\,f_1=\mu (f),$ and  
$A_1' =\Sigma_{I_1,f_1}(A_1).$ 
Then $A'/J^{\rm st} \simeq A_1'.$

%%\smallskip

\noindent (b) Furthermore, set $X =$spec$\,A,$ and let $X_1$ 
be an irreducible closed subvariety of $X$
which is not contained in the support of $D_f.$ Denote by $J=I(X_1)$ 
the defining ideal of $X_1= $spec$\,A_1,$ and set
$X' =\Sigma_{I,f}(X),\,\,\,X_1' =\Sigma_{I_1,f_1}(X_1)$ where $I_1,\,\,f_1$ 
are as in (a).
Then the strict transform $J^{\rm st}$ of the ideal $J$ in $A'=\C [X']$ coincides 
with the defining ideal of the strict transform $X^{\rm st}_1$ of $X_1$ in $X',$ 
and the variety $X^{\rm st}_1$ is isomorphic to $X_1'.$}

\label{pr2.2}
\smallskip

\noindent {\it Proof.~}  (a)
Let $\mu' :A' \longrightarrow\!\!\!\to A_1'$ be a surjective extension of $\mu$ 
as in Proposition \ref{pr2.1}.
We have to show that ${\rm ker} \, \mu' =J^{\rm st}.$
Note that $J =  {\rm ker} \, \mu \subset {\rm ker} \, \mu',$ and
$f \notin {\rm ker} \, \mu',$ since $\mu'(f)=
\mu (f) \ne 0$ in $ A_1 \subset A_1'.$
For an element $a' \in {\rm ker }\, \mu'$ 
chose $k \in \N$ so that $f^k a' \in A.$ 
Then $0=\mu' (f^k a')=
{\mu} (f^k a'),$ that is, $f^ka' \in J,$ and hence, $a' \in J^{\rm st}.$ 
Thus, 
${\rm ker} \, \mu' \subset J^{\rm st}.$ 

Vice versa, let  $a' \in J^{\rm st},$ and let $k \in \N$ be such that 
$f^k a' \in J = {\rm ker} \, \mu.$ Thus, $\mu' (f^k a')=
{\mu} (f^k a')=0.$ But ${\mu} (f^k) \neq 0.$ 
Since $A_1$ and $A_1'$ are domains, the equality $\mu' (f^k)\mu'(a')=0$
implies $\mu'(a')=0,$ i.e. $a' \in {\rm ker} \, \mu'.$ Therefore, 
$J^{\rm st} \subset {\rm ker} \, \mu',$ and hence, $J^{\rm st} = {\rm ker} \, \mu'.$ 
This proves (a).

%%\smallskip

\noindent (b) Observe that $g' \in A' =\C [X']$ vanishes
on the strict transform $X^{\rm st}_1$ of $X_1$ iff
so does $f^kg'.$ If $k$ is large enough, then we have 
$f^kg' \in A=\C [X] \hookrightarrow \C [X'],$ and 
$f^kg'\,|\,X^{\rm st}_1 =0$ implies that
$f^kg'\,|\,X_1 =0.$ Thus, $g' \in I(X^{\rm st}_1)$ iff 
$f^kg'\in I(X_1) = J$  for large enough $k,$ i.e. iff $g' \in J^{\rm st}.$ 
Hence, 
$I(X^{\rm st}_1) \subset J^{\rm st}.$

Conversely, let $g' \in J^{\rm st},$ i.e. $f^kg' \in J = I(X_1)$ 
for large enough $k.$ 
Thus, $f^kg'\,|\,X_1 =0,$ which implies that $g'\,|\,X^{\rm st}_1 =0.$ 
Therefore, $J^{\rm st}  \subset I(X^{\rm st}_1),$ or 
$J^{\rm st} = I(X^{\rm st}_1).$
Now (a) provides an isomorphism 
$X^{\rm st}_1 = $spec$\,A'/J^{\rm st} \simeq $spec$\,A_1' = X_1'.$ 
The proof is completed.
\qed

\smallskip

\noindent \ncor {{\bf Restricting affine modification to a subvariety.}} 
{\it Let $X = $spec$\,A$ be an affine variety, and $X_1 = $spec$\,A_1$ be an
irreducible
closed subvariety of $X.$ Fix a proper ideal $I \subset A,$ and let the ideal
$I_1 \subset A_1$ consists of the restrictions to $X_1$ of the elements of $I.$
Fix also an element $f \in I$ such that $f_1 := f\,|\,X_1 \neq 0.$ 
Then there is a unique closed embedding $i'\,:\,X_1' \hookrightarrow X'$ 
making the following diagram commutative:

\begin{picture}(150,95)
\unitlength0.2em
\put(45,25){$X_1'$}
\put(60,25){$\hookrightarrow$}
\put(62,29){{$ i'$}}
\put(75,25){$X^\prime$}
\put(88,25){$\hookrightarrow$}
\put(101,25){$\C^r \times \C^s$}
\put(45,5){$X_1$}
\put(60,5){$\hookrightarrow$}
\put(62,9){{$ i$}}
\put(75,5){$X$}
\put(92,5){$\hookrightarrow$}
\put(108,5){$\C^r\,\,\,,$}
\put(47,22){$\vector(0,-1){11}$}
\put(50,16){{$\rm \sigma_1$}}
\put(77,22){$\vector(0,-1){11}$}
\put(80,16){{$\rm \sigma$}}
\put(110,22){$\vector(0,-1){11}$}
\put(113,16){{$\rm pr_1$}}
\end{picture}

\noindent where $\sigma\,:\,X' \to X$ resp. $\sigma_1\,:\,X_1' \to X_1$ 
is the blowup morphism of the affine modification $X' = \Sigma_{I,\,f} (X)$ resp. 
$X_1'= \Sigma_{I_1,\,f_1} (X_1),$ $i\,:\,X_1 \hookrightarrow X$ 
is the identical embedding, and pr$_1$ is the first projection.}

\label{cor2.1}
\smallskip
{%\footnotesize
\noindent \rem In particular, affine modifications commute 
with direct products. 

\smallskip

\noindent \nex {{\bf Modification of an affine hypersurface 
along a hyperplane section 
with center at a point.}} Let $g \in \C[x_1,\dots,x_n,\,y]$ 
be an irreducible polynomial such that $g({\overline 0})=0.$ 
Set $X = \C^{n+1}$ and
$X_1 = g^{-1}(0) \subset X.$ That is, $X = $spec$\,A$ resp. 
$X_1 = $spec$\,A_1$ where $A = \C^{[n+1]}$ resp. 
$A_1 := \C^{[n+1]}/J,\,\,\,J:=(g).$ 
Set $I = (x_1,\dots,x_n,\,y)\subset A,$ 
$I_1 = (x_1 + (g),\dots,x_n+ (g),\,y+ (g))\subset A_1,$ $f=y \in I$ and 
$f_1 = y+ (g) \in I_1.$ As in Example \ref{ex1.4} above consider 
the new coordinates $(y_1,\dots,y_n,\,y)$ in $X' \simeq \C^{n+1}$ where 
the blowup morphism $\sigma_I\,:\,X' \to X$ is given by 
$x_i = yy_i,\,\,i=1,\dots,n.$ 
Represent $g(yy_1,\dots,yy_n,\,y) = y^{\mu}g_1(y_1,\dots,y_n,\,y)$
where $g_1 \in \C[y_1,\dots,y_n,\,y]$ is not divisible by $y.$ 
Then the equation of the affine modification $X_1' = X^{\rm st}_1$ of 
$X_1$
along the divisor $D_{f_1}$ with center $I_1$ in 
$X'\simeq \C^{n+1} = \Sigma_{I,\,f}(X)$ is $g_1=0.$ 

\label{ex2.1}
\smallskip

\noindent \nex {{\bf Modification of an affine space
along a divisor with center at a point} (cf. \cite[Cor. 6.1]{Za 2}).} Let 
$p \in \C^{[k]}$ be a non-constant polynomial such that 
$p({\overline 0})=0.$ Represent $\C^k$ as the hypersurface $X$
in $\C^{k+1}$ with the equation $y-p(\bx)=0.$ Then as above, the
affine modification $X'= \Sigma_{({\overline 0}), p} \C^k = 
\Sigma_{({\overline 0}), y} X$ is isomorphic to the hypersurface  
in $\C^{k+1}$ with the equation ${p(y\bx)-y \over y}=0.$}

\smallskip

\noindent \ncor {{\bf Lifting automorphisms to an affine modification.}} 
\label{cor2.2}
\noindent {\it (a) Let $(A,\,I,\,f)$ be an affine triple, and 
let $\var \in $Aut$\,A$ 
be an automorphism such that $\var(I)=I$ and $\var(f)= \alpha f$ where 
$\alpha \in A$ is an invertible element. 
Then there exists a unique extension 
$\var' \in $Aut$\,A'$ of 
$\var$ to an automorphism of the affine modification 
$A' = \Sigma_{I,\,f} (A).$ 

%%\smallskip

\noindent (b) Furthermore, 
the induced automorphism $\hvar$ of $X = $spec$\,A$ preserves the
divisor $D = D_f$ and the center $C=V(I)$ of the blow up. 
It can be lifted in a unique way to an
automorphism ${\widehat \var'}$ of the affine modification 
$X' = \Sigma_{I,\,f} (X)$ 
which preserves the exceptional divisor $E' \subset X'$ of the blow up; in
the complement $X' \setminus E'$ we have ${\widehat \var'} = \sigma_I^{-1} 
\circ \var'\circ\sigma_I$ where $\sigma_I\,:\,X' \to X$ is 
the blowup morphism.}

\smallskip

{%\footnotesize
\noindent \rem \label{rm2.2'} In particular,
let $(X,\,C,\,D)$ be a triple as in Remark \ref{rm1.5} above, and 
let $\hvar \in $Aut$\,X$ be an automorphism of $X$ such that $\hvar(D) = D,\, 
\hvar(C) = C.$ Then there exists a unique automorphism $\hvar' \in $Aut$\,X'$
of the affine modification $X' = \Sigma_{C,\,D} (X)$ such that 
$\hvar'\,|\,(X' \setminus E') = \hvar\,|\,(X \setminus D)$ under the natural identification 
$\si_C\,|\,(X' \setminus E')\,:\,X' \setminus E' 
\stackrel{\simeq}{\longrightarrow} X \setminus D.$ }

\smallskip

\noindent {%\footnotesize 
\ex \label{ex2.2} 
Let $A = \C^{[3]} = \C[x,\,y,\,z],$ that is, 
$X = \C^3,\,\,\,I = (x,\,z),$ and $f = x.$ Then $X' \simeq \C^3,$ and 
$\sigma_I^*\,:\,A \to A'=\C[x',\,y',\,z']$ is given as 
$$
\sigma_I^*\,:\,\,
%\begin{array}{l}%
%x \longmapsto x'\\[2pt]
%y \longmapsto y'\\[2pt]
%z \longmapsto x'z'
%\end{array}
(x,\,y,\,z) \longmapsto (x',\,y',\,x'z')
$$
(cf. Examples \ref{ex1.4}, \ref{ex1.6}). Let $\mu \in $Aut$\,A$ be given by
$$
\mu\,:\,\,
%\begin{array}{l}%
%x \longmapsto x\\[2pt]
%y \longmapsto y+xg_1(x,\,y,\,z)\\[2pt]
%z \longmapsto z+xg_2(x,\,y,\,z)\,,
%\end{array}
(x,\,y,\,z) \longmapsto (x,\,y+xg_1(x,\,y,\,z),\,z+xg_2(x,\,y,\,z))
$$
where $g_1,\,g_2 \in \C^{[3]}.$ Then we have $\mu(I)=I,\,\,\,\mu(f)=f,$ 
and the extension $\mu' \in $Aut$\,A'$ can be given as 
$$
\mu'\,:\,\,
%\begin{array}{l}%
%x' \longmapsto x'\\[2pt]
%y' \longmapsto y'+x'g_1(x',\,y',\,x'z')\\[2pt]
%z' \longmapsto z'+g_2(x',\,y',\,x'z')\,.
%\end{array}
(x',\,y',\,z') \longmapsto (x',\,y'+x'g_1(x',\,y',\,x'z'),\,z'+g_2(x',\,y',\,x'z'))\,.
$$}
Recall the following notion.

\smallskip 

\noindent \defin Let $A$ be a commutative algebra over $\C,$ and let
$\partial$ be a derivation  of $A.$ It is called {\it locally nilpotent} 
if for any $a \in A$
we have $\partial^n(a)=0$ for some $n = n(\partial,\,a).$ 
For an algebra $A$ denote by 
LND$(A)$ the set of all locally nilpotent derivations 
(LND-s for short) of $A.$
Giving an LND $\partial$ on a Noetherian $\C-$algebra $A$ 
is the same as giving a regular $\C_+-$action\footnote{Here  
$\C_+$ stands for the additive group  of the complex number field.} 
$\vi_{\partial}$ on spec$\,A.$ For $a\in A$ and $t \in \C_+$ we put
$$\varphi_\partial(t,\,a)=\exp (t\partial)(a)=
\sum\limits_{i=0}^{n(\partial,\,a)-1} \frac{t^i\partial^ia}{i!} $$ 
\cite{Re}. The kernel ker$\,\partial$ coincides with the subalgebra 
$A^{\vi_{\partial}} \subset A$ of $\vi_{\partial}-$invariants. 

\label{df2.4}
\smallskip 

\noindent \ncor {{\bf Lifting LND-s and $\C_+-$actions to an affine modification.}}
{\it Let $(A,\,I,\,f)$ be an affine triple, and set 
$X =$spec$\,A.$
Let $\partial$ be an LND of the algebra $A = \C[X]$ 
such that $\partial(f) = 0$ and $\partial(I) \subset I.$ 
Then  $\partial$ can be
lifted in a unique way to an LND ${\partial}^{\prime}$ of the affine modification 
$A'= \C[X']= \Sigma_{I,\,f} (A).$ 
The corresponding $\C_+-$action $\vi_{\partial}$ on $X$ 
(which leaves invariant the subvarieties $D_f$ and $C=V(I)$) 
can be lifted in a unique way to a 
$\C_+-$action $\vi_{\partial}'$ on the affine modification 
$X'= \Sigma_{I,\,f}(X)$
which leaves invariant the exceptional divisor $E'.$}

\label{cor2.3}
\smallskip

\noindent {%\footnotesize
\ex \label{ex2.3} Let $A = \C[x,\,y],$ i.e. $X = \C^2,$ and $A' = \Sigma_{I,\,x}(A)$ where
$I = (x, \,q(y)),\,\,\,q \in \C[y].$ Then $X' = \Sigma_{I,\,x}(\C^2)$
is the surface given in $\C^3$ by the equation $xz-q(y) = 0,$
and the blowup morphism $\si_I\,:\,X \to \C^2$ is the restriction to $X$ 
of the standard projection 
$\C^3 \to \C^2,\,\,\,(x,\,y,\,z) \longmapsto (x,\,y)$ 
(see Example \ref{ex1.6}). 

Consider the triangular $\C_+-$action $\var(t,\,(x,\,y)) = (x,\,y+txp(x))$ 
on $\C^2$ where $p \in \C[x],$ and the corresponding LND 
$\partial\,:\,A \to A,$
$$
\partial\,:\,\,
%\begin{array}{l}%
%x \longmapsto 0\\[2pt]
%y \longmapsto xp(x)\,.
%\end{array}
(x,\,y) \longmapsto (0,\,xp(x))\,.
$$
Clearly, $\partial(I)\subset I.$ By Corollary \ref{cor2.3}, there exist
unique extensions $\partial'$ resp. $\var'$ of $\partial$ resp. $\var$
to $A';$ they are given, respectively, as follows:
$$
\partial'\,:\,\,
%\begin{array}{l}%
%x \longmapsto 0\\[2pt]
%y \longmapsto xp(x)\\[2pt]
%z = q(y)/x \longmapsto q'(y)p(x)\,
%\end{array}
(x,\,y,\,z= q(y)/x) \longmapsto (0,\,xp(x),\,q'(y)p(x))
$$
and
$$
\var'\,:\,\,
%\begin{array}{l}%
%x \longmapsto x\\[2pt]
%y \longmapsto y+txp(x)\\[2pt]
%z \longmapsto q(y+txp(x))/x = z + \sum\limits_{i=0}^{{\rm deg}\, q} 
%t^i{q^{(i)}(y)\over i!}x^{i-1}p^i(x)\,.
%\end{array}
(x,\,y,\,z) \longmapsto (x,\,y+txp(x),\,q(y+txp(x))/x 
= z + \sum\limits_{i=0}^{{\rm deg}\, q} 
t^i{q^{(i)}(y)\over i!}x^{i-1}p^i(x))\,.
$$}
\section{Topology of affine modifications}

%Hereafter $e(Y)$ denotes the Euler 
%characteristic of a topological space $Y;$ for a divisor $D,\,\,\,e(D)$ means
%$e(\mbox{supp}D).$ 
%We start with the next simple observation. 
%
%\smallskip
%
%\noindent \nlemma{{\bf The Euler 
%characteristic of an affine modification.}} 
%
%\noindent {\it Let $\si_I\,:\,X' \to X$ be the affine modification of 
%an affine variety $X$ along a divisor $D=D_f$ with center $I$ and with exceptional divisor 
%$E' \subset X'.$   
%Then $e(X') - e(X) = e(E') - e(D).$}
%
%\smallskip
%
%\noindent {\it Proof.~}  The isomorphism 
%$\si_I\,:\,X' \setminus E' \simeq X \setminus D_f$ 
%(see Proposition \ref{pr1.1}(a)) provides the equality
%$e(X' \setminus E') = e(X \setminus D_f).$ By
%the additivity of the Euler characteristic with respect to
%a disjoint constructive decomposition of a quasiprojective variety \cite{Du},
%we obtain $e(X') = e(X' \setminus E') + e(E')$ and 
%$e(X) = e(X \setminus D_f) + e(D_f),$ and the statement follows. \qed
%
%\label{lm3.1}
%\smallskip
%
%\noindent {%\footnotesize \ex Let $X = \C^n$ and $X' = \Sigma_{I,\,f}(X)$ 
%be as in Example \ref{ex1.6}. That is, 
%$I = (f,\,g)$ where $f,\,g \in \C^{[n]},$ and 
%$X' = \{f({\ol x})z - g({\ol x})=0\} \subset \C^{n+1}.$ 
%Put $D = \{f=0\} \subset \C^{n}$ and $C = V(I) = \{f=g=0\} \subset \C^{n}.$ 
%Then we have 
%$E' \simeq C \times \C,$ so that 
%$e(E') = e(C),$ and hence
%$$e(X') = 1 + e(C) - e(D).$$ In particular, $e(X') = 1$ iff $e(D) - e(C) = 
%e(D \setminus C) = 0.$ }
%
%\smallskip

%Till the end of 
In this section we appropriate the complex analytic point of view. 
Observe that, with the language of schemes, one can naturally extend 
the notion of the affine modification to quasiprojective varieties and 
to more general ring spaces and obtain results analogous to those of 
the previous sections. Instead, to simplify things, according
to Proposition \ref{pr1.1}(a) and Theorem \ref{thm1.1} we adopt
the following definitions.

\smallskip

\noindent \defin \label{def3.1} Consider a triple $(M,\,C,\,D)$ resp. 
a pair $(M',\,E')$ 
where $M$ resp. $M'$ is a reduced connected complex space, 
$D \subset M$ resp. $E' \subset M'$ is a  closed hypersurface,  
and $C \subset D$ is a closed analytic subvariety   
of codimension at least two in $M.$ Let $\si\,:\,M' \to M$ 
be a surjective
morphism  such that the restriction
$\si\,|\,(M' \setminus E')\,:\,M' \setminus E' \to M \setminus D$ 
is a biholomorphism, and $\si(E') \subset C.$ Then we say that the 
pair $(M',\,E')$ is a {\it pseudoaffine 
modification} (via $\si$) of the triple $(M,\,C,\,D)$ with a locus
subordinated to $(C,\,D)$ and with the exceptional divisor $E'.$ 

In particular, we consider the 
{\it pseudoaffine modification $M'=\Sigma_{C,\,D}(M)$ of $M$ 
along $D$  with center $C,$} where
$M' = {\widehat M} \setminus D', \,\,\,{\widehat M}$ 
is the blow up of $M$ at $C$ and $D'$ is the strict transform of 
$D$ in ${\widehat M}.$

{%\footnotesize
Anyhow, modifying $M$ we replace the divisor $D$ 
by the new one $E'.$ In the latter case $E' = E \setminus D'$ where 
$E=\si_C^{-1}(C) \subset {\widehat M}$ 
is the exceptional divisor of the blow up 
$\sigma_C\,:\,{\widehat M} \to M$ with center $C.$ 
While in the former case 
we blow up $M$ with center being an ideal sheaf supported by $C.$}

\smallskip

Recall that by a {\it vanishing loop} of a divisor $D$ in a complex manifold 
$M$ one means any loop in $\pi_1(M \setminus D)$ whose image 
in the homology group $H_1(M \setminus D)$ coincides with a fibre of 
the normal circle bundle of the smooth part reg$\,D.$

\smallskip

The next proposition 
provides a generalization of Lemma 3.4 in \cite{Ka 1}.

\smallskip

\noindent \nprop{{\bf Preserving the fundamental group under a modification.}}
{\it Let a pair $(M',\,E')$ be a pseudoaffine 
modification of a triple $(M,\,C,\,D)$ via a morphism 
$\si\,:\,M' \to M$ where $M,\,M'$ are complex manifolds, and let 
the divisors $D,\,E'$
admit finite decompositions into irreducible components
$D=\bigcup_{i=1}^n D_i$ resp. $E'=\bigcup_{j=1}^{n'} E'_j.$ 
Let $D_i^* = \sum_{j=1}^{n'} m_{ij} E'_{j}.$ 
Assume that 

%%\smallskip

\noindent (i) $\si(E_j') \cap {\rm reg}\,D_i \neq \emptyset$ as soon as 
$m_{ij} > 0.$

%%\smallskip

\noindent Then the following statements (a) and (b) hold.

%\smallskip

\noindent (a) If ($\alpha$) the lattice vectors 
$b_j = (m_{1j},\dots,m_{nj})\in \Z^n, \,\,j=1,\dots,n',$ 
generate the lattice $\Z^n,$
then ($\beta$) 
$\si_*\,:\,H_1(M';\,\Z) \to H_1(M;\,\Z)$ is an isomorphism.

%%\smallskip

\noindent The converse is true if  

%%\smallskip

\noindent (ii) $D_i$ is a principal divisor
defined by a holomorphic function $f_i$ on $M,$ i.e. 
$D_i=f_i^*(0),\,\,i=1,\dots,n.$ 

%\smallskip

\noindent (b) Assume further that

%%\smallskip

\noindent (iii) there is a disjoint partition 
$\{1,\dots,n'\} = J_1 \cup\dots\cup J_n$ 
such that 
$D_i^* = \sum_{j\in J_i} m_{ij} E'_{j} \neq 0,\,\,i=1,\dots,n.$

%%\smallskip

\noindent Set $d_i = $ g.c.d.$(m_{ij}\,|\,j\in J_i),\,\,i=1,\dots,n.$

%%\smallskip

\noindent Then ($\gamma$) 
$\si_*\,:\,\pi_1(M') \to \pi_1(M)$ is an isomorphism
if ($\delta$) $d_1=\dots=d_n=1.$

%%\smallskip

\noindent Under the condition (ii) the converse is also true.}

\label{pr3.1}
\smallskip

\noindent {\it Proof.~}   
Let $\alpha_{i}$ be a vanishing loop of $D_i$ in $M,\,\,i=1,\dots,n,$ and 
$\beta_{j}$ be a vanishing loop of $E'_{j}$ 
in $M',\,\,j=1,\dots,n'.$ Then the kernel of the natural homomorphism
$\pi_1(M\setminus D) \to \pi_1(M)$ coincides with 
the minimal normal subgroup $H :=\, <<\alpha_1,\,\dots,\alpha_n>>\,$ 
of the group $G:=\pi_1(M\setminus D)$ generated by 
$\alpha_1,\,\dots,\alpha_n,$ i.e. with the subgroup generated by
the conjugacy classes of $\alpha_1,\,\dots,\alpha_n$ 
(see e.g.  \cite[(2.3.a)]{Za 2}). Similarly, the kernel of the natural homomorphism
$\pi_1(M'\setminus E') \to \pi_1(M')$ is the normal subgroup
$H' = \,<<\beta_{1},\dots,\beta_{n'}>>\,$ of the group
$ \pi_1(M'\setminus E') \stackrel{\si_*}{\simeq} \pi_1(M\setminus D)=G.$ 
Thus, $\pi_1(M)\simeq G/H$ and $\pi_1(M')\simeq G/H'.$ Moreover, $\si_*$
identifies $H'$ with a subgroup of $H,$ the surjection 
$\si_*\,:\,\pi_1(M') \longrightarrow\!\!\!\to \pi_1(M)$ 
coincides with the canonical surjection 
$G/H \longrightarrow\!\!\!\to G/H',$ and hence, ker$\,\si_* \simeq H/H'.$ 
It follows that 
$\si_*\,:\,\pi_1(M') \to \pi_1(M)$ is an isomorphism iff $H = H'.$ 

%%\smallskip

\noindent {\it Proof of (a).} Denote $G' = [G,\,G],$ and let
$\rho\,:\,G \longrightarrow\!\!\!\to G/G' \simeq H_1(M \setminus D;\,\Z)$
be the canonical surjection. Set 
${\widetilde H} = \rho(H),\,\,\,{\widetilde H'} = \rho(H')$ and 
${\widetilde \alpha_i} = \rho(\alpha_i),\,\,\,{\widetilde \beta_j} = 
\rho(\beta_j).$ Clearly, ker($\,\si_*\,:\,H_1(M';\,\Z) \to H_1(M;\,\Z)$) 
$\simeq {\widetilde H}/{\widetilde H'}.$ 
In view of the condition ($i$) we have 
${\widetilde \beta_j} = \sum_{i=1}^n m_{ij} {\widetilde \alpha_i},\,\,j=1,\dots,n'.$ 
Under the condition ($\alpha$), the elements
${\widetilde \alpha_1},\dots,{\widetilde \alpha_n}$ can be expressed in 
terms of ${\widetilde \beta_1},\dots,{\widetilde \beta_{n'}},$ and so, 
${\widetilde H'}={\widetilde H}.$ This proves the implication 
($\alpha$) $\Longrightarrow$ ($\beta$). 

Assuming the condition ($ii$) we have an isomorphism
$f_*\,:\,{\widetilde H} \stackrel{\simeq}{\longrightarrow} \Z^n$ 
where 
$f :=(f_1,\dots,f_n)$
which sends 
${\widetilde \alpha_1},\dots,{\widetilde \alpha_n}$ 
to the standard generators of the lattice $\Z^n$ and 
identifies the subgroup ${\widetilde H'}$ with the sublattice 
$\Lambda \subset \Z^n$ spanned by
the vectors $b_j=(f\circ \si)_*({\widetilde \beta_j}),\,\,j=1,\dots,n'.$
Therefore,  
${\widetilde H}={\widetilde H'}$ iff $\Lambda = \Z^n,$ that is, iff the
condition ($\alpha$) holds. This yields the converse implication 
($\beta$) $\Longrightarrow$ ($\alpha$).

%%\smallskip

\noindent {\it Proof of (b).}   Under the assumptions ($i$) and ($iii$) 
the image $\si_*(\beta_J) \in G$ where $j \in J_i$ 
is conjugate with the element 
$\alpha_i^{m_{ij}}.$ Thus, 
$H' = <<\beta_1,\dots,\beta_{n'}>>$ 

\noindent $= <<\alpha_i^{m_{ij}}\,|\,j\in J_i,\,\,i=1,\dots,n>> = 
<<\alpha_i^{d_i}\,|\,i=1,\dots,n>>.$ 

\noindent Therefore, the condition ($\delta$) 
implies the coincidence 
$H=H',$ and so, it implies ($\gamma$). If ($i$)-($iii$) hold, then, clearly, 
we have the implications 
($\gamma$) $\Longrightarrow$ ($\beta$) $\Longrightarrow$ 
($\alpha$) $\Longleftrightarrow$ ($\delta$), 
and hence, ($\gamma$) $\Longrightarrow$ ($\delta$). 
\qed

\smallskip

The next theorem and its corollary generalize Theorem 3.5 in 
\cite{Ka 1} (see also \cite[Thm. 5.1]{Za 2}).

\smallskip

\noindent \nthm{{\bf Preserving the homology  under a modification.}}

\noindent {\it Let a pair $(M',\,E')$ be a pseudoaffine 
modification of a triple $(M,\,C,\,D)$ via a morphism 
$\si\,:\,M' \to M$ (see Definition \ref{def3.1}). Suppose that

%%\smallskip

\noindent (i) $M,\,M'$ are complex manifolds and $D,\,E'$ 
are topological manifolds admitting finite 
decompositions into irreducible components
$D=\sum_{i=1}^n D_i$ resp. 
$E' = \sum_{j=1}^{m} E'_{j}$ where $m=n$ and 
$E_i' = \si^*(D_i),\,\,i=1,\dots,n;$

%\smallskip

\noindent (ii) $\si(E_i) \cap $ reg$\,D_i \neq \emptyset,\,\,i=1,\dots,n.$ 

%%\smallskip

\noindent 
Put $\tau = i\circ(\si\,|\,E')\,:\,E' \to D$ where $i\,:\,C 
\hookrightarrow D$ is the identical embedding. 
Then 
$\si_*\,:\,H_*(M';\,\Z) \longrightarrow H_*(M;\,\Z)$
is an isomorphism  iff 

%%\smallskip

\noindent (iii) 
$\tau_*\,:\,H_*(E_i';\,\Z) \longrightarrow H_*(D_i;\,\Z)$ is an isomorphism
for all $i=1,\dots,n.$}

\label{thm3.1} 
\smallskip

\noindent {\it Proof.~}  Set $\breve{M} = M \setminus D,\,\,\,\breve{M'} = 
M' \setminus E'$ and 
$\breve{\si} = \si\,|\,\breve{M'}\,:\,\breve{M'}
\stackrel{\simeq}\longrightarrow \breve{M}.$ 
Consider the following commutative diagram where the horizontal 
lines are exact homology sequences of pairs with $\Z-$coefficients:

{\begin{picture}(500,95)

\put(-20,65){$\dots\longrightarrow H_{j+1}(M',\,\breve{M'})
\longrightarrow H_{j}(\breve{M'})\longrightarrow H_{j}(M')
\longrightarrow H_{j}(M',\,\breve{M'}) \longrightarrow 
H_{j-1}(\breve{M'})\longrightarrow \dots$}

\put(61,50){$\vector(0,-1){25}$}
\put(66,35){$(\si,\,\breve{\si})_*$}
\put(137,50){$\vector(0,-1){25}$}
\put(142,35){$\breve{\si}_*$}
\put(203,50){$\vector(0,-1){25}$}
\put(208,35){$\si_*$}
\put(277,50){$\vector(0,-1){25}$}
\put(282,35){$({\si},\,\breve{\si})_*$}
\put(350,50){$\vector(0,-1){25}$}
\put(355,35){$\breve{\si}_*$}
\put(420,35){$(*)$}

\put(-17,5){$\dots\longrightarrow H_{j+1}(M,\,\breve{M})
\longrightarrow \,\,H_{j}(\breve{M})\,\longrightarrow H_{j}(M)
\longrightarrow \,H_{j}(M,\,\breve{M}) \,
\longrightarrow \,\,H_{j-1}(\breve{M})\,\,\longrightarrow \dots$}

\end{picture}

\noindent Due to the Thom isomorphism, 
it can be replaced by the following one:

\begin{picture}(200,95)

\put(0,65){$\dots\longrightarrow H_{j-1}(E')
\longrightarrow H_{j}(\breve{M'})\longrightarrow H_{j}(M')
\longrightarrow H_{j-2}(E') \longrightarrow 
H_{j-1}(\breve{M'})\longrightarrow \dots$}

\put(63,50){$\vector(0,-1){25}$}
\put(69,35){$\tau_*$}
\put(137,50){$\vector(0,-1){25}$}
\put(122,35){$\simeq$}
\put(142,35){$\breve{\si}_*$}
\put(199,50){$\vector(0,-1){25}$}
\put(204,35){$\si_*$}
\put(270,50){$\vector(0,-1){25}$}
\put(275,35){$\tau_*$}
\put(342,50){$\vector(0,-1){25}$}
\put(327,35){$\simeq$}
\put(347,35){$\breve{\si}_*$}
\put(420,35){$(**)$}

\put(9,5){$\dots\longrightarrow H_{j-1}(D)
\longrightarrow H_{j}(\breve{M})\longrightarrow H_{j}(M)
\longrightarrow H_{j-2}(D) 
\longrightarrow H_{j-1}(\breve{M})\longrightarrow \dots$}

\end{picture}

\noindent which is still commutative. Indeed, let $\Delta_i'$ 
be a small complex disc in $M'$ which meets $E_i'$ transversally 
at its origin, 
and let this origin be a generic point of $E_i'.$ 
Then the image
$\Delta_i = \si(\Delta_i')$ is a complex disc in $M$ transversal to reg$\,D_i$ 
at a generic point of $C_i.$ Hence, for the Thom classes 
$u_i' \in H^{2}(M',\,M'\setminus E'_i)$ of $E'_i$ resp. 
$u_i \in H^{2}(M,\,M \setminus D_i)$ of $D_i,$ uniquely defined by the conditions 
$u_i'(\Delta_i') = 1$
resp.  $u_i(\Delta_i) = 1,$ we have $\si^*(u_i) = u_i',\,i=1,\dots,n.$ 
Recall that the 
Thom isomorphism $H_{j+2}(M',\,\breve{M'}) \simeq H_{j}(E')$ resp. 
$H_{j+2}(M,\,\breve{M}) \simeq H_{j}(D')$ is defined as the cap-product 
with the 
Thom class $u' = \sum_{i=1}^n u_i'$ resp. $u = \sum_{i=1}^n u_i\,:$ 
$\,\,\eta' \mapsto u'\cap \eta'$ resp. $\eta \mapsto u\cap \eta$
(see e.g. \cite[VIII.11.21]{Do}, \cite{MilSta}). 
Now it is easily seen that the diagram $(**)$ is commutative. 

It remains to note that by the Five Lemma, $\tau_*$ in $(**)$ 
yields an isomorphism for all $j$
iff  $\si_*\,:\,H_{j}(M') \to H_{j}(M)$ is an isomorphism for all $j.$
The proof is completed. \qed

\smallskip

{%\footnotesize 
\noindent \rems 

\noindent \nrem \label{rm3.1}
If
$\si=\si_C$ is 
the blowup morphism with a 
smooth center $C$ contained in reg$\,D,$ or otherwise, if 
$M'$ is the affine modification $\Sigma_{I,\,f}(M)$ 
with the ideal $I$ admitting 
a regular system of generators started with $f$ (see Proposition \ref{pr1.1}(d)),
then 
$\tau\,:\,E' \to C$ is a smooth fibration over $C$ with a fibre $\C^k,$ 
where $k =$codim$_D C,$ and so, the contraction $\si\,|\,E'\,:\,E' \to C$
is a homotopy equivalence. Thus, in this case
$\tau_*\,:\,H_*(E';\,\Z) \longrightarrow H_*(D;\,\Z)$ 
is an isomorphism iff 
$i_*\,:\,H_*(C;\,\Z) \longrightarrow H_*(D;\,\Z)$ is.

\smallskip

\noindent \nrem \label{rm3.1'} From the biholomorphism 
$\si\,:\,M' \setminus E' \to M 
\setminus D$ we get the equality for the Euler characteristics 
$e(M') - e(M) = e(E') - e(D)$ (see \cite{Du}). 
In the situation as in the above remark we have $e(E') = e(C),$ and hence 
$e(M') - e(M) = e(C) - e(D).$}

\smallskip

\noindent \ncor{{\bf Preserving contractibility under a modification.}}

\noindent {\it Suppose that the conditions ($i$), ($ii$) and ($iii$) 
of Theorem \ref{thm3.1} are fulfilled. Then
the pseudoaffine modification $M'$ of $M$ is acyclic 
resp. contractible iff $M$ is.}

\label{cor3.2}
\smallskip

\noindent {\it Proof.~}  The equivalence of acyclicity of $M'$ and of $M$ 
follows immediately from Theorem \ref{thm3.1}. 

By the Hurewicz and Whitehead Theorems \cite[(2.11.5), (2.14.2)]{FoFu},  
contractibility of $M$ resp. $M'$ is equivalent to
acyclicity and simply connectedness of $M$ resp. $M'.$
Notice that  
the conditions ($i$), ($iii$) and ($\delta$) of Proposition \ref{pr3.1} are
fulfilled. By (b) of this proposition, $\si\,:\,M' \to M$
induces an isomorphism of the fundamental groups. Thus, $M'$ is 
simply connected iff $M$ is, and the statement follows. \qed 

\smallskip

%\noindent \rem \label{rm3.2} Suppose that $\pi_1(M) = {\bf 1}$ or $\pi_1(M') = {\bf 1},$ 
%the homomorphism
%$\si_*\,:\,\pi_2(M') \to \pi_2(M)$ is surjective and the conditions
%($i$) - ($iii$) of Theorem \ref{thm3.1} are fulfilled. Then by the Whitehead
%Theorem \cite[(2.14.5)]{FoFu},  
%$\si\,:\,M' \to M$ is a homotopy equivalence. 

\smallskip

{%\footnotesize
We give below several examples of applications of Corollary \ref{cor3.2}.

\smallskip 

\noindent \nex{{\bf The tom Dieck--Petrie surfaces} \cite{tDP}.} \label{ex3.1} 
These are the smooth surfaces $X_{k,\,l}$ in $\C^3$
defined by the polynomials
$$p_{k,l}= {(xz + 1)^k - (yz + 1)^l - z \over z} \in \C[x,\,y,\,z]\,$$
where $k,\,l \ge 2,\,\,\,$gcd$(k,\,l)=1.$ 
As in Example \ref{ex2.2} above (see also \cite[Example 6.1]{Za 2}) 
we have $X_{k,\,l} = \Sigma_{\rho,\,\Gamma_{k,\,l}} (\C^2)$
where $\Gamma_{k,\,l} = \{x^k - y^l = 0\} \subset \C^2$ and 
$\rho = (1,\,1) \in \Gamma_{k,\,l}.$  
By Corollary \ref{cor3.2}, the surfaces $X_{k,l}$ are contractible. 

\smallskip

\noindent \nex{{\bf Russell's cubic is contractible.}} \label{ex3.2} 
Recall (see Example \ref{ex1.7}) that the Russell cubic threefold 
$X=\{x + x^2y + z^2 + t^3=0\}\subset\C^4$ 
is the affine modification of $\C^3$ along the divisor $2D,$
where $D := \{x=0\} \subset \C^3,$ with center at the ideal 
$I = (-x^2, \,x + z^2 + t^3) \subset \C^{[3]}$ 
supported by the plane curve
$C = \Gamma_{2, \,3} =\{x=z^{2}+t^{3}=0\}\subset D \simeq \C^2.$ 
The exceptional divisor $E'$ coincides with the book-surface 
$B=\{x=0\} \subset X',\,\,\,B \simeq \C\times \Gamma_{2, \,3}.$
Therefore, the condition ($iii$) of Theorem \ref{thm3.1} is fulfilled. 
As well, the conditions ($i$) and ($ii$) hold, 
and so, by Corollary \ref{cor3.2}, $X$ is contractible.  Moreover, by 
the Dimca-Ramanujam Theorem \cite{Di 1, Ram} (see also \cite[Thm 4.2]{Za 2}),
$X$ is diffeomorphic to $\R^6.$

\smallskip 

\noindent \nex {{\bf Modifying an affine space along a scroll.}}
\label{ex3.3} 
Let $M = f^*(0)$ 
be a smooth reduced analytic hypersurface in $\C^n$ 
where $f \in {\cal O}(\C^n),$ 
and let 
$f_1,\dots,f_k \in {\cal O}(\C^n)$ be holomorphic functions without 
common zeros on $M.$ 
Consider the smooth analytic subvarieties
$$D := M \times \C^k \subset \C^{n+k}_{(\bx,\,\bu)}=\C^n_{\bx} 
\times \C^k_{\bu}\,\,\,\,\,\,\,\,
{\rm  and}\,\,\,\,\,\,\,
C := \{f(\bx) = 0 = g(\bx,\,\bu)\} \subset \C^{n+k}_{(\bx,\,\bu)}$$ 
where $g(\bx,\,\bu) := \sum_{i=1}^k u_if_i \in {\cal O}(\C^{n+k}).$ 
The natural embeddings 
$M \times {\overline 0} \hookrightarrow C \hookrightarrow D$ being 
homotopy equivalences, by Corollary \ref{cor3.2} and Remark \ref{rm3.1}, 
the pseudoaffine modification $X' = \Sigma_{C,\,D} (\C^{n+k})$ of $\C^{n+k}$
is a smooth contractible analytic hypersurface given in 
$\C^{n+k+1}_{(\bx,\,\bu,\,v)}$ by the equation 
$f(\bx)v - g(\bx,\,\bu) = 0$ (cf. Example \ref{ex1.6}). 

Since $M$ is supposed being smooth one may take, for instance,  
$k=n$ and $f_i = {\partial f \over \partial x_i}, \,\,i=1,\dots,n;$
then $C \simeq TM$ and $D \simeq T\C^n\,|\,M.$ 

\smallskip

The next example shows that in general, Corollary \ref{cor3.2} does not hold 
if some of the conditions ($i$) -- ($iii$) of Theorem \ref{thm3.1} 
are violated. 

\smallskip 

\noindent \nex{{\bf Modifying with center at the singular locus 
of the divisor} \cite[Remark on p. 418]{Ka 1}.} \label{ex3.4}
Let 
$\si\,:\,X \to  \C^2$ be a pseudoaffine modification of $\C^2$ 
along the cuspidal cubic $\Gamma:=\Gamma_{2,3} = \{x^2 - y^3 = 0\}$ 
with center at the cusp ${\overline 0} \in \Gamma,$ that is, $\si$ 
is the restriction to $X := \hX \setminus  \Gamma'$ of the blowup morphism  
$\hsi\,:\,\hX \to \C^2$ with center at the origin, whereas 
$\Gamma' \subset \hX$ is the strict transform of $\Gamma$ in 
$\hX.$ Then $X$ is neither simply connected nor acyclic. Indeed, 
set $X^* =  X \setminus E' = \hX \setminus (E \cup \Gamma')$ where 
$E\subset \hX$ is the exceptional divisor of $\hsi,$ and let  
$\beta \in H_1(X^*;\,\Z)$ be a vanishing loop of $\Gamma'$ in $X^*.$ 
It is easily seen that $\si_*(\beta) = 2\alpha$ where  
$\alpha \in H_1(\C^2 \setminus \Gamma;\,\Z) \simeq \Z$ 
is a vanishing loop of 
$\Gamma$ in $\C^2.$ 
In virtue of the isomorphism 
$\si\,|\,X^*\,:\,X^* \stackrel{\simeq}{\longrightarrow} 
\C^2 \setminus \Gamma$
and of the exact sequence 
$${\bf 0} \longrightarrow \,<\beta> \,\longrightarrow H_1(X^*;\,\Z) 
\longrightarrow H_1(X;\,\Z) \longrightarrow {\bf 0}$$
(cf. Lemma \ref{lm4.3} below)
we have $$H_1(X;\,\Z) \simeq H_1(X^*;\,\Z) / \,<\beta> \, \simeq
H_1(\C^2 \setminus \Gamma;\,\Z) / \,<2\alpha> \, \simeq \Z / 2\Z.$$

\noindent \rem \label{rm3.3} 
Notice that the proper transform $\Gamma^{\rm pr}$ of $\Gamma$ in $\hX$
(see Remark \ref{rm1.3'}) coincides with the union 
$E \cup \Gamma',$ and so, it differs from the strict transform $\Gamma'.$ 
In turn, $X$ differs from the affine modification 
$\Sigma_{I,\,f}(\C^2) \simeq \C^2 \setminus \Gamma$ where 
$I := (x^2-y^3,\,x,\,y).$
Indeed, the blow up $Z = \hX$ can be given in 
$\C^2_{(x,\,y)} \times \pr^2_{(u : v : w)}$ by the equations 
$xv - yu =0,\,\, w = xu-y^2v,$ and hence, $\Gamma^{\rm pr} = \{w=0\}$ 
contains the exceptional divisor\footnote{More generally, 
one can prove the following:

{\it If the center $C$ of a blow up of a manifold $X$
is smooth, reduced and it is contained in the singular 
locus of a reduced divisor $D_f,$ 
then the proper transform $D_f^{\rm pr}$ contains
the exceptional divisor $E$ of the blow up.}} 
$E = \{x=y=0\}.$}

\section{Topology of the hypersurfaces $uv=p(x_1,\dots,x_k)$}

\noindent {\bf Notation.} Let $X = X(p) \subset \C^{k+2}$ be the irreducible
hypersurface given by the equation  $uv=p(\bx)$ where 
$\bx = (x_1,\dots,x_k) \in \C^k,\,\,k\ge 1,$
and $p \in \C^{[k]},$ deg$\,p > 0.$  For $c \in \C$ denote 
$U_c = X \cap \{u = c\},\,\,V_c = X \cap \{v = c\},$ and $X_0 = p^{-1}(0)
\subset \C^k.$
We also regard $X_0$ as the subvariety of $X$ given in $\C^{k+2}$ by the 
equations $u=v=p(\bx)=0.$

{%\footnotesize
\smallskip

\noindent \rem \label{rm5.1}
The variety $X(p)$ is the affine modification of 
$\C^{k+1}_{\bx,\,u}$ along the hyperplane $D_u = \{u=0\}$ with center 
$I = (p, u) \subset \C^{[k+1]}$ and with the exceptional divisor 
$U_0 \subset X(p)$
(see Example \ref{ex1.6}). If the divisor $p^*(0)$ is reduced, then 
$X(p) =\Sigma_{X_0,\,D_u} (\C^{k+1}).$}

\smallskip

The following lemma is easy, and we omit the proof. 

\smallskip

\noindent \lemma \label{lm5.1}
{\it (a) There are natural isomorphisms 
$X \setminus U_0 \simeq X \setminus V_0
\simeq \C^k \times \C^*,\,\,\,U_0 \simeq V_0 \simeq X_0 \times \C,$ 
and $U_0 \cap V_0 = X_0.$  In particular, $e(X) = e(X_0).$

%%\smallskip

\noindent (b) $X$ is smooth iff 
$p^*(0)$ is a smooth reduced (not necessarily irreducible) divisor in $\C^k$
(supported by $X_0$).
Furthermore, if the divisor $p^*(0)$ is reduced, then  we have
$${\rm sing}\,X ={\rm sing}\,U_0 \cap {\rm sing}\,V_0 = 
{\rm sing}\,X_0\quad \mbox{and}\quad  
{\rm sing}\,U_0 \simeq ({\rm sing}\,X_0) \times \C \simeq {\rm sing}\,V_0$$}
\prop \label{pr5.1}
{\it Assume that $X$ is smooth. Then $X$ is simply connected, and 
the embedding $X_0 \hookrightarrow X$ induces an isomorphism of 
the reduced homology groups $\tH_*(X;\,\Z) \simeq \tH_{*-2}(X_0;\,\Z).$}

\smallskip 

\noindent {\it Proof.~}  Since $X$ is assumed being smooth, 
by Lemma \ref{lm5.1}(b),
$X_0$ is smooth and reduced. By Lemma \ref{lm5.1}(a), we have that 
$\pi_1(X \setminus U_0) \simeq \pi_1(\C^k \times \C^*)\simeq \Z.$ 
The kernel of the natural surjection $i_*\,:\, \pi_1(X \setminus U_0) 
\longrightarrow\!\!\!\to \pi_1(X)$ 
is generated,
as a normal subgroup, by vanishing loops, say $\alpha_i,$ of the irreducible 
components $U^{(i)}_0 \simeq X^{(i)}_0 \times \C$ of $U_0$ where 
$X^{(i)}_0 \subset \C^k,\,\,i=1,\dots,l,$ are the irreducible components of $X_0$
(see e.g. \cite[(2.3.a)]{Za 2}). 

Denote $\rho = (u\,|\,(X \setminus U_0))_*\,:\,\pi_1(X \setminus U_0) \to 
\pi_1(\C^*) \simeq \Z.$ Fix a generator $t$ of the group 
$\pi_1(X \setminus U_0) \simeq \Z$ such that $\rho(t) = 1.$ Clearly,   
$\rho(\alpha_i) = \pm 1 \in \Z = \pi_1(\C^*);$ that is, 
$\alpha_i = \pm t,\,\,i=1,\dots,l,$ and hence, Ker$i_* = 
\langle \alpha_1,\dots,\alpha_l 
\rangle = \Z.$
Thus, $\pi_1(X) = {\bf 1},$ as claimed.

Therefore, also $H_1(X;\,\Z) = 0.$ Due to Lemma \ref{lm5.1}(a), from the exact
homology sequence of the pair $(X,\,X \setminus U_0)$ we obtain the 
isomorphisms of the $\Z-$homology groups: 
$H_i(X) \simeq H_i(X,\,X \setminus U_0)$ for all $i \ge 3.$ 
In view of 
the Thom isomorphism (denoted below by $\tau$) this 
yields  the isomorphisms
$$H_i(X) \simeq H_i(X,\,X \setminus U_0)\stackrel{\tau}{\simeq} 
H_{i-2}(U_0)\simeq H_{i-2}(X_0)
\,\,\,\forall \,\,\,i\ge 3\,.$$
For $i=2$ we get
$$0 = H_2(X \setminus U_0) \to H_2(X) \to H_2(X,\,X \setminus U_0)
\stackrel{\tau}{\simeq} H_{0}(U_0)$$ 
$$\simeq H_{0}(X_0) \stackrel{\partial_*}{\longrightarrow} \Z  
\simeq H_1(X \setminus U_0) \to 0 = H_1(X)\,.$$
Since $X_0$ is a smooth reduced divisor, the number $l= $deg$\,p$ of 
its irreducible components coincides with the number $b_0(X_0)$ of 
its connected components. Recall (see e.g. \cite{MilSta}) that the 
Thom isomorphism 
$\tau\,:\,H_0(U_0) \stackrel{\simeq}{\longrightarrow} 
H_2(X,\,X \setminus U_0)$ 
to each point $P \in U_0$ associates the relative homology class 
$[\Delta_P] \in H_2(X,\,X \setminus U_0)$ of a disc
$\Delta_P \subset X$ centered at $P$ and transversal 
to $U_0;$  furthermore, $\partial_*[\Delta_P] = [\partial \Delta_P] 
\in H_1(X \setminus U_0).$ 
Thus, for  
$\beta = (\beta_1,\dots,\beta_l) \in H_2(X,\,X \setminus U_0) 
\simeq H_0(X_0) \simeq \Z^l$ 
we have 
$\partial_*(\beta) = \sum_{i=1}^l \beta_i \in H_1(X \setminus U_0) \simeq \Z.$
It follows that 
$H_2(X) \simeq $ ker$\,\partial_* \simeq \Z^{b_0(X_0)-1} \simeq \tH_{0}(X_0);$ 
in particular, $H_2(X) = 0$ iff $X_0$ is irreducible. 
Since $X$ is irreducible, we conclude that 
$\tH_*(X) \simeq \tH_{*-2}(X_0),$ 
and the second statement of the proposition follows.
\qed

\smallskip

\noindent \cor {\it Suppose that $X$ is smooth. Then  
$X$ is contractible iff $X_0$ is acyclic; 
actually, in this case $X$ 
is diffeomorphic to $\R^{2k+2}.$}

\label{cor5.2}
\smallskip

The last statement follows from 
the Dimca-Ramanujam Theorem \cite{Di 1, Ram} 
(see also \cite[Thm 4.2]{Za 2}) for $k \ge 2,$ 
and for $k=1$ it follows from the next proposition.  

\smallskip

\noindent \prop {\it (a) Let $k=1.$ 
The following conditions 
are equivalent:

%%\smallskip

(i) $X \simeq \C^2;$

%%\smallskip

(ii) $X\subset \C^3$ is a smooth acyclic surface;

%%\smallskip

(iii) deg$\,p = 1.$

%\smallskip 

\noindent (b) Let $k=2.$ The following conditions are equivalent:

(i) the embedding $X \hookrightarrow \C^4$ is rectifiable;

%%\smallskip

(ii) $X \simeq \C^3;$

%%\smallskip

(iii) $X\subset \C^4$ is a smooth acyclic 3-fold;

%%\smallskip

(iv) the divisor $p^*(0)$ is reduced and $X_0 \simeq \C;$ 

%%\smallskip

(v) the polynomial $p \in \C^{[2]}$ is equivalent to a linear one,
i.e. for some $\alpha \in $Aut$\,\C^2$ we have $p\circ \alpha = q$ 
where $q(x_1,\,x_2) = x_1.$}

\label{pr5.2}
\smallskip 

\noindent {\it Proof.~}  (a) Let $k = 1.$ The implications ($iii$) 
$\Longrightarrow$ 
($i$) $\Longrightarrow$ ($ii$) are clear. Conversely, 
by Lemma \ref{lm5.1}(a),(b), under the assumption 
($ii$) we have $1= e(X) = e(X_0) =$ deg$\,p.$ This proves (a).

%%\smallskip 

\noindent (b) The theorem of Abhyankar-Moh and Suzuki
\cite{AM, Suz} asserts the equivalence ($iv$) $\Longleftrightarrow$ 
($v$). The implications ($v$) $\Longrightarrow$ 
($i$) $\Longrightarrow$ ($ii$) $\Longrightarrow$ ($iii$) are easy; 
($iii$) $\Longrightarrow$ ($iv$) 
follows from Lemma \ref{lm5.1}(b) and Proposition \ref{pr5.1}, 
and so, we are done.
\qed

\smallskip
{%\footnotesize
\noindent \rem \label{rm5.2}
However, starting with $k=3$ 
the equivalences of Proposition \ref{pr5.2} fail; 
see Examples \ref{ex7.1} and  \ref{ex7.2} below. }

\smallskip

\noindent {\bf A generalization.} More generally, consider the variety 
$X = X(p_1,\dots,p_m) \subset \C^{k+m+1}$ 
given by a system of equations $uv_i=p_i(\bx),\,\,i=1,\dots,m$ where 
$\bx = (x_1,\dots,x_k) \in \C^k,\,\,\,k \ge m \ge 1,$ and 
$p_i \in \C^{[k]} \setminus \C,\,\,i=1,\dots,m.$
This variety $X$ is the affine modification of $\C^{k+1}_{(\bx,\,u)}$ 
along the hyperplane $D_u = \{u=0\}$ with center at the ideal
$I = (p_1,\dots,p_m,\, u) \subset \C^{[k+1]}$ supported 
by the affine variety
$X_0 := \{p_1(\bx) =\dots =p_m(\bx) =0\} \subset D_u \simeq \C^{k}.$ 
Clearly, $X$ is irreducible iff $X_0$ is a set-theoretic complete 
intersection, i.e. all its irreducible components are codimension $m$ 
subvarieties of $\C^{k},$ which will be always assumed in the sequel. 
Under this assumption most of the results proved above 
for $m = 1$ remain true
in this more general setting. 

\smallskip

\noindent \prop {\it (a) Denote by $U_0$ the hypersurface in 
$X= X(p_1,\dots,p_m)$ given by the equation $u=0.$ Then we have 
$U_0 \simeq X_0 \times \C^m$ and 
$X \setminus U_0 \simeq \C^* \times \C^m.$ In particular, 
$e(X) = e(X_0).$

%%\smallskip

\noindent (b) Assume further that 

%%\smallskip

\noindent (i)
$X_0$ is the ideal-theoretic complete intersection of the divisors 
$p_i^*(0),\,\,i=1,\dots,m,$ that is, $I(X_0) = (p_1,\dots,p_m).$ 

%%\smallskip

\noindent 
Then we have 
${\rm sing}\,X_0 \subset  {\rm sing}\,X \subset {\rm sing}\,U_0 \simeq 
({\rm sing}\,X_0) \times \C^m.$
In particular, $X$ is smooth iff $X_0$ is.

%%\smallskip

\noindent (c) Suppose that the assumption (i) is fulfilled and $X$ is smooth.
Then we have $\pi_1(X) = {\bf 1}$ and 
$\tH_*(X;\,\Z) \simeq \tH_{*-2}(X_0;\,\Z).$
In particular, $X$ is contractible (and, moreover, diffeomorphic to 
$\R^{2k +2}$ for $k \ge 2$) iff $X_0$ is acyclic.}

\label{pr5.3}
\smallskip 

The proof goes exactly in the same way as before, and so, 
we leave it to the reader.

{%\footnotesize
\smallskip 

\noindent \rem \label{rm5.3} Presumably, even being diffeomorphic to the 
affine space, the variety $X= X(p_1,\dots,p_m)$ is not, in general, 
isomorphic to $\C^{k+1}$ for $k \ge 3,$ and so, 
it should provide an exotic algebraic structure on $\C^{k +1}$ 
(see e.g. \cite{Za 2}). But at present we have no invariant available
to distinguish $X$ from the affine space (see Remark \ref{rm6.1} below). 

Suppose that, indeed, for a certain smooth acyclic 
complete intersection $X_0 \subset \C^k$ the variety $X= X(p_1,\dots,p_m)$
is not isomorphic to $\C^{k+1}.$ Then we would have an example 
showing that Miyanishi's characterization of the affine 3-space ${\bf A}^3_k$ \cite{Miy}
does not hold any more in higher dimensions. 
Indeed, by Proposition \ref{pr5.3}, the varieties 
$X$ and $U_0 \simeq X_0 \times \C^m$ being smooth and acyclic 
we have $e(X) = e(U_0)=1,\,\,\,X \setminus U_0 \simeq \C^* \times \C^k,$ 
the algebras $\C[X]$ and  $\C[U_0]$ are UFD  
and have only constants as the units (see e.g. \cite[Prop. 3.2]{Ka 1}).
Thus, all the assumptions of the Miyanishi Theorem are fulfilled, whereas 
$X \not\simeq \C^{k+1}.$   

Further, if for a certain codimension $m$ smooth acyclic 
complete intersection 
$X_0 \subset \C^k$ non-isomorphic to $\C^{k-m}$  
the variety $X$ were isomorphic to $\C^{k+1},$
this would answer in negative, alternatively, either 
to the Zariski Cancellation Problem\footnote{In the particular case 
when $k=3$ and $X_0 \simeq \C$ 
is a complete intersection given by $p_1(\bx) = p_2(\bx)=0$ in $\C^3,$ 
the smooth contractible 4-folds $X=X(p_1,\,p_2)\subset \C^5$ 
were studied, in algebraic fashion, in \cite{As} as potential 
counterexamples to the Zariski Cancellation Problem. Indeed, in
\cite{As} an isomorphism $X \times \C \simeq \C^5$ was established.} or 
to the Abhyankar-Sathaye Embedding Problem\footnote{Cf. 
also Remark \ref{rm7.3}
below for another conjectural counterexample 
to the Abhyankar-Sathaye Embedding Problem.}. 
Indeed, the former happens if 
the hypersurface $U_0 = X_0 \times \C^m \subset X$ is isomorphic to 
$\C^k$ (which is only possible if $X_0$ was contractible). 
Otherwise, the latter takes place since 
$U_0$ is the zero fibre of the polynomial 
$u\,|\,X \in \C[X] \simeq \C^{[k+1]}$
with all other fibres $U_c,\,\,c \neq 0,$ isomorphic to $\C^k.$
Observe that due to the Miyanishi-Sugie and Fujita Cancellation Theorem, 
for $k-m \le 2$ only the second possibility might happen. 

Formally, there is also a possibility that $X \not\simeq \C^{k+1}$ whereas
$U_0 \simeq \C^k.$ In that case we would have an example of an exotic 
$\C^n,\,\, n = k+1 \ge 4,$ fibered by the affine spaces 
$U_c \simeq \C^{n-1}.$ }

\smallskip

We may enlarge our collection of contractible affine varieties 
passing to ramified cyclic coverings over $X= X(p_1,\dots,p_m),$ as follows. 

\smallskip

\noindent \prop {\it 
(a) Suppose that the variety $X = X(p_1,\dots,p_m)$ as in Proposition 
\ref{pr5.3}(c) above is smooth and contractible. Then for
any $n \in \N$ the variety $X_{n} \subset \C^{k+m+1}$ given by the 
system of equations $u^nv_i=p_i(\bx),\,\,i=1,\dots,m,$ 
is smooth and contractible, too.

%%\smallskip

\noindent (b) For a sequence of integers $s_0,\dots,s_m \in \N$ such 
that gcd$(s_i,\,s_j) = 1$ for all $i \neq j,$ consider the variety 
$$Y = Y_{s_0,\dots,s_m}(p_1,\dots,p_m):=\{u^{s_0}v_i^{s_i} = 
p_i(\bx),\,\,i=1,\dots,m\} \subset \C^{k+m+1}\,$$
where $p_i \in \C^{[k]},\,\,i=1,\dots,m.$
Suppose that the following conditions are fulfilled:

%%\smallskip

\noindent (i) $p := p_1\cdot\dots\cdot p_m$ is a prime decomposition, and
$D := p^*(\overline 0)$ is a reduced simple normal crossing divisor in
$\C^k;$

%%\smallskip

\noindent (ii) the divisor $D_i:=  D_{p_i}$ is $\Z_q-$acyclic for 
any prime divisor $q$ of $s_i,,\,\,i=1,\dots m,$ and
$X_0 = \cap_{i=1}^m D_i \subset \C^k$ 
is a smooth acyclic complete intersection; 

%%\smallskip

\noindent (iii) the group $\pi_1(\C^k \setminus D)$ is abelian 
(and hence, isomorphic to $\Z^m$).

%%\smallskip

\noindent Then $Y$ is a smooth contractible variety; 
it is diffeomorphic to $\R^{2k+2}$ if $k \ge 2.$}

\label{pr5.4}
\smallskip

\noindent {\it Proof.~}  $(a)$ The variety $X_n$ is a cyclic covering of $X$ 
ramified to order $n$ on $U_0$ with the covering morphism 
$\rho\,:\, X_n \to X,\,\,\,\rho\,:\,(\bx,\,u,\,\bv) \longmapsto (\bx,\,u^n,\,\bv).$  
By Proposition \ref{pr5.3}(c), the ramification divisor 
$U_0 = X_0 \times \C^m$ is acyclic,
the fundamental group of its complement $\pi_1(X \setminus U_0) \simeq \Z$ 
is abelian, and the regular function $u\,|\,(X \setminus U_0)$ where
$X \setminus U_0 \simeq \C^* \times \C^k,$ is a quasi-invariant of weight $1$
of the natural $\C^*-$action. Now the assertion of $(a)$ 
follows from Theorem A in \cite{Ka 3} (see also \cite[Thm. 7.1]{Za 2}).

%%\smallskip

\noindent $(b)$ Denote $V^{(i)}_0 = X \cap \{v_i = 0\},\,\,i=1,\dots,m$ where
$X = X(p_1,\dots,p_m).$
The morphism $\C^{k+m+1} \to \C^{k+m+1},\,\,\,(\bx,\,u,\,v_1,\dots,v_m) 
\longmapsto (\bx,\,u^{s_0},\,v_1^{s_1},\dots,v_m^{s_m}),$ restricted to $Y$
makes $Y$ a multi-cyclic covering of $X$ branched to order 
$s_0$ over $U_0$ resp. to order $s_i$ over $V^{(i)}_0,\,\,i=1,\dots,m.$
Thus, we may use Theorem 8.1 in \cite{Za 2} which provides conditions 
to guarantee contractibility of a multi-cyclic covering over a 
contractible manifold. To see that these conditions are satisfied, 
first of all, we observe
that the function $u$ resp. $v_i,\,\,i=1,\dots,m,$ is a quasi-invariant 
of weight $1$ resp. $-1$ of the $\C^*-$action 
$(\lambda,\,(\bx,\,u,\,\bv)) \longmapsto (\bx,\,\lambda u,\,\lambda^{-1} \bv)$ 
on $\C^{k+m+1}$ which leaves the variety $X$ invariant. 

Further, the hypersurface $U_0 \simeq X_0 \times \C^m$ is smooth and acyclic 
since $X_0$ is. For each $i=1,\dots,m$ the hypersurface $V^{(i)}_0$ is the 
affine modification of the smooth variety $D_i \times \C$ along the divisor 
$D_i$ with center $X_0 \subset D_i$ and with the exceptional divisor 
$E_i := U_0 \cap V^{(i)}_0 \simeq X_0 \times \C^{m-1}$ via the morphism 
$\si_i\,:\,V^{(i)}_0 \to D_i \times \C$ which is the restriction to 
$V^{(i)}_0$ of the projection 
$\pi\,:\,\C^{k+m+1} \to \C^{k+1},\,\,\,\pi(\bx,\,u,\,\bv)=(\bx,\,u)$ 
(see Corollary \ref{cor2.1}). 

Notice that the proof of Theorem \ref{thm3.1} on preservation of 
the homology under a modification goes equally for the $\Z_p-$homology groups.
In virtue of the condition $(ii)$ above, by this Theorem,   
the smooth hypersurface $V^{(i)}_0 \subset X$ is 
$\Z_q-$acyclic for any prime divisor $q$ of $s_i,\,\,i=1,\dots,m.$

We have an isomorphism
$X^* := X \setminus (U_0 \cup V^{(i)}_0 \cup \dots \cup V^{(i)}_0) \simeq
(\C^k \setminus D) \times \C^*.$ Hence by the condition $(iii),$
the fundamental group $\pi_1(X^*)$ is abelian. Now all the 
assumptions of Theorem 8.1 in \cite{Za 2} are verified. 
By this theorem, $Y$ is contractible. 
\qed 

\smallskip

\noindent \cor {\it Let $X_0 = p^*(0),\,\,p \in \C{[k]},$ 
be a smooth reduced acyclic hypersurface in $\C^k.$ 
Then for any $n \in \N$ the hypersurface 
$X_n := \{u^nv = p(\bx)\}$ in $\C^{k+2}$ is smooth and contractible. 

If, furthermore, $\pi_1(\C^k \setminus X) \simeq \Z,$ 
then for any relatively prime integers $s_0,\,s_1 \in \N$ 
the hypersurface $Y_{s_0,\,s_1} :=\{u^{s_0}v^{s_1} = p(\bx)\}$ 
in $\C^{k+2}$ is smooth and contractible.}

\label{cor5.3}
\smallskip

{%\footnotesize
\noindent \rem \label{rm5.4} For instance, one may take as $X_0$ 
the tom Dieck-Petrie surface 
$X_{k,l} \subset \C^3$ (see Example \ref{ex3.1}). Indeed, 
it is easily seen that it satisfies all the conditions of 
Corollary \ref{cor5.3}.}

\section{$\C_+-$actions on the hypersurfaces $uv=p(x_1,\dots,x_k)$}

This section is devoted to the proof of the following theorem.

\smallskip

\noindent {\bf The Transitivity Theorem.} {\it Let $X = \{uv-p(\bx) = 0\} 
\subset \C^{k+2}$ where
$k \ge 2$ and $p \in \C^{[k]} \setminus \C.$ Then the automorphism group
Aut$\,X$ acts $m-$transitively on $X \setminus $sing$X$ for any $m \in \N.$}

\smallskip
{%\footnotesize
We keep all the notation from Section 4. Set 
$\sigma_i = \pi_i\,|\,X,\,\,i=1,\,2$ where 
$\pi_i\,:\,\C^{k+2} \to \C^{k+1},\,\,\,\pi_1\,:\,(\bx,\,u,\,v) 
\longmapsto (\bx,\,v),$ and $\pi_2\,:\,(\bx,\,u,\,v) 
\longmapsto (\bx,\,u),$
are the canonical projections. Then 
$\sigma_1\,:\,X \to \C^{k+1}=\C^k_{\bx} \times \C_v$ resp. 
$\sigma_2\,:\,X \to \C^{k+1}=\C^k_{\bx} \times \C_u$ 
is the affine modification of $\C^{k+1}$ along the hyperplane 
$D_v = \{v=0\}$ resp. $D_u = \{u=0\}$ 
with center\footnote{with center $C_1=X_0 \subset D_v$ resp. 
$C_2=X_0 \subset D_u$ if the divisor $p^*(0)$ is reduced.} 
$I_1 = (p,\,v) \subset \C^{[k+1]}$ resp. 
$I_2 = (p,\,u) \subset \C^{[k+1]}$
and with the exceptional divisor $V_0 \subset X$ resp. 
$U_0 \subset X$ 
(see Example \ref{ex1.6}). }

\smallskip

Concretizing Corollaries \ref{cor2.2} and \ref{cor2.3} 
in our setting we obtain the following statement.

\smallskip 

\noindent \lemma {\it Let $\var\,:\,\C_+ \times \C^{k+1} \to \C^{k+1}$ 
resp.
$\var\,:\, \C^{k+1} \to \C^{k+1}$ 
be a regular $\C_+-$action on $\C^{k+1} = \C^k_{\bx} \times \C_v$ 
resp. an automorphism of $\C^{k+1} = \C^k_{\bx} \times 
\C_v.$
Suppose that $\var$ leaves the subvarieties $D_v = \{v=0\}$ and  
$X_0$ invariant. 
Then there is a unique regular $\C_+-$action $\hvar$ on $X$ resp. 
an automorphism $\hvar$ of $X$
which leaves the hypersurface $V_0$ invariant and such that
the restriction $\hvar\,|\,(X \setminus V_0)$ coincides with 
$\sigma_1^{-1}\var\sigma_1.$ }

\label{lm6.1}
\smallskip 

\noindent {\bf Notation.} Let $G_1$ resp. $G_2$ 
be the subgroup of the group Aut$\,\C^{k+1}$ where 
$\C^{k+1}=\C^k_{\bx} \times \C_v$ resp. 
$\C^{k+1}=\C^k_{\bx} \times \C_u,$
generated by all the $\C_+-$subgroups $T$ of Aut$\,\C^{k+1}$ 
such that the function $v$ resp. $u$ is a $T-$invariant, and 
the restriction of $T$ to the invariant hyperplane $D_v$ resp. 
$D_u$ leaves the subvariety $X_0$ invariant.  
Denote $\hG_i,\,\,i=1,\,2,$ the subgroup of the group Aut$\,X$ 
which 
corresponds to $G_i$ in view of Lemma \ref{lm6.1}, and let 
$\hG \subset $Aut$\,X$ be the subgroup generated by $\hG_1$ and $\hG_2.$ 
Notice that 
$\hG_2 = \varepsilon \hG_1 \varepsilon^{-1}$ where 
$\varepsilon \in $ Aut$\,X,\,\,\varepsilon\,:\,(\bx,\,u,\,v) 
\longmapsto (\bx,\,v,\,u).$

\smallskip 

The Transitivity Theorem can be precised as follows.

\smallskip 

\noindent \thm {\it The group $\hG$ acts 
$m-$transitively on $X \setminus $sing$X$ for any $m \in \N.$}

\label{thm6.1}
\smallskip
{%\footnotesize
\noindent \rems.
\noindent \nrem \label{rm6.1'}
For $k=1$ a description of the automorphism group Aut$\,X$ 
of the surface $X = X_p \subset \C^3$ was given in \cite{ML 2}.
It follows from that description that the action of Aut$\,X$ 
on $X$ is transitive, but not 2-transitive. 

\noindent \nrem \label{rm6.1} 
Recall \cite[(9.2)]{KaML 1, Za 2} that the 
{\it Makar-Limanov invariant} of an algebra 
$A$ over $\C$ is the subalgebra ML$(A) \subset A$ which consists of 
all the elements invariant under every $\C_+-$subgroup of the 
automorphism group Aut$\,A;$ or, which is the same, 
ML$(A) = \bigcap_{\partial \in {\rm LND}(A)}$ ker$\,\partial$ 
(see Definition \ref{df2.4}). 
From Theorem \ref{thm6.1} it follows that 
for the algebra $A = \C[X]$ 
where $X$ is as above, this invariant is trivial,  
i.e. ML$(A) = \C.$ 
This is also true for $k=1;$ see \cite{ML 2}.
The problem arises to find a substitution of the 
Makar-Limanov invariant which would permit to distinguish 
the varieties $X=X(p)$ up to isomorphism, especially 
those diffeomorphic to the affine spaces.}

\smallskip 

The proof of Theorem \ref{thm6.1} is based on Lemmas 
\ref{lm6.2} - \ref{lm6.6}
below. In the next lemma for a class of $\C_+-$actions 
$\var = \var_{\partial}$ 
on $\C^{k+1}$ which preserve the decomposition 
$\C^{k+1} = \C^k_{\bx} \times \C_v$
we specify the lifts $\hvar \in $Aut$\,X$ of $\var,$ cf. Lemma \ref{lm6.1}. 

\smallskip

\noindent \lemma {\it Let $\delta$ be an LND of the polynomial algebra 
$\C^{[k]} = \C[x_1,\dots,x_k],$ and let $q \in \C[z]$
be a degree $d$ polynomial with the roots $z_1=0,\,z_2,\dots,z_m$ where 
$z_i \neq z_j$ for $i \neq j.$  
Then 

%%\smallskip

\noindent (a) the formulas
\be\partial (x_i) = q(v) \delta (x_i),\,\,i=1,
\dots,k,\,\,\,\,\,\,\partial (v) =0\,\ee
define an LND $\partial = \partial_{\delta,\,q}$ on $\C^{[k+1]}=
\C[x_1,\dots,x_k,\,v].$ The associated $\C_+-$action 
$\var_{\partial} \subset G_1$ on $\C^{k+1}$ fixes each point of 
the hyperplanes $D_j=\{v=z_j\},\,\,j=1,\dots,m.$

%%\smallskip

\noindent (b)  
The corresponding lifted LND ${\widehatpar}$ on $A = \C[X]$ is 
the restriction
to $X$ of the LND on $\C^{[k+2]} = \C[\bx,\,u,\,v]$ (denote it also by 
${\widehatpar}$) given by the formulas
\be {\widehatpar}(x_j) = \partial (x_j),\,\,\,{\widehatpar}(v) = 
\partial (v) = 0,\,\,\,{\widehatpar}(u) = 
{1\over v}\partial (p(\bx)) = {q(v)\over v} \delta p(\bx) \,.\ee
Furthermore, the fixed point set Fix$\,\hvar_{\partial}$ of the associated 
lifted $\C_+-$action $\hvar_{\partial} \subset \hG_1$ on $X$ contains the
union of the hypersurfaces $V_{z_i}=\{v=z_i\}$ for all the nonzero roots 
$z_i,\,\,i=2,\dots,m,$ of the polynomial $q;$ if $q'(0)=0,$ 
then it contains $V_0,$ too.}

\label{lm6.2}
\smallskip

\noindent  {\it Proof.~}  It is easy to check that $\partial$ resp. 
${\widehatpar}$ defined by the 
formulas (4) resp. (5) is, indeed, an LND of the algebra $\C^{[k+1]}$ resp.  
$\C^{[k+2]}=\C[\bx,\,u,\,v].$ 
From (5) it follows that
$${\widehatpar}(uv-p(\bx)) = 
v\cdot {q(v) \over v} \delta(p(\bx)) - q(v) \delta(p(\bx)) = 0\,,$$
and so, ${\widehatpar}(J) \subset J$ where $J \subset \C^{[k+2]}=\C[\bx,\,u,\,v]$ 
is the principal ideal generated by the polynomial
$uv-p(\bx).$ Therefore, ${\widehatpar}$ induces an LND of the quotient 
$\C[X] = \C[\bx,\,u,\,v]/J.$  
The other statements of the lemma can be verified without difficulty. 
\qed

\smallskip 
{%\footnotesize
\noindent \ex \label{ex6.1} In particular, putting in Lemma \ref{lm6.2}
$\delta =\delta_{i}:= {\partial \over \partial x_i} \in $
LND$\,(\C^{[k]})$ we obtain 
$\partial=\partial_i =\partial_{i,\,q}:= q(v)
{\partial \over \partial x_i} \in $ LND$\,(\C^{[k+1]})$ where 
LND$\,(A)$ denotes 
the set of all LND's of an algebra $A.$ 

For the `lifted' LND  ${\widehatpar}\in $ LND$\,(\C^{[k+2]})$ we have 
$${\widehatpar}_i(x_i) = q(v),\,\,\,{\widehatpar}_i(x_j) =0 \,\,\,\,
if \,\,\,\,j\neq i,\,\,\,{\widehatpar}_i(v) =0,\,\,\,
{\widehatpar}_i(u) =
v^{-1}q(v){\partial \over \partial x_i}p(\bx)\,.$$
It follows that ${\widehatpar}_i^{2}(x_i) =0$ and ${\widehatpar}_i^{n}(u) =
v^{-1}q^n(v){\partial^n \over \partial x_i^n}p(\bx)\,.$ The associated 
$\C_+-$action $\hvar_i:=\hvar_{\partial_i}$ on 
$\C^{[k+2]}$ is given as
$$\hvar^t_i (x_i) = x_i + tq(v),\,\,\,\hvar^t_i(x_j) 
= x_j\,\,\,\,
if \,\,\,\,j\neq i,\,\,\,\hvar^t_i(v) = v,\,\,\,{\rm and}$$
\be \hvar^t_i(u) = u + {1\over v} \sum_{n=1}^{{\rm deg}\, p} 
{t^n \over n!}q^n(v) {\partial^n \over \partial x_i^n}p(\bx)
,\,\,\,\,\,\,\,\,\,{\rm i.e.} \ee
$$\hvar^t_i (\bx_i,\,u,\,v) = 
\bigl(\bx + tq(v){\bar e}_i,\,u + {1\over v}[p(\bx + 
tq(v){\bar e}_i) - p(\bx)],\,v \bigr)\,,$$
where $({\bar e}_1,\dots,{\bar e}_k)$ is the standard basis of $\C^k.$}

\smallskip

In the following lemma we keep all the notation from Lemma \ref{lm6.2} 
and Example \ref{ex6.1}. We also denote by ${\cal O}_{\hvar}(P)$ 
the orbit of a point $P \in X$ under a $\C_+-$action $\hvar$ on $X.$

\smallskip

\noindent \lemma {\it Put $q(v) = v.$ Then:

%%\smallskip 

\noindent (a) For any point $P \in U_0 \setminus {\rm sing}\,U_0$ 
there exists $i\in \{1,\dots,k\}$ such that ${\cal O}_{\hvar_i}(P)
\not\subset U_0.$

%%\smallskip 

\noindent (b) For any point $P \in {\rm sing}\,U_0 \setminus {\rm sing}\,X$ 
there exists $\delta \in {\rm LND}\,(\C^{[k]})$ such that for the associated 
lifted $\C_+-$action $\hvar_{\partial}$ on $X$ one has 
${\cal O}_{\hvar_{\partial}}(P) \not\subset {\rm sing}\,U_0.$} 

\label{lm6.3}
\smallskip

\noindent {\it Proof.~}  (a) As follows from Lemma \ref{lm5.1}(b), $P = (\bx^0,\,0,\,v^0) 
\not\in {\rm sing}\,U_0$ iff grad$\,_{\bx^0}p \neq 0,$ i.e. 
${\partial p \over \partial x_i}(\bx^0)\neq 0$ for some $i\in \{1,\dots,k\}.$
By (6), the $u-$coordinate $u_i(t)$ of the point 
$\hvar_i^t(P_0) \in {\cal O}_{\hvar_i}(P)$ 
is a polynomial in $t$ with the non-zero linear term 
$({\partial p \over \partial x_i} (\bx^0)) t.$ Hence, 
${\cal O}_{\hvar_i}(P)\not\subset U_0,$ as claimed.

%%\smallskip

\noindent (b) Let 
$P = (\bx^0,\,0,\,v^0) \in {\rm sing}\,U_0 \setminus {\rm sing}\,X,$ that is, 
grad$\,_{\bx^0}p = 0,$ and $v^0 \neq 0.$ Fix a line $l$ through $\bx^0$ in 
$\C^k$ such that the restriction $p\,|\,l$ is non-constant. Chose a new 
coordinate system in $\C^k$ with the first coordinate axis being parallel to $l.$ 
Take 
$\delta = \delta_1 = {\partial \over \partial x_1} \in $ LND$\,(\C^{[k]})$ 
with respect to the new coordinates, so that $\hvar_{\partial} = \hvar_1.$
By (6), the projection $\bx^t = \bx^0 + tv^0{\bar e}_1$ of the point 
$\hvar^t_1(P) \in {\cal O}_{\hvar_1}(P)$ to $\C^k_{\bx}$ runs over $l.$ 
Since grad$\,p$ does not vanish identically on $l,$ we obtain that 
${\cal O}_{\hvar_1}(P)\not\subset {\rm sing}\,U_0.$ The proof is completed.
\qed

\smallskip 

\noindent \lemma {\it For any set of $m$ distinct points 
$P_1,\dots, P_m \in X \setminus {\rm sing}\,X$
there exists an automorphism $\hvar \in \hG_1$ of $X$
close enough to the identity such that 
$\hvar(P_i) \notin U_0,\,\,i=1,\dots,m.$}  

\label{lm6.4}
\smallskip

\noindent {\it Proof.~}  Starting with $m=0$ assume, by induction, that the statement 
is true for any set of $m$  distinct points in $X\setminus {\rm sing}\,X.$ 
Take an arbitrary set of $m+1$  distinct points 
$P_1,\dots, P_{m+1} \in X \setminus {\rm sing}\,X.$ 
By the inductive hypothesis, we may suppose that $P_1,\dots, P_m \not\in U_0$ and 
$P_{m+1} \in X \setminus {\rm sing}\,X.$ After applying to $P=P_{m+1},$ 
if necessary, an automorphism $\hvar^t_{\delta} \in \hG_1$ of 
Lemma \ref{lm6.3}(b) with $|t|$ small enough, we can achieve that 
$P_{m+1} \not\in {\rm sing}\,U_0,$
while still keeping $P_1,\dots, P_m \not\in U_0.$ If it occurs that 
$P_{m+1} \in U_0 \setminus {\rm sing}\,U_0,$ then applying Lemma \ref{lm6.3}(a)
to $P=P_{m+1}$ in the same way as above, we are done.
\qed

\smallskip 

\noindent {\bf Notation.}
Denote by $G_0$ the subgroup of the group Aut$\,\C^k$ generated by all the
$\C_+-$subgroups of Aut$\,\C^k.$ 

The proof of the following lemma can be found in 
\cite{Je, Ka 2}; for the sake of completeness we reproduce it here.  

\smallskip 

\noindent \lemma {\it For any $k \ge 2$ and any $m \in \N$ 
the group $G_0$ acts $m-$transitively on $\C^k.$ }

\label{lm6.5}
\smallskip

\noindent {\it Proof.~}  The proof goes by induction on $k,\,\,k \ge 2.$ 
Fix two arbitrary sets of $m$ distinct points $P_1,\dots, P_m$ and 
$Q_1,\dots, Q_m$ in $\C^k$ where $P_i = (\bx^{(i)},\,x_n^{(i)})$ and 
$Q_i = (\by^{(i)},\,y_n^{(i)})$ for certain 
$\bx^{(i)},\,\by^{(i)} \in \C^{k-1}$ and 
$x_n^{(i)},\,y_n^{(i)} \in \C,\,\,i=1,\dots,m.$
Choosing a generic coordinate system in $\C^k$ we may suppose that 
$\bx^{(i)} \neq \bx^{(j)},\,\,x_n^{(i)} \neq x_n^{(j)},\,\,
\by^{(i)} \neq \by^{(j)}$ and $y_n^{(i)} \neq y_n^{(j)}$ for all $i \neq j.$

If $n > 2$ then by the inductive hypothesis, we can find an automorphism 
$\alpha' \in G_0(\C^{k-1})$ such that 
$\alpha'(\bx^{(i)}) = \by^{(i)},\,\,i=1,\dots,m.$ 
After applying the automorphism 
$\alpha = (\alpha',\,{\rm id}_{\C})  \in G_0$ we may suppose that
$\by^{(i)} = \bx^{(i)},\,\,i=1,\dots,m.$ Let $p \in \C^{[k-1]}$ 
be a polynomial such that 
$p(\bx^{(i)}) = y_n^{(i)} - x_n^{(i)},\,\,i=1,\dots,m.$
Consider the triangular $\C_+-$action 
$\beta^t$ on $\C^k$ given as 
$\beta^t\,:\,(\bx,\,x_n) \longmapsto (\bx,\,x_n + tp(\bx)).$ 
Then the automorphism $\beta:= \beta^1 \in G_0$ sends $P_i$ to 
$Q_i,\,\,i=1,\dots,m.$ This provides the induction step.

For $k=2$ we start with the triangular automorphism  
$\alpha \in G_0, \,\,\,\alpha\,:\,
(x_1,\,x_2) \longmapsto (x_1 + q(x_2),\,x_2)$ where $q \in \C[z]$ is 
a polynomial such that $q(x_2^{(i)}) = y_1^{(i)} - x_1^{(i)},\,\,i=1,\dots,m,$
and then we apply $\beta \in G_0,$ as above. This completes the proof.
\qed

\smallskip

\noindent {\bf Notation.}
For any set of $n+1$ distinct non-zero
complex numbers $c_0,\,c_1,\dots,c_n \in \C 
\setminus \{0\}$ denote by 
${\rm Stab}_{c_1,\dots,c_n}\,(V_{c_0})$ 
the subgroup of the group 
$\hG_1$ which consists of the automorphisms 
of $X$ leaving the hypersurface 
$V_{c_0}$ invariant and fixing each point of 
the hypersurfaces $V_{c_i},\,\,i=1,\dots,n.$ 

\smallskip 

\noindent \lemma {\it (a) The group 
${\rm Stab}_{c_1,\dots,c_n}\,(V_{c_0})$ acts $m-$transitively 
on $V_{c_0}$ for any $m \in \N.$

%%\smallskip

\noindent (b) For any set of $m$ distinct points 
$P_1,\dots, P_m \in X \setminus {\rm sing}\,X$
there exists an automorphism $\hvar \in \hG$ of $X$
such that $\hvar(P_i) \in U_1,\,\,i=1,\dots,m.$}

\label{lm6.6}
\smallskip

\noindent {\it Proof.~}  (a) Fix two arbitrary sets of $m$ distinct points 
$P_1,\dots, P_m $ and $Q_1,\dots, Q_m $ in $V_{c_0}$ where
$P_i = (\bx^{(i)},\,u^{(i)},\,c_0),\,\,u^{(i)}=p(\bx^{(i)})/c_0,$ and  
$Q_i = (\by^{(i)},\,\tu^{(i)},\,c_0),\,\,\tu^{(i)} = p(\by^{(i)})/c_0,\,\,i=1,\dots,m.$ 
By Lemma \ref{lm6.5}, there exists an
automorphism $\alpha \in G_0$ such that 
$\alpha (\bx^{(i)}) = \by^{(i)},\,\,i=1,\dots,m.$

Decompose $\alpha = \psi_1^{t_1} \circ \dots \circ \psi_l^{t_l}$ 
into a product of elements of $\C_+-$subgroups of the group Aut$\,\C^k,$ 
and let $\tdelta_j \in $ LND$\,(\C^{[k]})$ be the infinitesimal generator 
of the subgroup $\{\psi_j^{t}\}_{t \in \C_+},\,\,j=1,\dots,l.$ 

Let also $q \in \C[z]$ be the degree $n$ polynomial with the roots $c_1,\dots,c_n$ 
such that $q(c_0)=1.$ Denote by
$\tpartial_j = \partial_{\tdelta_j,\,q}$ the LND of the algebra 
$\C^{[k+1]} = \C[\bx,\,v]$ defined as in (5) of Lemma \ref{lm6.2}(a),
and by $\tvar_j$ the corresponding $\C_+-$action on 
$\C^{k+1} = \C^k_{\bx} \times \C_v,\,\,j=1,\dots,l.$ Then we have:
$\tvar_j\,:\,(\bx,\,c_0) \longmapsto (\psi_j(\bx),\,c_0),$ and $\tvar_j$
fixes each point of the hyperplanes 
$D_{c_i} = \{v = c_i\},\,\,i=1,\dots,n.$
The composition $\talpha = \tvar_1^{t_1} \circ \dots \circ \tvar_l^{t_l}$
also fixes each point of the union $\bigcup_{i=1}^n D_{c_i},$ 
stabilizes the hyperplane $D_{c_0},$ and $\talpha\,|\,D_{c_0} = \alpha,$ i.e.
$\talpha(\bx,\,c_0) = (\alpha(\bx),\,c_0).$
Therefore, $\talpha (\bx^{(i)},\,c_0) = (\by^{(i)},\,c_0),\,\,i=1,\dots,m.$

By Lemma \ref{lm6.2}(b), the lift 
$\halpha = \hvar_1^{t_1} \circ \dots \circ \hvar_l^{t_l} \in \hG_1$ 
fixes each point of the union 
$\bigcup_{i=1}^n V_{c_i}$ and stabilizes the hypersurface $V_{c_0},$
that is, $\halpha \in {\rm Stab}_{c_1,\dots,c_n}\,(V_{c_0}).$ 
Moreover, since $\halpha\,|\,V_{c_0} = \si_1^{-1} \talpha \si_1 \,|\,V_{c_0}$
(see Lemma \ref{lm6.1})
we have $\halpha(P_i) = Q_i,\,\,i=1,\dots,m.$ This proves (a). 

%%\smallskip

\noindent (b) By Lemma \ref{lm6.4}, we may assume that  
$P_i = (\bx^{(i)},\,u^{(i)},\,v^{(i)}) \in X \setminus V_0,\,\,i=1,\dots,m.$ 
Reordering, if necessary, the points $P_i$ we may also suppose that 
$v^{(1)}=\dots=v^{(m')}=:c_0\neq v^{(j)}, \,\,j=m'+1,\dots,m.$ 
By (a), there exists an automorphism $\hvar \in \hG_1$ such that 
$\hvar(P_j) \in V_{c_0} \cap U_1,\,\,j=1,\dots,m',$ and 
$\hvar(P_j) = P_j,\,\,j=m'+1,\dots,m.$ Hence, proceeding by induction, 
we are done. 
\qed  

\smallskip

\noindent {\bf Proof of Theorem \ref{thm6.1}}. Fix an arbitrary set of $m$ 
distinct points $P_1,\dots, P_m$ in $X \setminus {\rm sing}\,X$ and another 
such set $Q_1,\dots,Q_m$ in $U_1.$ By Lemma \ref{lm6.6}(b), after applying, 
if necessary, to the points $P_1,\dots, P_m$ an automorphism $\hvar_1 \in \hG$ 
we may suppose that also $P_1,\dots, P_m \in U_1.$ Then, exchanging in 
Lemma \ref{lm6.6}(a) the roles of $U_1$ and $V_1,$ by this lemma, 
we can find an automorphism $\hvar_2 \in \hG_2$ such that 
$\hvar_2(P_i) = Q_i,\,\,i=1,\dots,m.$  The proof is completed.
\qed

\section{Examples of acyclic surfaces in $\C^3$ and of smooth contractible 4-folds 
$uv=p(x,\,y,\,z)$ in $\C^5$}

In the examples of smooth acyclic surfaces in $\C^3$ with 
big fundamental groups (see Example \ref{ex7.1} below) we use the following simple 
lemma (cf. \cite[Lemma 7.2]{Za 2}).

\smallskip

\noindent \lemma {\it Let $Y$ be a connected simply connected 
complex manifold, $F$ be a smooth irreducible hypersurface in 
$Y,$ and $p\,:\,X \to Y$ be a branched cyclic covering over $Y$
ramified to order $s$ on $F.$ 
Denote $G_X = \pi_1(X^*)$ resp. $G_Y = \pi_1(Y^*),$
where $X^* := X \setminus p^{-1}(F)$ and $Y^* := Y \setminus F.$ 
Identifying $G_X$ with the index $s$ subgroup 
$p_*(G_X) \subset G_Y$ we 
denote by $\beta_F$ a vanishing loop of $F$ in 
$Y$ and by $<<(\beta_F)^s>>$  
the normal closure in $G_Y$ of 
the cyclic subgroup generated by $(\beta_F)^s.$ 
Then
$\pi_1(X) \simeq G_X / <<(\beta_F)^s>>\,,$ and it 
is an index $s$ 
subgroup of the group
$G_Y / <<(\beta_F)^s>>\,$ with a cyclic quotient.}

\label{lm7.1}
\smallskip

\noindent {\it Proof.~}  We have the following commutative diagram: 

\begin{picture}(200,195)

\put(40,125){${\bf 1}\longrightarrow N_X 
\longrightarrow G_X = \pi_1(X^*) \stackrel{(i_X)_*}
{\longrightarrow} \pi_1(X) \longrightarrow {\bf 1}$}

\put(74,182){$\bf 1$}
\put(77,170){$\vector(0,-1){25}$}
\put(122,182){$\bf 1$}
\put(124,170){$\vector(0,-1){25}$}

\put(124,50){$\vector(0,-1){25}$}
\put(112,5){$\Z / s\Z$}

\put(78,110){$\vector(0,-1){25}$}
\put(83,95){$p_*$}
\put(124,110){$\vector(0,-1){25}$}
\put(130,95){$p_*$}
\put(220,110){$\vector(0,-1){25}$}
\put(225,95){$p_*$}

\put(40,65){${\bf 1}\longrightarrow N_Y \stackrel
{\simeq}{\longrightarrow} 
G_Y = \pi_1(Y^*) \stackrel{(i_Y)_*}{\longrightarrow}
\pi_1(Y)= {\bf 1} \,\,\,\,\,\, \,\,$}

\end{picture}

\noindent where 
$i_X\,:\, X^* \hookrightarrow X$ resp. $i_Y\,:\, Y^* \hookrightarrow Y$
denotes the identical embedding and $N_X := $ker$\,(i_X)_*$ resp. 
$N_Y := $ker$\,(i_Y)_*.$ Thus, $\pi_1(X) \simeq G_X / N_X,\,\,\,{\bf 1} = 
\pi_1(Y) \simeq G_Y / N_Y$ where $N_X = << \alpha_F >> \,$
and $N_Y = << \beta_F >>\,$
with $\alpha_F \in N_X$ resp. $\beta_F \in N_Y$ being a vanishing loop of 
$p^{-1}(F) \simeq F$ in $X$ resp. of $F$ in $Y;$ see e.g. 
\cite[(2.3.a)]{Za 2}. Moreover, identifying 
$G_X$ with the subgroup $p_*(G_X) \subset G_Y$ we may assume that 
$\alpha_F = (\beta_F)^s.$

Since the quotient $G_Y/G_X \simeq \Z/s\Z$ is abelian we have
$K_Y := [G_Y,\,G_Y] \subset G_X.$ Since $G_Y = N_Y = \,<<\beta_F>>\,$ 
the abelianization $(G_Y)_{\rm ab} := G_Y/K_Y$ 
is a cyclic group generated by the class $K_Y \beta_F.$ 
Hence, any element $g \in G_Y$ can be written as $g = g'(\beta_F)^t$ where
$g' \in K_Y \subset G_X$ and $t \in \Z.$ Thus, we have 
$g(\beta_F)^s g^{-1} = g'(\beta_F)^s g'^{-1}$ 
for any $g \in G.$ Therefore, the normal closure $N_X$ of the cyclic subgroup 
$<(\beta_F)^s = \alpha_F>$ in the group $G_X$ coincides with its normal 
closure in the bigger group $G_Y,$ i.e. $N_X \subset G_Y$ is a normal subgroup  
and it coincides with the subgroup $\,<<(\beta_F)^s>>\,.$ Hence, we have
$\pi_1(X) \simeq G_X / N_X \simeq  G_X / <<(\beta_F)^s>>\,,$ as required.

The chain of normal subgroups $N_X \subset G_X  \subset G_Y$ 
yields the short exact sequence  
$${\bf 1} \longrightarrow \pi_1(X) \simeq G_X / N_X  \longrightarrow 
G_Y / N_X \longrightarrow G_Y / G_X \simeq \Z / n\Z 
\longrightarrow {\bf 1}\,,$$ and the last assertion of the lemma follows.  
\qed

\smallskip
{%\footnotesize
\noindent \ex \label{ex7.1} Let  
$X_{k,l,s,m} = p_{k,l,s,m}^{-1}(0) \subset \C^3$ be the surface
defined by the polynomial 
$$p_{k,l,s,m}= {(xz^m + 1)^k - (yz^m + 1)^l - z^s \over z^m} 
\in \C[x,\,y,\,z]\,$$ 
where $0 \le m \le s.$ It is smooth if $m > 0,$ 
and it has at most one singular point $P_0 = (1,\,1,\,0)$ if $m=0.$
For gcd$(k,\,l) = $gcd$(k,\,s) = $gcd$(l,\,s) = 1$ and $m = s$ the surface
$Y_{k,l,s} := X_{k,l,s,s}$ is acyclic. 
Indeed, it can be presented as a cyclic $\C^*-$covering over the contractible 
tom Dieck-Petrie surface
$X_{k,l} = X_{k,l,1,1}\subset \C^3$ (see Example \ref{ex3.1}) 
branched to order $s$ along 
the line 
$L_{k,l} := X_{k,l} \cap \{z=0\} \simeq \C$ in $X_{k,l},$ and the 
acyclicity follows 
as in the proof of Theorem A in \cite{Ka 1} 
(see also \cite[\S 5]{Za 2}). 

However, in general the acyclic surface $Y_{k,l,s}$ is not contractible and 
possesses quite a big fundamental group.  Indeed, let 
$\sigma\,:\,\C^3 \to \C^3,\,\,\,\si(x,\,y,\,z) = (xz,\,yz,\,z),$ be the affine 
modification of $\C^3$ along the plane $z=0$ with center at the origin
(see Example \ref{ex1.4}). Then 
the restriction $\sigma\,|\,X_{k,l,s,m}\,:\,X_{k,l,s,m} \to X_{k,l,s,m-1}$ 
is the affine modification of $X_{k,l,s,m-1}$ 
along the line $D:= X_{k,l,s,m-1} \cap \{z = 0\}$ with center at the origin
(see Example \ref{ex2.1}). 

Furthermore, the surface $X_{k,l,s,1}$ coincides with the pseudoaffine 
modification of the smooth quasiprojective surface 
$X'_{k,l,s,0} := X_{k,l,s,0} \setminus \{P_0\}$ 
along the smooth curve $D^* := D \setminus \{P_0\}$ with center 
${\overline 0} \in D^*$ (see Definition \ref{def3.1}). 

By Corollary \ref{cor3.2}, the surface 
$Y_{k,l,s} = X_{k,l,s,s}$ being acyclic also the surfaces $X_{k,l,s,m}$ 
are acyclic for all $m=1,\dots,s.$ 
The repeated application of Lemma 3.4 in \cite{Ka 1} (or of Proposition 
\ref{pr3.1} above) yields the isomorphisms 
$$\pi_1(Y_{k,l,s}) = \pi_1(X_{k,l,s,s}) \simeq \pi_1(X_{k,l,s,s-1})
\simeq \ldots \simeq \pi_1(X_{k,l,s,1}) \simeq \pi_1(X'_{k,l,s,0})\,.$$  
The surface $X_{k,l,s,0} \simeq X_{k,l,s}:= \{x^k - y^l - 
z^s = 0\} \subset \C^3$ 
is homotopy equivalent 
to the cone over the Pham--Brieskorn 3-manifold 
$M_{k,l,s} := X_{k,l,s}\cap S^5,$ 
that is, over the link of the surface singularity of $X_{k,l,s}$
in the sphere $S^5.$ In turn, 
$X'_{k,l,s,0} \simeq X_{k,l,s} \setminus\{{\overline 0}\}$ 
is homotopy equivalent to the link $M_{k,l,s},$ and thus 
$\pi_1(Y_{k,l,s}) \simeq \pi_1(M_{k,l,s});$ denote the latter group by 
$ G'_{k,l,s}.$ 

The structure of the group $ G'_{k,l,s}$ is well known \cite{Mil 2}. 
It is finite iff $1/k + 1/l + 1/s > 1,$
infinite nilpotent iff $1/k + 1/l + 1/s = 1.$ If  $1/k + 1/l + 1/s \neq 1,$ 
then 
$G'_{k,l,s} = [G_{k,l,s},\,G_{k,l,s}]$  where 
$$G_{k,l,s}:=<\gamma_1,\,\gamma_2,\,\gamma_3\,|
\,\gamma_1^k=\gamma_2^l=\gamma_3^s=\gamma_1\gamma_2\gamma_3>$$ is a 
central extension of the Schwarz triangular group 
$$T_{k,l,s} := <b_1,\,b_2,\,b_3\,|\,b_1^2=b_2^2=b_3^2=1,\,\,(b_1b_2)^k = 
(b_2b_3)^l = (b_3b_1)^s =1>\,.$$
Note that for $1/k + 1/l + 1/s < 1$ the Schwarz triangular group $T_{k,l,s}$
contains a free subgroup with two generators. Therefore,  
the group $G'_{k,l,s}$ also contains such a subgroup;
in particular, it is not solvable. It is known \cite{Bri} 
that the Pham--Brieskorn manifold 
$M_{k,l,s}$ is a homology 3-sphere iff gcd$(k,\,l) = $gcd$(k,\,s) = 
$gcd$(l,\,s) = 1.$  
Under this condition the group $G'_{k,l,s}$ is perfect,
i.e. coincides with its commutator subgroup; indeed, its abelianization
$H_1(Y_{k,l,s};\,\Z)$ is trivial. Notice that the condition of relative 
primeness never holds in the Euclidean case $1/k + 1/l + 1/s = 1;$ 
in the spherical one 
$1/k + 1/l + 1/s > 1$ it holds only for the Kleinian icosahedral triple 
$(k,\,l,\,s) = (2,\,3,\,5).$ 

The isomorphism  
$X_{k,l}^*:=X_{k,l} \setminus L_{k,l} \simeq \C^2 \setminus \Gamma_{k,l}$ 
where $\Gamma_{k,l}:=\{x^k - y^l=0\}\subset \C^2$ \cite{tDP} provides the 
presentation     
$$B_{k,l}:=\pi_1(X_{k,l}^*) = <a,\,b\,|\,a^k = b^l>$$
(see e.g. \cite{Di 2}).
Since 
$Y_{k,l,s}^* := Y_{k,l,s} \setminus \{z=0\} \to X_{k,l}^*,\,\,\,
(x,\,y,\,z) \longmapsto (x,\,y,\,z^s),$ is a non-ramified cyclic covering, 
the group $\pi_1(Y_{k,l,s}^*)$ is isomorphic to 
an index $s$ subgroup, say, $\tC_{k,l,s}$ of the group $B_{k,l}$ 
with the cyclic quotient $B_{k,l}/\tC_{k,l,s} \simeq \Z/s\Z.$
We have 
${\rm ker}\,\left(i_{*}\,:\,\pi_1 (Y_{k,l,s}^*) \to 
\pi_1 (Y_{k,l,s})\right) = \,<<\alpha^s>>\,$ 
where $\alpha\in B_{k,l}$ is a vanishing loop of the line 
$L_{k,l} \subset X_{k,l}$ \cite[(2.3.a)]{Za 2}.  
It can be shown that  
$\alpha = a^qb^p \in B_{k,l}$ where $p,\,q \in \Z$ are such that 
$kp + lq =1.$ 

Therefore,  
for gcd$(k,\,l) = $gcd$(k,\,s) = $gcd$(l,\,s) = 1$ 
the group $G'_{k,l,s} \simeq \pi_1(Y_{k,l,s})$
is isomorphic to an index $s$ subgroup $C_{k,l,s}$ of the quotient 
$$B_{k,l,s}:= B_{k,l}/\,<<\alpha^s>>\, = 
 <a,\,b\,|\,a^k = b^l,\,\,(a^qb^p)^s = 1>\,$$
(see Lemma \ref{lm7.1} above).
In particular, for $k = 2,\,l=3$ in view of the isomorphism 
$X_{2,3}^* \simeq \C^2 \setminus \Gamma_{2,3}$
we have that $B_{2,3} = B_3$
is the 3-braid group with the generators 
$\sigma_1,\,\sigma_2\in B_3$ being vanishing loops of $L_{2,3}$ in $X_{2,3},$  
$a =  \sigma_1\sigma_2\sigma_1,\,\,\,b = \sigma_1\sigma_2.$ 
Therefore, $G'_{2,3,s}$ 
is isomorphic to an index $s$ subgroup $C_{2,3,s}$ of the group 
$$B_{2,3,s} = B_3/<<\sigma_1^s>> =
<\sigma_1,\,\sigma_2\,|\,\sigma_1\sigma_2\sigma_1=\sigma_2\sigma_1\sigma_2,
\,\,\sigma_1^s = \sigma_2^s = 1>\,$$
which consists of the elements of algebraic length\footnote{Recall 
that the  algebraic length of an element $\prod_{n=1}^m \si_{i_n}^{a_n}$ 
is the integer $\sum_{n=1}^m a_n.$} divisible by $s.$ 

%\smallskip

%\noindent \rem Observe that $X_{k,l,s,m}$ being an acyclic surface of 
%logarithmic Kodaira dimension 1, it  can be
%obtained starting with a Hirzebruch surface under the broken chains
%construction (see \cite{FlZa} for terminology). Actually, the construction of
%$X_{k,l,s,m}$ needs only three broken chains.
%More generally, every acyclic surface of logarithmic Kodaira dimension 1
%with three broken chains can be presented as a hypersurface
%in $\C^3$ (cf. \cite{KaML 2}). 

\smallskip 

\noindent \ex \label{ex7.2} 
Consider further
the 4-fold 
$X=X^{k,l,s,m}$ in $\C^5$ given by the equation 
$uv-p_{k,l,s,m}(x,\,y,\,z) = 0.$
By Corollary \ref{cor5.2}, the surface $X_{k,l,s,m}$ being 
acyclic implies that the hypersurface $X \subset \C^5$ 
is diffeomorphic to the Euclidean space $\R^8.$ 
But in general, as we have seen above, the surface
$X_0=X_{k,l,s,m}\subset \C^3$ is not contractible and possesses quite 
a big fundamental group. This shows that
Proposition \ref{pr5.2} cannot be extended to $k=3.$ 

\smallskip

Next we give an example of a polynomial $p=p_{2,\,3}^0 \in \C^{[3]}$ 
with a smooth acyclic (even contractible) zero fibre $F_0 = p^{-1}(0)$ 
and non-acyclic generic 
fibres $F_c = p^{-1}(c).$ In fact, the surface $F_0$ in this example 
is isomorphic to the tom Dieck-Petrie surface $X_{2,3}$ (see Example \ref{ex7.1}
above), but the embedding 
$X_{2,3} \stackrel{\simeq}{\longrightarrow} F_0 \subset \C^3$ 
is not equivalent to the standard one up to the action 
on $\C^3$ of the automorphism group Aut$\,\C^3.$ This provides also examples 
of non-equivalent embeddings of an exotic $\C^n$ into $\C^{n+1}.$ 

\smallskip  

\noindent \ex \label{ex7.3} Consider the affine modification 
$\si=\si_{{\overline 0},\,H}\,:\,\C^4 \to \C^4,\,\,\,\si\,:\,(x,\,y,\,z,\,t) 
\longmapsto (x,\,xy,\,xz,\,xt),$ of $\C^4$ along the hyperplane $H_0 = \{x=0\}$ 
with center at the origin (cf. Example \ref{ex1.4}). Consider also the Russell 
cubic $X \subset \C^4$ with the equation $$-x + x^2y + (z+1)^2 - (t+1)^3=0$$
(cf. Example \ref{ex1.7}). By Corollary \ref{cor2.1}, the restriction of $\si$
to the strict transform $X'$ of $X$ yields the affine modification 
$\si_{{\overline 0},\,B}\,:\,X' = 
\Sigma_{{\overline 0},\,B}(X) \to X$ of $X$ along the book-surface 
$B=H_0 \cap X \simeq \Gamma_{2,3} \times \C$ 
with center ${\overline 0} \in B \setminus {\rm sing}\ B.$ The hypersurface 
$X' \subset \C^4$ is  given by the equation 
$$-1 + x^2y + {(xz+1)^2 - (xt+1)^3 \over x} = 0\,$$
(cf. Example \ref{ex2.1}). The isomorphism  
$\C^3 \simeq X'$ as in Example \ref{ex1.7} provides 
an embedding of $\C^3$ into $\C^4.$ 
A direct computation shows that this embedding is rectifiable.\footnote{It is 
rectifiable e.g. via the composition 
$\gamma \circ \beta \circ \alpha$ of the triangular automorphisms 

\noindent $\alpha\,:\, (x,\,y,\,z,\,t) \longmapsto 
(x,\,y,\,u,\,t)\quad
{\rm where}\quad u =z + f(x,\,t) = 
z - xt^2 (xt + 3)/2\,,$ 

\noindent $\beta\,:\, (x,\,y,\,u,\,t) \longmapsto (x,\,v,\,u,\,t)
\quad
{\rm where}\quad v = y + g(x,\,u,\,t) = y + 
u t^2 (xt+3) + xt^4 (xt + 3)^2/4\quad {\rm and}$

\noindent $\gamma \,:\, (x,\,v,\,u,\,t) \longmapsto (x,\,v,\,u,\,w) 
\quad {\rm where}
\quad w = -3t + h(x,\,u,\,v) = -3t + x^2 v + x u^2 + 2u - 1 \,.$}

The hyperplane section $X_{2,3}^0 := X' \cap D_0$ where $D_0 = \{y=0\},$  
is isomorphic to the tom Dieck-Petrie surface $X_{2,3}$ 
(see Example \ref{ex7.1}). But the embedding 
$X_{2,3} \stackrel{\simeq}{\longrightarrow} X_{2,3}^0 
\hookrightarrow X' \simeq \C^3$ 
of the tom Dieck-Petrie surface into $\C^3$ is not equivalent to the 
standard one $X_{2,3}\hookrightarrow D_0 = \C^3.$ Indeed, the latter one
is defined by the polynomial 
$p_{2,3} = {(xz+1)^2 - (xt+1)^3 \over x} - 1 \in \C^{[3]}$ 
with all the fibres being contractible surfaces; see e.g. 
\cite[Example 6.1]{Za 2}. 
On the other hand, it is easily seen that
for a generic $c \in \C$ the fibre $F_c = p^{-1}(c)$ of the regular 
function
$p=p_{2,\,3}^0 := y\,|\,X' \in \C[X'] \simeq \C^{[3]}$ 
which defines the surface $X_{2,3}^0$
in $X'$ has the Euler characteristic $e(F_c) = 5$ ({\it hint}: 
use the fibration $F_c \to \C$ defined by the restriction $x\,|\,F_c$). 
In particular, the surfaces $F_c$ for a generic $c \in \C$ are not acyclic. 

Consider further the exotic product-structure $X_{2,3,n} := X_{2,3} \times 
\C^{n-2}$ on $\C^n,\,\,\,n \ge 3$ (see \cite[\S4]{Za 2}). By the similar 
arguments as above, two realizations $X_{2,3,n} := X_{2,3} \times \C^{n-2}
\hookrightarrow \C^{n+1} = \C^3 \times \C^{n-2}$ and $X_{2,3,n}^0 := 
X_{2,3}^0 \times \C^{n-2} \hookrightarrow \C^{n+1} $ of this exotic $\C^n$ 
as a hypersurface in $\C^{n+1}$ are not equivalent modulo the action on 
$\C^{n+1}$ of the automorphism group Aut$\,\C^{n+1}.$

\smallskip
%{%\footnotesize
\noindent \rems 

\noindent \nrem \label{rm7.2} By Corollary \ref{cor5.2}, 
the zero fibres of the polynomials $uv - p_{2,\,3}(x,\,y,\,z) \in \C^{[5]}$ 
and 
$uv - p_{2,\,3}^0 (x,\,y,\,z) \in \C^{[5]},$ as well as the generic fibres 
of the first one, are smooth contractible hypersurfaces in $\C^5,$ whereas 
the generic fibres of the second one have Euler characteristic $5$ 
(see Lemma \ref{lm5.1}(a)).

\smallskip

\noindent \nrem \label{rm7.3} More generally, 
we have presented above a collection of examples of smooth
acyclic hypersurfaces in $\C^n$ (see also \cite{Za 2}). 
Most of them are not rectifiable; 
actually, their defining polynomials have a 
fibre non-isomorphic to $\C^{n-1}.$ By Corollary \ref{cor3.2}
(see also Remark \ref{rm3.1}), performing 
the affine modification of $\C^n$ along such a hypersurface 
$D=p^*(0),\,\,p \in \C^{[n]},$ with a smooth reduced acyclic center 
$C \subset D$ leads to a smooth contractible affine 
$n-$fold $X.$  
The question arises when $X$ itself is isomorphic to the affine
space $\C^n$ (cf. Question A in sect. 7 below). 

Suppose that this is the case, and that, moreover, 
$C \simeq \C^k.$
Then the exceptional divisor $E = \si_C^{-1}(C) = q^*(0) \subset X$ where
$\si_C\,:\,X \to \C^n$ is the blowup morphism and 
$q := p\circ \si_C \in \C[X] \simeq \C^{[n]},$ would be isomorphic to $\C^{n-1}.$ 
But $q$ having a fibre non-isomorphic to $\C^{n-1},\,$
the hypersurface $E \simeq \C^{n-1}$ in $X \simeq \C^n$ 
cannot be rectifiable. 
Thus,
this would answer in negative to the Abhyankar-Sathaye 
Embedding Problem (cf. Remark \ref{rm5.3})\footnote{Recall 
that the varieties discussed in Remark \ref{rm5.3}
arise as affine modifications of the affine space 
$\C^n$ along a hyperplane with non-linear acyclic centers. 
Whereas here we consider, in particular, affine modifications of $\C^n$ 
along acyclic hypersurfaces with center at a linear subspace.}.  }

\section{On the theorems of Sathaye and Wright}

We propose the following

\smallskip

\noindent {\bf Questions}

%%\smallskip

\noindent {\bf A}) {\it When an affine modification 
$X = \Sigma_{I,\,f}(\C^n)$
of $\C^n$ is isomorphic to $\C^n?$}

%%\smallskip

\noindent {\bf B}) {\it Let $X  = \Sigma_{I,\,f}(\C^n)$ be an 
affine modification of $\C^n$ canonically embedded into $ \C^{n+k}$ 
(see Definition \ref{def1.2}).
If $X \simeq \C^n,$ is it necessarily rectifiable in $ \C^{n+k}$?}

%%\smallskip

\noindent {\bf C}) {\it Which embeddings $\C^n \hookrightarrow \C^{n+k}$
appear in this way?}

Question {\bf B} is a specialization of the 
Abhyankar-Sathaye Embedding Problem
mentioned in the Introduction. This problem is 
known to be answered in affirmative for $k \ge n+2$ (the 
Jelonek-Kaliman-Nori-Srinivas Theorem 
\cite{Je, Ka 2, Sr}) and for n=k=1 (the Abhyankar-Moh 
and Suzuki Embedding
Theorem \cite{AM, Suz}). In a special case when $n=2$ and $k=1$
the positive answer to Question {\bf B} is provided by 
a theorem of A. Sathaye\footnote{which corresponds to 
the case 
$n=1$ in Theorem \ref{thm4.1} below.}  \cite{Sat}, 
generalized by D. Wright\footnote{See also 
\cite{Ru 2, RuSat} for some generalizations.} \cite{Wr 2} as follows.

\smallskip

\noindent \nnthm{{\bf (Sathaye--Wright).}} 

\noindent {\it Let $X=X_{n,\,f,\,g}$ be a surface 
in $\C^3$ given by the equation 
$f(x,\,y)z^n + g(x,\,y) = 0$ where 
$f,\,g \in \C[x,\,y],\,\,\,n \in \N.$ 
Suppose that $X \simeq \C^2.$ Then $X$ is 
rectifiable, i.e. there exists an automorphism 
$\alpha \in $Aut$\,\C^2$ which transforms $X$ 
into a coordinate plane. }

\label{thm4.1}
\smallskip

We give below a new proof of Theorem \ref{thm4.1}, 
as well as some generalization. The proof is easy if one of the 
polynomials $f,\,g$ 
is constant; in the sequel we do not consider this possibility. 
Observe that for $n = 1$ the surface 
$X=X_{1,\,f,\,g}$ as in Theorem \ref{thm4.1}
is the affine modification $\Sigma_{I,\,f}(\C^2)$ of 
$\C^2$ along the divisor $D_f = f^*(0)$ with center $I = (f,\,g)$ 
(see Example \ref{ex1.6}) whereas 
the surface $X_{n,\,f,\,g}$ can be presented as a 
cyclic covering of $X_{1,\,f,\,g}$ ramified to order $n$ on 
$D_g=g^*(0).$ So, in the proof
we use affine modifications.

\smallskip

\noindent \thm  {\it Let $X=X_{n,\,f,\,g}$ be an irreducible 
smooth surface in $\C^3$ given by the equation 
$f(x,\,y)z^n + g(x,\,y) = 0$ where 
$f,\,g \in \C[x,\,y],\,\,n \in \N.$
The following conditions (i) - (iv) 
are equivalent\footnote{Hereafter 
$e(X)$ denotes the Euler characteristic of $X.$}:

%%\smallskip

\noindent (i) $e(X)=1$ and $H_1(X;\,\Z) = 0.$

%%\smallskip

\noindent (ii) $X$ is acyclic, i.e. ${\tilde H}_*(X;\,\Z) = 0.$

%%\smallskip

\noindent (iii) $X \simeq \C^2.$

%%\smallskip

\noindent (iv) $X$ is rectifiable.

%%\smallskip

\noindent If $n > 1,$ then the above conditions 
are equivalent to the following one:

%%\smallskip

\noindent (v) The pair $(f,\,g)$ is rectifiable in the following sense: 
there exists an automorphisms $\alpha \in $Aut$\,\C^{[2]}$ such that 
$(\alpha(f),\,\alpha(g)) = (p(x),\,y).$}

\label{thm4.2}
\smallskip

The implications $(v) \Longrightarrow (iv) \Longrightarrow (iii) 
\Longrightarrow (ii) \Longrightarrow (i)$ are easy; 
in the sequel we only prove $(i) \Longrightarrow (iv)$ in the case $n=1$ and 
$(i) \Longrightarrow (v)$ in the case $n>1.$

\smallskip

{%\footnotesize
\noindent \rems

\noindent \nrem There are examples of acyclic or even contractible
smooth algebraic surfaces in $\C^3$ non-isomorphic to $\C^2,$ see e.g.
\cite{tDP, KaML 2} and Example \ref{ex7.1} below. Moreover, 
any smooth contractible affine algebraic surface of logarithmic 
Kodaira dimension $1$ admits such an embedding into $\C^3$ \cite{KaML 2}. 
However, Theorem \ref{thm4.2} shows that the image of this embedding 
cannot be given by a `binomial' equation $f(x,\,y)z^n + g(x,\,y) = 0.$ 

\label{rem4.1}
%%\smallskip 

\noindent \nrem As a corollary we obtain that if the zero fibre of the polynomial 
$f(x,\,y)z^n + g(x,\,y)$ is acyclic (resp. isomorphic to $\C^2$) 
then so is every fibre. 
Notice that in general, a polynomial $p \in \C^{[3]}$ with a smooth acyclic, 
or even contractible, zero fibre may have non-acyclic 
generic fibres; see Example \ref{ex7.3} above. 

\label{rem4.2}
%%\smallskip 

\noindent \nrem The next simple observation will be useful in what follows. 
Let $X_1 = X_{1,\,f,\,g}$ be a surface as in 
Theorem \ref{thm4.2} with $n=1.$ Since $X_1$ 
is supposed being irreducible, the divisors $D_f$ 
and $D_g$ have no irreducible component in common. Thus, the center 
$D_f \cdot D_g$
of the blow up 
$\sigma_I\,:\,X_1 \to \C^2,\,\,\,\sigma_I(x,\,y,\,z) = (x,\,y),$
is supported by the finite set supp$\,D_f \cap $supp$\,D_g$ (see Example \ref{ex1.6}). 
The exceptional curve 
$E' \subset X_1,\,\,\,E' = \{f(x,\,y)=g(x,\,y)=0\},$ 
is isomorphic to the product $\C \times ($supp$\,D_f \cap $supp$\,D_g),$ 
and hence, it consists of $\kappa$ vertical lines in $X_1$ where
$\kappa:= $card (supp$\,D_f \cap $supp$\,D_g$).  

\label{rem4.3}}
\smallskip

The proof of Theorem \ref{thm4.2} is based on the following 
Lemmas \ref{lm4.1} - \ref{lm4.6}. 

\smallskip

\noindent \lemma {\it Let $X$ be a smooth irreducible affine surface. Then   

%%\smallskip

\noindent (a) The following conditions ($i$) and ($ii$) 
(resp. ($i'$) and ($ii'$)) are equivalent:

%%\smallskip

\noindent ($i$) $e(X)=1$ and $H_1(X;\,\Z) = 0$ resp. 
($i'$) $e(X)=1$ and $b_1(X) = 0;$

%%%\smallskip

\noindent ($ii$) $X$ is acyclic, i.e. ${\tilde H}_*(X;\,\Z) = 0,$ resp.
($ii'$) $X$ is $\Q-$acyclic, i.e. ${\tilde H}_*(X;\,\Q) = 0.$

%%\smallskip

\noindent (b) If ($i$) holds then the algebra $A = \C[X]$ is UFD.}

\label{lm4.1}
\smallskip

\noindent {\it Proof.~}  (a) The implications ($ii$) $\Longrightarrow$ ($i$) resp. 
($ii'$) $\Longrightarrow$ ($i'$) are straightforward.
To prove the converse ones, notice that by the Lefchetz 
Hyperplane Section Theorem \cite[Thm. 7.2]{Mil 1}, $X$ 
has homotopy type of a finite cell complex 
of real dimension at most two. 
Hence, $H_3(X;\,\Z) = H_4(X;\,\Z)=0$ and $H_2(X;\,\Z)$
is a free abelian group. Therefore, $e(X) = 1 -b_1(X) + b_2(X),$ and so, 
if $e(X) = 1$ and $b_1(X)= 0,$ then $b_2(X) = 0,$ and moreover, 
$H_2(X;\,\Z)= 0.$ Thus, ($i$) implies that ${\tilde H}_*(X;\,\Z)=0,$ i.e. 
$X$ is acyclic; in turn, ($i'$) implies that ${\tilde H}_*(X;\,\Q)=0,$ i.e. 
$X$ is $\Q-$acyclic. This proves (a).  

In view of (a), (b) follows from \cite[(1.17)--(1.20)]{Fuj} (see also
\cite[Prop. 3.2]{Ka 1}). \qed

\smallskip

\noindent \lemma {\it Let $f,\,g \in \C^{[2]} \setminus \C$ 
be two non-constant polynomials without common 
factor\footnote{For instance, this is so if the surface 
$X_n=X_{n,\,f,\,g}$ is irreducible.}. Then the surface 
$X_n=X_{n,\,f,\,g}=\{fz^n+g=0\} \subset \C^3$ is smooth 
iff the following two conditions are fulfilled:

\noindent (i) 
For any point $P_0 \in $supp$\,D_f \cap $supp$\,D_g$
the divisor $D_g$ is non-singular and reduced at $P_0.$ 
If so is the divisor $D_f$ at $P_0$ too, then $D_f$ and $D_g$ 
are transversal at $P_0.$ 

\noindent (ii) If $n > 1,$ then $D_g$ is a smooth reduced divisor.}

\label{lm4.2}
\smallskip

\noindent {\it Proof.~}  The statement easily follows from the equality 
$$\qquad\qquad\mbox{grad}\,(fz^n-g) = (z^n\mbox{grad}\,f - 
\mbox{grad}\,g,\,nz^{n-1}f)\,.\qquad \qquad \qquad \mbox{\qed}$$
\noindent \lemma {\it Let $X$ be a connected complex 
manifold, and let $D_1,\dots,D_s$ be reduced irreducible 
principal divisors in $X.$ Set
$D = \bigcup_{i=1}^s D_i$ and $X^* = X \setminus D.$ Then there is an exact sequence 
$${\bf 0} \longrightarrow \Z^s \stackrel{\mu}{\longrightarrow} H_1(X^*;\, \Z)
\stackrel{i_*'}{\longrightarrow} H_1(X;\, \Z) \longrightarrow {\bf 0}\,\,,$$
where $i\,:\,X^* \hookrightarrow X$ is the identical embedding, and $\mu$ sends
the standard basis vectors $(e_1,\dots,e_s)$ of the lattice $\Z^s$ into the
the vanishing loop classes $\alpha_1,\dots,\alpha_s$ of 
$D_1,\dots,D_s,$ respectively. Moreover, this sequence splits.}

\label{lm4.3}
\smallskip

\noindent {\it Proof.~} 
We have the following exact sequence of fundamental groups:
$${\bf 1} \longrightarrow \,<<\alpha_1,\dots,\alpha_s>>\, 
\longrightarrow \pi_1(X^*)\stackrel{i_*'}{\longrightarrow} \pi_1(X) 
\longrightarrow {\bf 1}$$
where $<<\alpha_1,\dots,\alpha_s>>$ denotes the normal subgroup 
of the group $\pi_1(X^*)$ generated by these vanishing loops 
(see e.g. \cite[(2.3.a)]{Za 2}). Passing to the abelianizations shows that 
the above homology sequence is exact besides, possibly, at the second term. 
The divisor $D_i'$ being principal we have $D_i'= g_i^*(0)$ where 
$g_i$ is a holomorphic function on $X,\,\,i=1,\dots,s.$ The morphism
$\varphi = (g_1,\dots, g_s)\,:\,X_1^* \to (\C^*)^s$ yields a surjection 
$$\varphi_*\,:\,H_1(X^*;\,\Z) \longrightarrow\!\!\!\to \Z^s,\,\,
\varphi_*(\alpha(D_i')) = e_i,\,\,i=1,\dots,s.$$
This provides the exactness at the second term, the splitting 
$H_1(X^*;\, \Z) \simeq H \oplus $ ker$\,\varphi_*$ where 
$H = <\alpha_1,\dots,\alpha_s>,$ and an isomorphism 
$H_1(X;\,\Z)\simeq $ ker$\,\varphi_*.$
\qed

\smallskip

\noindent \lemma {\it (a) If the surface $X_n=X_{n,\,f,\,g}$ is smooth 
resp. irreducible resp. $\Q-$acyclic, then so is $X_1=X_{1,\,f,\,g}.$

%\smallskip

\noindent (b) Suppose that the surface $X_n$ where $n > 1,$ is 
smooth and irreducible, and
$H_1(X_n;\,\Z) =0.$ Then also  $H_1(X_1;\,\Z)=0,$ and the divisor $D_g$ is  
irreducible. Hence, if $X_n$ is acyclic, then so is $X_1.$}

\label{lm4.4}
\smallskip

\noindent {\it Proof.~}  (a) 
The first two statements of (a) easily follow from Lemma \ref{lm4.2}. 
As for the third one, consider the cyclic ramified covering $\C^3 \ni
(x,\,y,\,z) \longmapsto (x,\,y,\,z^n) \in \C^3,$ which restricts to 
a cyclic covering $\si_n\,:\,X_n \to X_1$ branched to order $n$ over the curve 
$D_g' = \{z=0\}$ in $X_1.$ For any prime $p$ 
which does not divide
$n$ the transfer provides an isomorphism of the homology group 
$H_*(X_1;\,\Z_p)$ and the subgroup of the homology group $H_*(X_n;\,\Z_p)$
fixed by the monodromy action \cite[III(2.4)]{Bre}. 
Thus, $\Z_p-$acyclicity of $X_n$ implies 
$\Z_p-$acyclicity of $X_1.$ If $X_n$ is $\Q-$acyclic, then it is 
$\Z_p-$acyclic for all but finite number of 
the primes $p,$ and the same holds for $X_1.$ Therefore, 
$X_1$ is $\Q-$acyclic, too. This proves (a).

%%\smallskip

\noindent (b)   
Let $D_g'$ resp. $D_g''$ be the divisor in $X_1$ resp. in $X_n$ 
given by the equation $z=0.$ Set
$X_1^*=X_1\setminus D_g'$ and $X_n^*=X_n \setminus D_g''.$
Then $\si_n\,|\,X_n^* :\,X_n^* \to X_1^*$ is an $n-$sheeted unramified 
covering, and hence, $(\si_n)_* (H_1 (X_n^*;\,\Z))$ is  
a subgroup of index at most $n$ of the group $H_1 (X_1^*;\,\Z).$ 

Let $\alpha(D_i') \in H_1(X_1^*;\, \Z)$ resp. 
$\alpha(D_i'') \in H_1(X_n^*;\, \Z)$ be
the vanishing loop class  of the irreducible component 
$D_i'$ of $D_g'$ resp. $D_i''$ of $D_g'',\,\,\, i=1,\dots,s.$
By Lemma \ref{lm4.3}, we have the following commutative diagram 
where the horizontal lines are exact sequences:

\begin{picture}(200,95)

\put(40,65){${\bf 0}\longrightarrow 
\Z^s \stackrel{\simeq}{\longrightarrow} 
H_1(X_n^*;\, \Z) \stackrel{i_*''}{\longrightarrow} 
H_1(X_n;\, \Z)= {\bf 0} \longrightarrow {\bf 0}$}

\put(77,50){$\vector(0,-1){25}$}
\put(81,35){$\tau$}
\put(137,50){$\vector(0,-1){25}$}
\put(141,35){$(\si_n)_*$}
\put(215,50){$\vector(0,-1){25}$}
\put(219,35){$(\si_n)_*$}

\put(40,5){${\bf 0}\longrightarrow \Z^s \stackrel{\mu}
{\longrightarrow} H_1(X_1^*;\, \Z) \stackrel{i_*'}
{\longrightarrow} H_1(X_1;\, \Z)
\longrightarrow {\bf 0}\,\,.$}

\end{picture}
 
\noindent Clearly, $(\si_n)_*(\alpha(D_i'')) = n \alpha(D_i'),$ and so,
$\tau(e_i) = ne_i, \,\,i=1,\dots,s.$ 
Thus, the image $(\si_n)_*(H_1(X_n^*;\,\Z))$ is a subgroup of 
index at least $n^s$ 
of the group ker$\,i'_* \subset H_1(X_1^*;\,\Z).$ 
On the other hand, it should be a subgroup of index at most $n.$ 
Therefore,
$s=1,$ i.e. $D_g$ is irreducible, and, furthermore, 
ker$\,i'_* = H_1(X_1^*;\,\Z).$ Hence, $H_1(X_1;\,\Z)=0.$ 

The last statement of (b) follows from (a) and Lemma \ref{lm4.1}(a).
\qed
 
\smallskip

\noindent \lemma {\it If the surface $X_n=X_{n,\,f,\,g}$ is smooth, 
irreducible and $\Q-$acyclic, then the following assertions hold. 
 
\smallskip

\noindent (a) $ e(D_f\setminus D_g) = 0, \,\, $
or, which is equivalent, $e(D_f) = \kappa := e(D_f \cap D_g).$

\noindent Furthermore, if $n>1,$ then $e(D_g) = e(D_f \cup D_g)=1.$ 
If, in addition, $H_1(X_n;\,\Z)=0,$
then the divisor $D_g$ is smooth, reduced, 
irreducible and isomorphic to $\C.$

\smallskip

\noindent (b)
The irreducible components $D_f^{(i)},\,\,i=1,\dots,k,$ of the divisor $D_f$ 
are disjoint simply connected curves smooth outside of $D_g;$ each of them 
meets supp$\,D_g$ at one point, which is smooth and reduced on $D_g.$ 
In particular, $k=e(D_f) = e(D_f \cap D_g) = \kappa \ge 1.$ 

\noindent If, in addition, $H_1(X_1;\,\Z)=0,$ then each of the curves 
$D_f^{(i)},\,\,i=1,\dots,k,$ is smooth and meets supp$\,D_g$ 
transversally.}

\label{lm4.5}
\smallskip

\noindent {\it Proof.~}  (a) Denote by $E'' $ resp. by $D_g''$ the curve 
in $X_n$ given by the equations $f=g=0$ resp. $z = 0.$ 
Consider the disjoint constructive decompositions 
$$X := X_n = X' \cup X'' \cup X''' \qquad \mbox{and} \qquad  Y :=\C^2= Y' 
\cup Y'' \cup Y'''\,,$$
where $$X' = X \setminus (E'' \cup D_g''),\quad X'' = D_g'' 
\simeq \mbox{supp} D_g, \quad X''' = E'' \setminus D_g''$$
(so, $X'''$ is a disjoint union of curves isomorphic to $\C^*$), and
$$Y' = \C^2 \setminus (D_f \cup D_g),\quad Y'' = D_g,\quad Y''' = 
D_f \setminus D_g\,.$$
Clearly, $\si_I\,|\,X'\,:\,X' \to Y'$ is a non-ramified $n-$sheeted 
covering, and hence, $e(X') = ne(Y').$ By the additivity of the 
Euler characteristic \cite{Du},
we have
\be 1= e(X) = e(X') + e(X'') + e(X''') = ne(Y') + e(D_g)\,,\ee and
\be 1= e(Y) = e(Y') + e(Y'') + e(Y''') = e(Y') + e(D_g)+ 
e(D_f \setminus D_g)\,.\ee
Subtracting (2) from (1) we obtain
\be (n-1)e(Y') = e(D_f \setminus D_g)=e(D_f )-\kappa\,.\ee
Putting here $n=1$ we obtain the equalities in (a). Since by 
Lemma \ref{lm4.4}(a), the surface
$X_1:=X_{1,\,f,\,g}$ is still $\Q-$acyclic, the case $n>1$ 
can be reduced to the case $n=1.$ 

Now, if $n >1$ we obtain from (3) the equality 
$e(Y') = 0.$ By the definition of $Y',$ this yields the equality
$e(D_f \cup D_g) =1,$ and also, by (1), we have $e(D_g)=1.$ 

If, in addition, $H_1(X_n;\,\Z)=0,$
then by Lemmas \ref{lm4.2} and \ref{lm4.4}(b), the divisor $D_g$ 
is smooth, reduced and irreducible. Since $e(D_g)=1$ it is isomorphic to 
$\C.$ This proves (a).

%%\smallskip

\noindent (b) is proven in Claims 1-6 below. 

\smallskip 

\noindent {\bf Claim 1.} {\it Each irreducible component of $D_f$ 
meets supp$\,D_g.$}

%%\smallskip

\noindent {\it Proof.~}  
Assuming that there exists an irreducible component of $D_f$
which does not meet supp$\,D_g$ we would get a decomposition 
$f = f_1f_2$ where $f_1 \not= $const and
$(f_1,\,g)=\C^{[2]}.$ Let $\xi,\,\eta \in \C^{[2]}$ be 
polynomials such that
$\xi f_1 + \eta g = 1.$ Replacing $\eta g = 1 - \xi f_1$ we obtain 
the relation $f_1(\xi + \eta f_2z^n)\,|\,X = 1,$ i.e. $f_1\,|\,X \in \C[X]$ 
is a non-constant invertible regular function, which is impossible. 
Indeed, since $b_1(X)=0$ 
the regular function $f_1\,|\,X$ can be expressed as $\exp(\psi)$ 
where $\psi$ is a non-constant holomorphic function on $X.$ 
But then $f_1\,|\,X$
cannot be regular, a contradiction.  \qed

\smallskip 

\noindent {\bf Claim 2.} {\it Each connected component of the curve
supp$\,D_f$ is simply connected and meets $D_g$ at one point only.}

%%\smallskip

\noindent {\it Proof.~}  Let $D_1,\dots,D_s$ be the connected components 
of the curve supp$D_f.$
For a connected affine curve $\Gamma$ we always have
$e(\Gamma) \le 1,$ and $e(\Gamma) = 1$ iff $b_1(\Gamma) =0,$ i.e. 
iff $\Gamma$ is simply connected (see e.g. \cite{Za 1}).    
From this observation and Claim 1 it follows that 
$e(D_i \setminus D_g) \le 0,\,\,i=1,\dots,s.$ Since we have by (a),
$$0 = e(D_f \setminus D_g) = \sum_{i=1}^s e(D_i \setminus D_g)\,,$$
we obtain that all the summands in the latter sum vanish. 
Thus, $0 = e(D_i \setminus D_g) = e(D_i) - e(D_i \cap D_g).$ 
Together with Claim 1 and the above observation this 
yields the inequalities $1 \le e(D_i \cap D_g) = e(D_i) \le 1.$ 
Therefore,
$e(D_i \cap D_g) = e(D_i) = 1,$ and thus $b_1(D_i) = 0,$ i.e. 
$D_i$ is simply connected and meets $D_g$ at one point only, 
$i=1,\dots,s.$ \qed

\smallskip

\noindent {\bf Claim 3.} {\it $k \le \kappa$ where $k$ is 
the number of irreducible components of $D_f.$ }

%%\smallskip

\noindent {\it Proof.~}  
Denote $D_f^{(i)},\,\,i=1,\dots,k,$ resp. $D_g^{(j)},\,\,j=1,\dots,l,$
the irreducible components of the divisor $D_f$ resp. $D_g.$ 
In the notation as in the proof of (a) above, consider 
the non-ramified $n-$sheeted covering $\si := \si_I\,|\,X'\,:\,X' \to Y'.$ 
It is easily seen that 
$H_1(Y';\,\Z) \simeq \Z^{k + l}\,.$ As in the proof of 
Lemma \ref{lm4.3} we have the exact sequence
$$\Z^{\kappa + l} \longrightarrow H_1(X', \Z)
\longrightarrow H_1(X, \Z) \longrightarrow {\bf 0}\,\,.$$
Since $H_1(X, \Z)$ 
is a torsion group it follows that $H_1(X', \Z)$ contains a subgroup of a 
finite index which is a homomorphic image of $\Z^{\kappa + l}.$
The image $\si_*(H_1(X';\,\Z))$ is a finite index subgroup of the group
$H_1(Y';\,\Z).$ 
Henceforth, $\kappa + l \ge k + l,$ or $\kappa  \ge k.$ \qed

\smallskip

\noindent {\bf Claim 4.} {\it $k = \kappa,$ and the connected components of 
the divisor $D_f$ coincide with its irreducible components.}

%%\smallskip

\noindent {\it Proof.~}  Indeed, by (a) and 
Claims 2, 3 we have $\kappa = e(D_f) = s \le k \le \kappa.$ Hence, 
$s=k=\kappa,$ and the claim follows. \qed

\smallskip

\noindent {\bf Claim 5.} {\it $D_f\setminus D_g$ is 
a smooth divisor. Furthermore, if $H_1(X_1;\,\Z)=0,$ then
the divisor $D_f$ is smooth.} 

%%\smallskip

\noindent {\it Proof.~}  Let $D=D_f^{(i)}$ be an irreducible component of 
the divisor $D_f.$ It is simply connected and meets $D_g$ 
at a unique point, say, $P.$ 
Let $\nu\,:\,\C \to D$ be a normalization map such that $\nu(0)=P.$ 
Then the polynomial $r:=g\circ\nu \in \C^{[1]}$ vanishes only at zero. 
This implies that 
the curve $D \setminus \{P\}$ is smooth (indeed, as $r(t) = ct^m$ for some 
$c \in \C^*$ and $m \in \N,$
the derivative $r'$ vanishes only at the origin). 

It remains to show that $D$ is smooth at $P$ providing that 
$H_1(X_1;\,\Z)=0.$
Assume on the contrary that $P$ is a singular point of $D.$ Denote 
$A_1 = \C[X_1].$ Set 
$\var = f_1 \circ \si_I \in A_1$ where $f_1 \in \C^{[2]}$ is an 
irreducible polynomial which defines $D,$ and 
$\si_I\,:\,X_1 \to \C^2,\,\,\,I = (f,\,g) \subset \C^{[2]},$ 
is the blowup morphism. Then $\var^{-1}(0) =
\si_I^{-1}(P)=: E_P'$ is the irreducible component over $P$
of the exceptional divisor $E' \subset X_1.$
We claim that $E_P'$ is a multiple fibre of $\var$ whereas 
its generic fibres are irreducible. 
Indeed, since the generic fibres of $f_1\,|\,(\C^2 \setminus D_f)$ 
are irreducible, in view of the isomorphism $\si_I\,|\,(X_1 \setminus E')\,:\,
X_1 \setminus E' \stackrel{\simeq}\longrightarrow \C^2 \setminus D_f$ 
the same is true for $\var\,|\,(X_1 \setminus E').$ Further, 
we have the equalities 
$${\rm grad}\,(f_1\circ\pi)\,|\,E_P' = ({\rm grad}_P\,f_1,\,0) = 
{\ol 0}\,,$$  
where $\pi\,:\,\C^3 \to \C^2,\,\,\,\pi(x,\,y,\,z)=(x,\,y).$
It follows that (grad$\,\var)\,|\,E_P' = 0$ in local coordinates in $X_1.$ 
 
By Lemma \ref{lm4.4}(a), the surface $X_1$ is $\Q-$acyclic. Since 
$H_1(X_1;\,\Z)=0,$ by Lemma \ref{lm4.1}, actually it is acyclic, 
and the algebra $A_1$ is UFD. Thus, in view of $\var^*(0) = mE_P',$ 
where $m > 1,$ 
we have that $\var = \var_1^m$ for a certain $\var_1 \in A_1.$ 
Therefore, the generic fibres of $\var$ cannot be irreducible, 
which is a contradiction. \qed

\smallskip

\noindent {\bf Claim 6.} {\it If $H_1(X_1;\,\Z)=0,$ then 
the curves supp$\,D_f$ and supp$\,D_g$ meet transversally.}

%%\smallskip

\noindent {\it Proof.~}  We keep all the notation from the proof of Claim 5.
Assume that an irreducible component 
$D := $supp$\,D_f^{(i)}$ is tangent to $D_g$ at their unique 
intersection point $P.$
By Lemma \ref{lm4.2}, the divisor $D_g$ is smooth and reduced at $P,$ and
$D_f = mD+\dots$ for some $m > 1.$ 
Thus, we have $f = f_1^mf_2$ where $f_1 \in \C^{[2]}$ is an irreducible 
polynomial which defines $D.$ As above, the generic fibres of the regular function 
$\var = f_1 \circ \pi\,|\,X_1 \in A_1=\C[X_1]$ are irreducible.
We claim that $E_P' :=\var^{-1}(0)$ is a multiple fibre of $\var,$    
which contradicts to the fact that $A_1$ is UFD (see the proof of Claim 5). 

Indeed, by our assumption, we have grad$\,_P\,f = \gamma$ grad$\,_P\,g$ 
for some $\gamma \in \C.$ 
This implies that (grad$\,\var)\,|\,E_P' = 0$ in local coordinates in $X,$ or,
which is the same, that ${\rm grad}\,(f_1\circ\pi)\,|\,E_P'$ is proportional to
(grad$\,F)\,|\,E_P'$ where $F := f_1^mf_2z - g \in \C^{[3]}$ is 
the defining polynomial of the surface $X_1$ in $\C^3.$ 
The latter follows from the equalities
$$({\rm grad}\,F)\,|\,E_P' = (mzf_1^{m-1}f_2 \,{\rm grad}\,f_1 + 
zf_1^m \,{\rm grad}\,f_2 - {\rm grad}\,g,\,f_1^mf_2)\,|\,E_P' $$
$$= (-{\rm grad}_P\,g,\,0)=(-\gamma\,{\rm grad}_P\,f_1,\,0)=
-\gamma\,{\rm grad}\,(f_1\circ\pi)\,|\,E_P'\,. $$
Henceforth, (grad$\,\var)\,|\,E_P' = 0,$ as claimed. \qed

Now the proof of (b) is completed. \qed

\smallskip

\noindent \lemma {\it (a) Suppose that the surface $X_n$ where $n > 1,$ 
is smooth and irreducible. Then the implication $(i) \Longrightarrow (v)$ 
of Theorem \ref{thm4.2} holds. 

%%\smallskip

\noindent (b) Suppose further that the surface $X_1$ is smooth, irreducible 
and $\Q-$acyclic. Assume also that the curves supp$D_f$ and supp$D_g$ meet 
transversally\footnote{In particular, as follows from Lemma \ref{lm4.5}(b), supp$D_f$ 
is smooth.}. Then, by an appropriate automorphism $\alpha' =
(\alpha,\,{\rm id}) \in {\rm Aut}\,(\C^2 \times \C)$ of $\C^3,$ the equation of 
$X$ can be reduced to the following form:} 
$$p(x)z = y + xh(x,\,y) + {\rm const}.$$ \label{lm4.6} 
{\it Proof.~}  (a) If the condition ($i$) of Theorem \ref{thm4.2} is fulfilled then 
by Lemmas \ref{lm4.1}(a) and \ref{lm4.4}(b), both surfaces $X_1$ and $X_n$ are acyclic. 
Thus, by Lemma \ref{lm4.5}, $D_g \simeq \C$ and $D_f^{(i)} \simeq \C$ 
for each $i=1,\dots,\kappa.$ 

Let $f_1 \in \C^{[2]}$ be an irreducible polynomial which 
defines the irreducible component $D_f^{(1)}$ of $D_f.$ Since 
by Lemma \ref{lm4.5}(b), the components 
$D_f^{(i)},\,\,i=1,\dots,\kappa,$ of $D_f$ are disjoint and simply connected, 
for each $i=1,\dots,\kappa$ we have 
$f_1\,|\,D_f^{(i)}\equiv {\rm const} = c_i,$ 
and therefore, $f=p(f_1)$ for some polynomial $p \in \C^{[1]}.$
In view of Lemma \ref{lm4.5}(b), the curves 
$D_f^{(1)}\simeq \C$ and $D_g \simeq \C$ meet transversally at a unique point.
By the Abhyankar-Moh and Suzuki Embedding Theorem, 
it follows that
$\alpha := (f_1,\,g)\,:\,\C^2 \to \C^2$ 
is an automorphism which transforms the pair 
$(f_1,\,g)$ into a pair of coordinate functions $(x,\,y),$ and 
the pair $(f,\,g)$ into the pair $(p(x),\,y),$ as desired. This proves 
the implication ($i$) $\Longrightarrow$ ($v$) 
of Theorem \ref{thm4.2}.  

%%\smallskip

\noindent (b) 
By the Abhyankar-Moh and Suzuki Embedding Theorem, the smooth simply 
connected curve
$D_f^{(1)} \simeq \C$ can be transformed into a coordinate line 
$x = 0$ by an automorphism of $\C^2.$ 
Thus, we may assume that $f_1(x,\,y)=x,$ and so, 
as above, $f(x,\,y) = p(x)$ and $X = \{p(x)z = g(x,\,y)\} \subset \C^3.$

If $x_i$ is a root of $p,$ that is, 
$D_f^{(i)} = \{x = x_i\},$ then by Lemma \ref{lm4.5}(b), the polynomial
$g(x_i,\,y) \in \C[y]$ has only one root, say, $c_i,$ and this root is  
simple, that is,
$g(x_i,\,y) = \gamma_i(y-c_i).$ Hence, $g(x,\,y) =  \gamma_i(y-c_i) +
(x - x_i)h_i(x,\,y)$ for a certain $h_i \in \C^{[2]}.$ Plugging here
$x_1 = 0$ and replacing $z$ by $z/\gamma_1,$ we obtain 
the desired presentation.
\qed

\smallskip

This lemma and the next proposition yield the implication 
$(i) \Longrightarrow (iv)$ of Theorem \ref{thm4.2} for $n = 1.$ 
Denote by $G$ 
the subgroup of the group Aut$\,\C^3$ which consists of the automorphisms
of the type $(x,\,y,\,z) \longmapsto (x,\,\gamma_1y+xg_1,\,\gamma_2z+xg_2)$ 
where $g_i \in \C^{[3]}$ and $\gamma_i \in \C^*,\,\,i=1,\,2.$

\smallskip

\noindent \prop {\it Let $X$ be an irreducible smooth $\Q-$acyclic surface 
given in $\C^3$ by the equation $p(x)z = g(x,\,y)$ where 
$p \in \C^{[1]}\setminus \C$ and $g \in \C^{[2]}\setminus \C.$  
Suppose that the curves supp$D_p$ and supp$D_g$ meet 
transversally\footnote{By Lemma \ref{lm4.5}(b), this is true if $X$ 
satisfies the condition ($i$) of 
Theorem \ref{thm4.2}.}. 
Then $X$ can be transformed into 
a plane $L_c :=\{y=c\},\,\,c \in \C,$ by an automorphism $\alpha \in G.$}

\label{pr4.1}
\smallskip

\noindent The proof proceeds by induction on deg$\,p.$ 

\smallskip

\noindent {\bf Claim 1.} {\it The statement is true if ${\rm deg}\,p=1.$}

%%\smallskip

\noindent {\it Proof.~}  We may assume that $p(x) = x,$ and by Lemma \ref{lm4.5}(b),
that $g(x,\,y) = y + xh(x,\,y) - c,\,\,c \in \C.$ Then the automorphism
$\alpha \in $Aut$\,\C^3,\,\,\,\alpha(x,\,y,\,z)=(x,\,y,\,z+h),$ transforms 
$X$ into the surface $X' := \{xz = y - c\}.$ Furthermore, the automorphism
$\beta \in G,\,\,\,\beta(x,\,y,\,z)=(x,\,y+xz,\,z),$ transforms 
$X'$ into the plane $L_{c}.$ The resulting automorphism 
$\alpha' := \beta\circ\alpha,$
$\alpha'(x,\,y,\,z)=(x,\,y+x(z+h_1),\,z+h_1)$ where $h_1 = h(x,\,y+xz),$ 
does not belong, in general, to the group $G.$ 
But composing it further with the automorphism 
$\beta'\in $Aut$\,\C^3,\,\,\,\beta'(x,\,y,\,z)=(x,\,y,\,z-h_1(0,\,y)),$ 
we obtain a new one $\alpha''=\beta'\circ\alpha'$ which does belong to $G.$
It remains to note that $\beta'$ preserves the plane $L_c,$ hence,
$\alpha''(X)=\alpha'(X)=L_c,$ and we are done. \qed 

\smallskip

\noindent {\bf Claim 2.} {\it Let $p(x) =xq(x)$ where ${\rm deg}\,q \ge 1.$
Consider the surface $Y = \{q(x)z = g(x,\,y)\}$ in $\C^3.$ Then 
$X = \Sigma_{I_1,\,\xi}(Y)$ 
is the affine modification of $Y$ along the divisor $D_{\xi}$ with center 
$I_1 = (\xi,\,\eta) \subset \C[Y]$ where  
$\xi = x\,|\,Y$ and $\eta = z\,|\,Y.$ Furthermore, $Y$ is a smooth irreducible
$\Q-$acyclic surface; it is acyclic if
$X$ is. }

%%\smallskip

\noindent {\it Proof.~}  The modification $\Sigma_{I,\,x}(\C^3)$ along the plane
$D_x = \{x=0\}$ with center $I := (x,\,z) \subset \C^{[3]}$ is isomorphic to
$\C^3,$ and the blowup morphism $\si_I\,:\,\C^3 \to \C^3$ is given as 
$(x,\,y,\,z) \longmapsto (x,\,y,\,xz)$ (cf. Examples \ref{ex1.4}, 
\ref{ex2.1}). 
By Proposition \ref{pr2.2}(b), the strict transform $X=Y'$ of $Y$ under 
$\si_I$ coincides with the modification $\Sigma_{I_1,\,\xi}(Y).$ This proves 
the first assertion.

It is easily seen that under our assumptions, 
$X$ being smooth and irreducible implies the same for 
$Y$ (cf. Lemma \ref{lm4.2}).

Due to our assumptions and to Lemma \ref{lm4.5}(b), the irreducible component
$D_p^{1}:=\{x=0\}$ of the divisor $D_p = p^*(0)$ in $\C^2$ meets the divisor 
$D_g$ transversally at one point, say, $P_0 = (0,\,y_0),$ and $D_g$ is reduced 
at $P_0.$ Thus, we have $g(0,\,y) = \gamma (y-y_0).$ Hence, the 
polynomials $x,\,g(0,\,y),\,z \in \C^{[3]}$ generate the maximal ideal of the 
point $P_0':=(0,\,y_0,\,0) \in Y \subset \C^3.$ It follows that 
$I_1 = (\xi,\,\eta)=I(P_0')$ is a maximal ideal of the algebra $\C[Y].$
Thereby, the blowup morphism
$\si\,:\,X \to Y,\,\,\,\si:=\si_{I_1}=\si_I\,|\,X,$ 
consists of the usual blow up at $P_0'$
and deleting the strict transform of the curve 
$l_0 := \{\xi = 0\} \subset Y,\,\,l_0\simeq \C,$ 
passing through $P_0'$ (see Proposition \ref{pr1.1}). 
By Theorem \ref{thm3.1} (see also \cite[Thm. 3.5]{Ka 1}, 
\cite[Thm. 5.1]{Za 2}),
this modification preserves the homology, i.e. $\si_*\,:\,H_*(X;\,\Z)
\longrightarrow H_*(Y;\,\Z)$ is an isomorphism. Thus,
$Y$ is $\Q-$acyclic resp. acyclic iff $X$ is.
\qed 

\smallskip

\noindent {\bf Claim 3: The induction step.} {\it Let $Y$ be as in Claim 2.
Suppose that $Y$ can be transformed into a plane $L_c :=\{y=c\},\,\,c \in \C,$ 
by an automorphism $\alpha \in G.$ Then the same is true for $X.$}

%%\smallskip

\noindent {\it Proof.~}  Notice that the action of the group $G$ on $\C^{[3]}$ 
preserves 
the ideal $I = (x,\,z)$ and fixes the element $x \in I.$ 
Therefore, 
by Corollary \ref{cor2.2}, $\alpha$ can be lifted in a unique way
to an automorphism $\alpha'$ of $\C^3 = \Sigma_{I,\,x}(\C^3)$ such that 
$\alpha\circ \si_I = \si_I \circ \alpha'.$ The automorphism 
$\alpha'$ is of 
the form
$$\alpha'\,:\,(x,\,y,\,z) \longmapsto (x,\,\gamma_1y + 
xg_1(x,\,y,\,xz),\,\gamma_2z + g_2(x,\,y,\,xz))$$ 
(see Example \ref{ex2.2}). It sends the surface $X = Y'$
onto the strict transform $L_c'=L_c$ of the plane $L_c,$ 
i.e. again onto the plane $L_{c}$ (indeed, $\si_I^*(y) = y$). 
Now, composing  $\alpha'$ with the automorphism 
$\beta \in $Aut$\,\C^3,\,\,\,\beta(x,\,y,\,z) = (x,\,y,\,
z-g_2(0,\,\gamma_1^{-1}y,\,0)),$ we get an automorphism 
$\alpha'' :=\beta\circ\alpha'$ which, as it can be easily seen, 
belongs to
the group $G.$ Since $\beta$ preserves the plane $L_{c},$
$\alpha''$ still transforms $X$ into this plane. This proves 
Claim 3 and 
completes the proof of the proposition. \qed 

\smallskip

Now the proof of Theorem \ref{thm4.2} is completed. 
The next example shows that
Theorem \ref{thm4.2} cannot be extended to $\Q-$acyclic surfaces. 

{%\footnotesize
\smallskip

\noindent \nex{{\bf A Q-acyclic modification of the plane.}}
\label{ex4.1} 
Consider the surface $X \subset \C^3$ given by the equation 
$x^2z=x+y^2.$ It is easily seen that $X$
is $\Q-$acyclic but not acyclic. Indeed, in the notation as in 
Proposition 
\ref{pr4.1} we have $E' \simeq \C,$ 
$\,\,X \setminus E' \simeq \C^* \times \C$ and $\si^*(D_p) = 2 E',$ 
which provides that $H_i(X;\,\Z) = 0$ for $i \ge 2$ and 
$H_1(X;\,\Z) \simeq \Z/2\Z$ (cf. the proof of Theorem \ref{thm3.1}).  
Thus, the assumption of Proposition \ref{pr4.1} 
that the curves supp$D_p$ and supp$D_g$ meet transversally
is essential, as well as the condition $E' = \si^*(D)$ in (i) of 
Theorem 
\ref{thm3.1}.  }

\smallskip%\footnotesize
\renewcommand{\arraystretch}{0.8}
\hspace{2mm}\begin{tabular}{lccccl}

Shulim Kaliman    &&&&& Mikhail Zaidenberg\\
Department of Mathematics   &&&&& Institut Fourier de Math{\'e}matiques \\
and Computer Science   &&&&&  UMR 5582   \\
University of Miami  &&&&&  Universit{\'e} Grenoble I\\
Coral Gables   &&&&& BP 74\\
FL  33124 &&&&&  38402 St. Martin d'H{\`e}res--c{\'e}dex \\
U.S.A. &&&&& France\\
e-mail:&&&&& e-mail:\\
kaliman@math.miami.edu &&&&& zaidenbe@ujf-grenoble.fr

\end{tabular}


\begin{thebibliography}{999}
\addcontentsline{toc}{a}{References}

%%%% If there is no [] after (all) \bibitems, then LaTex numbers them. %%%%

%%\footnotesize

\bibitem[AM]{AM} S.S. Abhyankar, T.T. Moh, {\em Embedding of the line in the
plane}, J. Reine Angew. Math. {\bf 276} (1975), 148--166.
\bibitem[Akh]{Akh} D.N. Akhiezer, {\em Lie group actions in complex analysis}, 
Prospects in Mathem., Vieweg, 1995.
\bibitem[As]{As} T. Asanuma, {\em Non-linearizable
algebraic $k^*$-actions on affine spaces}, preprint, 1996, 23p.
 \bibitem[Bre]{Bre} G.E.Bredon,
{\em Introduction to compact transformation groups}, 
Academic Press, N.Y., 1972.
 \bibitem[Bri]{Bri} E. Brieskorn,
{\em Beispiele zur Differentialtopologie von Singularit\"aten},
Invent. Math. {\bf 2} (1966), 1-14.
 \bibitem[Da]{Da} E.D. Davis,  {\em Ideals of the principal class,
$R$-sequences and a certain monoidal transformation}, Pacific J. Math. 
{\bf 20} (1967), 197--205.
\bibitem[De]{De} H. Derksen, {\em
Constructive Invariant Theory and the Linearization Problem},
Ph.D. thesis, Basel, 1997.
 \bibitem[tDP]{tDP} T. tom Dieck, T. Petrie, {\em Contractible affine
surfaces of Kodaira dimension one}, Japan J. Math. {\bf 16} (1990), 147--169.
\bibitem[Di 1]{Di 1} A. Dimca, {\em Hypersurfaces in ${\C}^{2n}$ diffeomorphic
to ${\bf R}^{4n - 2} \,(n \geq 2)$}, Max-Planck Institute, preprint, 1991.
\bibitem[Di 2]{Di 2} A. Dimca, {\em Singularities and Topology of
Hypersurfaces}, Universitext, Springer, 1992.
\bibitem[Do]{Do} A. Dold, {\em Lectures on algebraic topology}, Springer,
Berlin e.a., 1974.
\bibitem[Du]{Du} A. Durfee, {\em Algebraic varieties which are a disjoint 
union of subvarieties,} In: Geometry and Topology. Proc. Conf. Athens/Ga, 
1985. Lect. Notes Pure Appl. Math. {\bf 105} (1987), 99-102.
 \bibitem[Ei]{Ei} D. Eisenbud, {\em Commutative algebra with a view towards 
 algebraic geometry}, Graduate Texts in Math., Springer, N.Y. e.a., 1994
 \bibitem[FlZa]{FlZa} H. Flenner, M. Zaidenberg,
{\sl Q-acyclic surfaces and their deformations}, Contempor. Mathem. {\bf 162}, 
Providence, RI, 1994, 143--208.
\bibitem[FoFu]{FoFu} A.T. Fomenko, B.D. Fuks, {\em A course in homotopic topology,} 
Textbook, Nauka, Moscow, 1989 (Russian).
\bibitem[Fuj]{Fuj} T. Fujita, {\em On the topology of non-complete algebraic
surfaces}, J. Fac. Sci. Univ. Tokyo, Sect.IA, {\bf 29} (1982), 503--566.
\bibitem[Ful]{Ful} W. Fulton, {\em Intersection Theory}, Ergebnisse der 
Mathematik, B.2, Springer, Berlin e.a., 1984.
 \bibitem[Ha]{Ha} R. Hartshorn, {\em Algebraic geometry}, Springer,
NY e.a., 1977.
 \bibitem[Hiro]{Hiro} H. Hironaka, {\em Resolution of singularities}, 
Ann. Math. {\bf 79} (1964), 109--326.
% \bibitem[HNK]{HNK} F. Hirzebruch, W.D. Neumann, and S.S. Koh, 
%{\em Differentiable manifolds and quadratic forms}, 
%Lect. Notes Pure Appl. Math. 
%{\bf 4}, M. Dekker Inc., New York, 1971.
 \bibitem [Je]{Je} Z. Jelonek, {\em The extension of regular and rational
embeddings}, Math. Ann. {\bf 113} (1987) 113--120.
 \bibitem [Ju]{Ju} H. W. E. Jung, {\em {\"U}ber ganze birationale
Transformationen der Ebene}, J. reine und angew. Math., {\bf 184} (1942),
161--174.
 \bibitem[Ka 1]{Ka 1} S. Kaliman, {\em Exotic analytic structures and
Eisenman intrinsic measures}, Israel Math. J. {\bf 88} (1994), 411--423.
 \bibitem[Ka 2]{Ka 2} S. Kaliman, {\em Extensions of isomorphisms between
affine algebraic subvarieties of $k^n$ to automorphisms of $k^n$}, Proc.
Amer. Math. Soc. {\bf 113} (1991), 325--334.
 \bibitem[Ka 3]{Ka 3} S. Kaliman, {\em Smooth contractible hypersurfaces in
${\C}^{n}$ and exotic algebraic structures on ${\C}^{3}$}, Math. Zeitschrift
{\bf 214} (1993), 499--510.
 \bibitem[KaML 1]{KaML 1} S. Kaliman, L. Makar-Limanov, {\em On
Russell--Koras contractible threefolds}, 
J. of Algebraic Geom. {\bf 6} (1997), 247-268.
\bibitem[KaML 2]{KaML 2} S. Kaliman, L. Makar-Limanov, {\em Affine algebraic
manifolds without dominant morphisms from Euclidean spaces}, 
Rocky Mount. J. Math. {\bf 27:2} (1997), 601 --609.
 \bibitem[ML 1]{ML 1} L. Makar-Limanov,
{\em On the hypersurface $x + x^2y + z^2 + t^3 = 0$
in ${\C}^{4}$ or a ${\C}^3$-like threefold which is not ${\C}^3$}, Israel 
J. Math. {\bf 96} (1996), 419--429.
 \bibitem[ML 2]{ML 2} L. Makar-Limanov,
{\em On groups of automorphisms of a class of surfaces},
Israel J. Math. {\bf 69} (1990), 250-256.
\bibitem[Mik]{Mik} A. Micali, {\em Sur les alg{\`e}bres universalles}, 
Annales Inst. Fourier, Grenoble {\bf 14} (1964), 33--88.
\bibitem[Mil 1]{Mil 1} J. Milnor, {\em Morse Theory}, Princeton Univ. Press,
Princeton, NJ, 1963. 
\bibitem[Mil 2]{Mil 2} J. Milnor, {\em On the 3-dimensional Brieskorn manifolds
$M(p,\,q,\,r)$}, in: Knots, groups, and 3-manifolds, L. P. Neuwirth, ed. 
Annals of Math. Stud., Princeton Univ. Press, Princeton, NJ, 1975, 175--225.
 \bibitem[MilSta]{MilSta} J. Milnor, J. Stasheff, {\em Characteristic classes},
Annals of Mathem. Studies {\bf 76}, Princeton Univ. Press and Univ. of Tokyo
Press, Princeton, NJ, 1974.
 \bibitem[Miy]{Miy} M. Miyanishi, {\em Algebraic characterization of the
affine 3--space}, Proc. Algebraic Geom. Seminar, Singapore, World Scientific,
1987, 53--67.
  \bibitem[Ram]{Ram} C.P. Ramanujam, {\sl A topological characterization of the
affine plane as an algebraic variety}, Ann. Math., {\bf 94} (1971), 69-88.
 \bibitem[Re]{Re} R. Rentschler, {\em Op{\'e}rations du groupe additif sur le
plane affine}, C.R. Acad. Sci. Paris, {\bf 267} (1968), 384--387.
 \bibitem[RoRu]{RoRu} J.-P. Rosay, W. Rudin, {\em Holomorphic maps from 
 $\C^n$ to $\C^n,$} Trans. Amer. Math. Soc. {\bf 310} (1988), 47--86.
 \bibitem[Ru 1]{Ru 1} P. Russell, {\em On a class of ${\C}^3$-like threefolds},
Preliminary Report, Berlin, 1992.
 \bibitem[Ru 2]{Ru 2} P. Russell, 
{\em Simple birational extensions of two dimensional affine rational domains},
Compos. Math. {\bf 33} (1976), 197-208.
 \bibitem[RuSat]{RuSat} P. Russell, A. Sathaye, {\em On finding and 
 cancelling variables in $k[X,\,Y,\,Z]$}, J. Algebra {\bf 57} (1979), 151--166.
  \bibitem[Sat]{Sat} A. Sathaye, {\em On linear planes}, Proc. Amer. Math. Soc.
{\bf 56} (1976), 1--7.
 \bibitem[Se]{Se} J.-P. Serre, {\em Trees}, Springer Verlag,
Berlin-Heidelberg-New York, 1980.
 \bibitem[Sr]{Sr} V. Srinivas, {\em On the embedding dimension of an affine
variety}, Math. Ann.  {\bf 289} (1991), 125-132.
 \bibitem[Suz]{Suz} M. Suzuki, {\em Propi{\'e}t{\'e}s topologiques des
polyn{\^o}mes de deux variables complexes, et automorphismes alg{\'e}brigue de
l'espace
${\C}^{2}$}, J. Math. Soc. Japan, {\bf 26} (1974), 241-257.
 \bibitem[vdK]{vdK} W. van der Kulk,
{\em On polynomial rings in two variables}, Nieuw Arch. Wisk. (3) {\bf 1}
(1953), 33--41.
 \bibitem[Va]{Va} W.V. Vasconselos, {\em Arithmetic of blowup algebras}, 
London Math. Soc. Lect. Notes Ser. {\bf 195}, Cambridge Univ. Press, 1994.
\bibitem[Wr 1]{Wr 1} D. Wright, {\em Abelian subgroups of 
${\rm Aut}_k(k[X,Y])$
and applications to actions on the affine plane}, Ill. J. Math. {\bf 23}
(1979), 579--634.
\bibitem[Wr 2]{Wr 2} D. Wright, {\em Cancellation of variables of the form 
$bT^n - a$}, J. Algebra {\bf 52} (1978), 94--100.
\bibitem[Za 1]{Za 1} M. Zaidenberg, {\sl Rational actions of the group 
$\bf C^*$ on ${\bf C}^{2}$, their quasi-invariants, and algebraic curves in 
${\bf C}^{2}$ with Euler characteristic 1}, Soviet Math. Dokl. {\bf 31:1} 
(1985), 57--60.
\bibitem[Za 2]{Za 2} M. Zaidenberg, 
{\em Lectures on exotic algebraic structures on affine spaces}, 
Schriftenreihe des Graduiertenkollegs Geometrie und Mathematische 
Physik, Heft 24, Institut F{\"u}r Mathematik, Ruhr-Universit{\"a}t 
Bochum, Juli 1997, 56p.
\bibitem[Zar]{Zar} O. Zariski, {\em Foundations of a general theory of 
birational correspondences}, Trans. Amer. Math. Soc. {\bf 53} (1943), 497--542

\end{thebibliography}
\end{document}